\documentclass[a4paper,10pt]{article}
\usepackage{amsmath,pictexwd}
\usepackage{enumerate}
\usepackage{amssymb,latexsym,color}
\usepackage{geometry}
\usepackage{xypic}
\xyoption{2cell}\UseTwocells
\UseHalfTwocells
\author{ Timothy Porter and Vladimir Turaev}

\title{Formal Homotopy Quantum Field Theories, I: \\ Formal Maps and Crossed $\mathcal{C}$-algebras.}
 
\newtheorem{theorem}{Theorem}
\newtheorem{proposition}[theorem]{Proposition} 
\newtheorem{corollary}[theorem]{Corollary}
\newtheorem{lemma}[theorem]{Lemma}
  
\begin{document}

\maketitle
\begin{abstract}
Homotopy Quantum Field Theories (HQFTs) were introduced by the second author to extend the ideas and methods of Topological Quantum Field Theories to closed $d$-manifolds endowed with extra structure in the form of homotopy classes of maps into a given `target' space, $B$. For $d = 1$, classifications of HQFTs in terms of algebraic structures are known when $B$ is a $K(G,1)$ and also when it is simply connected.  Here we study general HQFTs with $d = 1$ and target a general 2-type, giving a common generalisation of the classifying algebraic structures for the two cases previously known.  The algebraic models for 2-types that we use are crossed modules, $\mathcal{C}$, and we introduce a notion of formal $\mathcal{C}$-map, which extends the usual lattice-type constructions to this setting. This leads to a classification of `formal' 2-dimensional HQFTs with target $\mathcal{C}$, in terms of crossed $\mathcal{C}$-algebras.

{\noindent\bf A. M. S. Classification:} Primary:  18G50. Secondary: 55P99, 57R56, 81T45.  
\end{abstract}

\tableofcontents

\section{Introduction}
Homotopy Quantum Field Theories were introduced in \cite{turaev:hqft1} as an extension of the notion of a Topological Quantum Field Theory to $d$-manifolds and $(d+1)$-dimensional cobordisms endowed with extra structure in the form of `characteristic maps' to a fixed pointed `background' or `target' space, $B$. It is known that, for a given $d$ and $B$, these $(d+1)$-dimensional  HQFTs only use the $(d+1)$-type of $B$, 
that is, the structure of the homotopy type of $B$ up to and including the $d+1^{st}$ homotopy group, (see \cite{rodrigues}).

Restricting to the case $d=1$,  as we will in this paper, the corresponding HQFTs are variously refered to as being 2-dimensional or 1+1 dimensional. Of course, any closed connected 1-manifold\footnote{In general all our manifolds will be closed and orientable, so we will not always mention this explicitly.} over $B$ (1-$B$-manifold) is an oriented circle with a map to $B$, so will determine a homotopy class of maps $g : S^1\to B$ and hence an element of $\pi_1(B)$.  This allows a combinatorial model of the basic objects to be given (cf., the \emph{$\pi$-systems} of \cite{turaev:hqft1} \S 7.2), at least when $B$ is a $K(G,1)$.   A somewhat similar approach was used by Brightwell and Turner, \cite{B&T}, when $B$ is a simply connected 2-type, hence specified, up to homotopy, by its second homotopy group $\pi_2(B)$, which we will often write just as $A$.

 Various equivalent algebraic models for general 2-types are known, for instance, crossed modules, 2-groups, cat$^1$-groups, ... .  We have chosen to work with crossed modules as they are probably the simplest to use whilst being very near to the group theoretic methods that are well known for other cases.  (Crossed modules are strict 2-groups, for the reader used to such things, but their theory has been around a lot longer as it was initiated by Reidemeister, Peiffer and Whitehead in the 1940s and early 1950s.)  It is therefore natural to seek a common extension of the combinatorial systems used for the special cases in terms of such models.  Those methods allowed algebraic classifying objects, \emph{crossed $\mathcal{C}$-algebras}, to be identified that corresponded well to 2-dimensional HQFTs, namely,
 \begin{itemize}	\item  the crossed $G$-algebras, when $B\cong K(G,1)$, (cf., \cite{turaev:hqft1}, where the group is called $\pi$);
	\item  the $A$-Frobenius algebras when $B \cong K(A,2)$, (cf., \cite{B&T}),
\end{itemize}
for those special cases.

One slight complication arises, however, when we pass to the general case.  If $B$ is a 2-type, there will be a crossed module, $\mathcal{C}$, say, whose classifying space\footnote{cf., Brown-Higgins, \cite{B&H1991} and Porter, \cite{tqft1} for information on classifying spaces of crossed modules and crossed complexes.}, $B\mathcal{C}$, has the same homotopy 2-type as $B$ and therefore can provide an algebraic model for that, but there will be other crossed modules, not \emph{isomorphic}   to $\mathcal{C}$, for which this is also true.  It is \emph{weak equivalence} classes of crossed modules that correspond to 2-types not \emph{isomorphism classes}.  We therefore tackle one part only of the classification problem here.  Given a crossed module $\mathcal{C}$ providing an algebraic model for a 2-type $B$, we introduce a combinatorial (lattice gauge-like) model for the 1-dimensional $B$-manifolds and the corresponding cobordisms, then we give an analogue of HQFTs with target $B$.  We call these combinatorial gadgets \emph{formal $\mathcal{C}$-maps} and the resulting analogues of HQFTs, \emph{formal (2-dimensional) HQFTs}, although we will sometimes omit the `2-dimensional' as we will not be considering other cases here.  We classify formal HQFTs over a given $\mathcal{C}$ in terms of crossed $\mathcal{C}$-algebras, our promised common generalisation of crossed $\pi$-algebras and $A$-Frobenius algebras.  More precisely we will prove:

\

\textbf{Main Theorem}

\emph{There is a canonical bijection between isomorphism classes of formal 2-dimensional HQFTs based on a crossed module $\mathcal{C}$ and isomorphism classes of crossed $\mathcal{C}$-algebras.}

\

We start the paper  with a discussion of how to extend the notion of a $G$-colouring of a triangulation of a manifold, for $G$ a group,  to a colouring with values in a crossed module, $\mathcal{C}$. We introduce and briefly discuss crossed modules, before defining simplicial formal $\mathcal{C}$-maps as that extension of $G$-colourings, and also equivalence of formal $\mathcal{C}$-maps.  These ideas are related to ideas already explored to some extent in TQFTs.  There is a fairly obvious extension of these notions from simplicial complexes to CW-complexes. These are introduced after a section discussing the methods of simplifying formal maps within an equivalence class.  A detailed look at 2-dimensional formal $\mathcal{C}$-maps follows.

In section 3, we formally define formal HQFTs and explore some of the elementary consequences of that definition.  Crossed $\mathcal{C}$-algebras are introduced in section 4.  These generalise Frobenius algebras, as mentioned above.  Specific examples are postponed until section 6.  Section 5 contains the proof of the main theorem. This uses a combination of ideas from higher category theory, with material from TQFTs and geometric topology. Much of this follows the same track as for the simpler case of crossed $\pi$-algebras and HQFTs with background a $K(\pi,1)$ considered by the second author in \cite{turaev:hqft1}, but some parts involve the top group, $C$, of the crossed module as well and hence a new structural element.

We then embark, in section 6, on the detailed study of these crossed $\mathcal{C}$-algebras, introducing extensions of various constructions, pullback and pushforward, already known for crossed $\pi$-algebras.  Later these will be applied to give comparison results for change of the base $\mathcal{C}$, which will be essential when examining the way in which the algebras reflect the weak equivalence class of $\mathcal{C}$.

\

\textbf{Acknowledgements}

This work was supported by a grant, GR/S17635/01, from the EPSRC, for a visit by the second author to Bangor and Gregynog Hall.  The occasion was a regional meeting of the London Mathematical Society, locally organised by David Evans and Edwin Beggs. We would like to thank all of these, both organisations and individuals, and the staff of Gregynog Hall, for providing the opportunity for this work to be undertaken in very pleasant surroundings.  It has also benefited from discussions during a visit of the first author to the University of Ottawa, Summer 2007.

We would also like to acknowledge the help given by Ronnie Brown, who participated in many of the discussions both at Bangor and at Gregynog. His wealth of ideas, perspective and knowledge on crossed modules, crossed complexes and all the general `crossed menagerie' was invaluable.

\section{Crossed modules  and formal maps}
\subsection{A simplicial (lattice) approach}
In the construction of models for topological and homotopical quantum field theories, one often uses a (finite) group $G$, and a triangulation of the manifolds, $\Sigma$, etc., involved, and one assigns labels from $G$ to each (oriented)  edge of each (oriented) triangle, for example,\\

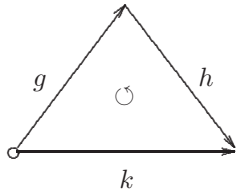
\begin{figure}[h]\begin{center}
\font\thinlinefont=cmr5
\begingroup\makeatletter\ifx\SetFigFont\undefined
\def\x#1#2#3#4#5#6#7\relax{\def\x{#1#2#3#4#5#6}}%
\expandafter\x\fmtname xxxxxx\relax \def\y{splain}%
\ifx\x\y   
\gdef\SetFigFont#1#2#3{%
  \ifnum #1<17\tiny\else \ifnum #1<20\small\else
  \ifnum #1<24\normalsize\else \ifnum #1<29\large\else
  \ifnum #1<34\Large\else \ifnum #1<41\LARGE\else
     \huge\fi\fi\fi\fi\fi\fi
  \csname #3\endcsname}%
\else
\gdef\SetFigFont#1#2#3{\begingroup
  \count@#1\relax \ifnum 25<\count@\count@25\fi
  \def\x{\endgroup\@setsize\SetFigFont{#2pt}}%
  \expandafter\x
    \csname \romannumeral\the\count@ pt\expandafter\endcsname
    \csname @\romannumeral\the\count@ pt\endcsname
  \csname #3\endcsname}%
\fi
\fi\endgroup
\mbox{\beginpicture
\setcoordinatesystem units <0.750000cm,0.750000cm>
\unitlength=0.750000cm
\linethickness=1pt
\setplotsymbol ({\makebox(0,0)[l]{\tencirc\symbol{'160}}})
\setshadesymbol ({\thinlinefont .})
\setlinear
%
%
\linethickness= 0.500pt
\setplotsymbol ({\thinlinefont .})
\plot  3.6 21.220  5.5 23.732 /
%
%
\plot  5.394 23.531  5.5 23.732  5.309 23.594 /
%
%
%
\linethickness= 0.500pt
\setplotsymbol ({\thinlinefont .})
\plot  5.50 23.732  7.436 21.220 /
%
%
\plot  7.265 21.355  7.436 21.220  7.349 21.420 /
%
%
%
\linethickness= 0.500pt
\setplotsymbol ({\thinlinefont .})
\putrule from  3.573 21.139 to  7.436 21.139
%
%
\plot  7.224 21.086  7.436 21.139  7.224 21.192 /
%
%
%
\put{$g$} [lB] at  3.863 22.3
%
%
\put{$h$} [lB] at  6.8 22.3
%
%
\put{$k$} [lB] at  5.4 20.5
%
%
\put{ $\circlearrowleft$} [lB] at  5.15  22
%
%
\linethickness= 5.500pt
\setplotsymbol ({\thinlinefont .})
\ellipticalarc axes ratio  0.093:0.093  360 degrees 
	from  3.613 21.114 center at  3.520 21.114
\linethickness=0pt
\putrectangle corners at  3.410 23.757 and  7.461 20.701
\endpicture}
\label{figure 1}\end{center}
\end{figure}
\noindent with the boundary/cocycle condition\footnote{Here the orientation is given as anticlockwise, which seems unnatural given the ordering, but this is necessary as we are using the `path order convention' on composition of labels on edges. The other convention also leads to some inelegance at times. We use both!}  that $kh^{-1}g^{-1} = 1$, so $k = gh$.

The geometric intuition behind this is that `integrating' the labels around the triangle  yields the identity.  This intuition corresponds to situations in which a $G$-bundle on $\Sigma$ is specified by charts and the elements $g$, $h$, $k$, etc. are transition automorphisms of the fibre.  The methods then use manipulations of the pictures as the triangulation is changed by subdivision, etc.

Another closely related view of this is to consider continuous functions $f : \Sigma \to BG$ to the classifying space of $G$.  If we triangulate $\Sigma$, we can assume that $f$ is a cellular map using a suitable cellular model of $BG$ and at the cost of replacing $f$ by a homotopic map and perhaps subdividing the triangulation.  From this perspective the previous model is a combinatorial description of such a continuous `characteristic' map, $f$.  The edges of the triangulation pick up group elements since the end points of each edge get mapped to the base point of $BG$, and $\pi_1BG \cong G$, whilst the faces give a realisation of the cocycle condition. Similarly we could use a labelled decomposition of the objects as  CW-complexes, cf. \cite{turaev:hqft1,lundell-weingram} and again the edges would pick up group elements, whilst the two cells give a cocycle condition.

The 1+1 homotopy quantum field theories, in general, work with objects,  1-dimensional $B$-manifolds, that are closed oriented 1-manifolds, $\Sigma$, with a `characteristic map' from $\Sigma$ to a general fixed background or target space, $B$, and it is known, cf. \cite{rodrigues}, that only the 2-type of $B$ contributes to the theory.  If $B$ is a $K(G,1)$, as above, then the existing theory and diagrams work well and yield a classification of the corresponding HQFTs, \cite{rodrigues,turaev:hqft1}.  If $B$ is a $K(A,2)$, so is simply connected, then Brightwell and Turner, \cite{B&T}, have related classification results, but what happens for a general 2-type, $B$?

Let $B$ be a CW-complex model for a 2-type (so $\pi_kB$ is trivial for $k>2$).  Assume it is reduced, so has a single vertex, then, denoting by $B_1$, the 1-skeleton of $B$,  the crossed module, $(\pi_2(B, B_1),\pi_1(B_1),\partial)$, will represent the 2-type of $B$ .  For any $B$-manifold, the characteristic map, $g : \Sigma \to B$, or for a $B$-cobordism,  the map, $F : M \to B$, can be replaced, up to homotopy, by a cellular map, so, in general, we can think of a combinatorial model for the $B$-manifolds and $B$-cobordisms  in terms of combining labelled triangles
\begin{figure}\begin{center}
\font\thinlinefont=cmr5
\begingroup\makeatletter\ifx\SetFigFont\undefined
\def\x#1#2#3#4#5#6#7\relax{\def\x{#1#2#3#4#5#6}}%
\expandafter\x\fmtname xxxxxx\relax \def\y{splain}%
\ifx\x\y   
\gdef\SetFigFont#1#2#3{%
  \ifnum #1<17\tiny\else \ifnum #1<20\small\else
  \ifnum #1<24\normalsize\else \ifnum #1<29\large\else
  \ifnum #1<34\Large\else \ifnum #1<41\LARGE\else
     \huge\fi\fi\fi\fi\fi\fi
  \csname #3\endcsname}%
\else
\gdef\SetFigFont#1#2#3{\begingroup
  \count@#1\relax \ifnum 25<\count@\count@25\fi
  \def\x{\endgroup\@setsize\SetFigFont{#2pt}}%
  \expandafter\x
    \csname \romannumeral\the\count@ pt\expandafter\endcsname
    \csname @\romannumeral\the\count@ pt\endcsname
  \csname #3\endcsname}%
\fi
\fi\endgroup
\mbox{\beginpicture
\setcoordinatesystem units <0.650000cm,0.650000cm>
\unitlength=0.650000cm
\linethickness=1pt
\setplotsymbol ({\makebox(0,0)[l]{\tencirc\symbol{'160}}})
\setshadesymbol ({\thinlinefont .})
\setlinear
%
%
\linethickness= 0.500pt
\setplotsymbol ({\thinlinefont .})
\plot  3.6 21.220  5.5 23.732 /
%
%
\plot  5.394 23.531  5.5 23.732  5.309 23.594 /
%
%
%
\linethickness= 0.500pt
\setplotsymbol ({\thinlinefont .})
\plot  5.50 23.732  7.436 21.220 /
%
%
\plot  7.265 21.355  7.436 21.220  7.349 21.420 /
%
%
%
\linethickness= 0.500pt
\setplotsymbol ({\thinlinefont .})
\putrule from  3.573 21.139 to  7.436 21.139
%
%
\plot  7.224 21.086  7.436 21.139  7.224 21.192 /
%
%
%
\put{$g$} [lB] at  3.863 22.3
%
%
\put{$h$} [lB] at  6.8 22.3
%
%
\put{$k$} [lB] at  5.4 20.5
%
%
\put{ $c  \circlearrowleft$} [lB] at  4.8  22
%
%
\linethickness= 5.500pt
\setplotsymbol ({\thinlinefont .})
\ellipticalarc axes ratio  0.093:0.093  360 degrees 
	from  3.613 21.114 center at  3.520 21.114
\linethickness=0pt
\putrectangle corners at  3.410 23.757 and  7.461 20.701
\endpicture}
\end{center}\end{figure}
with $g,h,k \in \pi_1(B_1)$ and $c\in\pi_2(B,B_1)$, and where the cocycle condition is replaced by a boundary condition of form
$$\partial c = kh^{-1}g^{-1} .$$
Usually $\pi_1(B_1)$ will be free and it will be useful to replace this particular crossed module by a general one.

\medskip

\textbf{Definition}

A \emph{crossed module}, $\mathcal{C} = (C,P, \partial)$, consists of groups $C$, $P$, a (left) action of
$P$ on $C$ (written $(p,c)\rightarrow {}^pc)$ and a homomorphism
$$\partial : C \rightarrow P$$
 such that

CM1 \quad $\partial({}^pc) = p\cdot\partial c\cdot p^{-1}$ \quad\quad for all $p 
\in P$, $c \in C$,\quad ($\partial$ is $P$-equivariant)
\\
and

 CM2 \quad ~ ${}^{\partial c} c^\prime = c\cdot c^\prime\cdot c^{-1}$ \quad\quad\quad for all
 $c, c^\prime \in C$, \quad (the Peiffer identitiy axiom).

\medskip

There are several well known examples of crossed modules.  We mention three:
\begin{enumerate}[(i)]
\item If $N$ is a normal subgroup of a group $P$, then $P$ acts by conjugation on $N$, ${}^pn = pnp^{-1}$, and the inclusion $\iota : N \to P$ is a crossed module. (Conversely, if $(C,P, \partial)$ is a crossed module, $\partial C$ is a normal subgroup of $P$.)
\item  If $M$ is a left $P$-module and we define $0 : M \to P$ to be the trivial homomorphism, $0(m) = 1_P$, for all $m \in M$, then  $(M,P, 0)$ is crossed module. (Conversely if $(C,P,\partial)$ is a crossed module, then $\ker \partial$ is a $P$-module, in fact a $P/\partial C$-module, as the image $\partial C$ acts trivially on the kernel.)
\item If $G$ is any group, $\alpha : G \to Aut(G)$, the canonical map sending $g \in G$ to the inner automorphism determined by $g$, is a crossed module for the standard action of $Aut(G)$ on $G$.
(This third example is of a generic type and later, (Lemma \ref{AutL}), we will see that for an algebra, $L$, $(U(L),Aut(L),\delta)$ is a crossed module, where $U(L)$ is the group of units of $L$ and $\delta$ maps a unit to the automorphism given by conjugation by it.)

\end{enumerate}


\textbf{Remark}

We recall that to any crossed module, $\mathcal{C}$, there is an associated (strict) 2-group\footnote{A strict 2-group is a 2-category with one object, for which both 1-cells and 2-cells are invertible} with a single object, having $P$ as its group of automorphism (1-cells) and the semidirect product, $C\rtimes P$ as its group of 2-cells. The Peiffer identity corresponds to the Interchange Law in that 2-category,  (see later for a bit more on this).

\

From now on, we fix a crossed module $\mathcal{C} = (C,P,\partial)$ as given.  Our formal $\mathcal{C}$-maps will initially be introduced via $\mathcal{C}$-labelled triangles as above, but will then be replaced by a cellular version as soon as the basic results are established confirming some basic intuitions.  The labelled triangles, tetrahedra, etc., will all need a base point as a `start vertex'.  The need for this can be seen in an elementary way as follows:\\
If we have the situation below, 
\begin{figure}[h]\begin{center}
\font\thinlinefont=cmr5
\begingroup\makeatletter\ifx\SetFigFont\undefined
\def\x#1#2#3#4#5#6#7\relax{\def\x{#1#2#3#4#5#6}}%
\expandafter\x\fmtname xxxxxx\relax \def\y{splain}%
\ifx\x\y   
\gdef\SetFigFont#1#2#3{%
  \ifnum #1<17\tiny\else \ifnum #1<20\small\else
  \ifnum #1<24\normalsize\else \ifnum #1<29\large\else
  \ifnum #1<34\Large\else \ifnum #1<41\LARGE\else
     \huge\fi\fi\fi\fi\fi\fi
  \csname #3\endcsname}%
\else
\gdef\SetFigFont#1#2#3{\begingroup
  \count@#1\relax \ifnum 25<\count@\count@25\fi
  \def\x{\endgroup\@setsize\SetFigFont{#2pt}}%
  \expandafter\x
    \csname \romannumeral\the\count@ pt\expandafter\endcsname
    \csname @\romannumeral\the\count@ pt\endcsname
  \csname #3\endcsname}%
\fi
\fi\endgroup
\mbox{\beginpicture
\setcoordinatesystem units <0.750000cm,0.750000cm>
\unitlength=0.750000cm
\linethickness=1pt
\setplotsymbol ({\makebox(0,0)[l]{\tencirc\symbol{'160}}})
\setshadesymbol ({\thinlinefont .})
\setlinear
%
%
\linethickness= 0.500pt
\setplotsymbol ({\thinlinefont .})
\plot  3.505 21.141  7.436 21.139 /
%
%
\plot  7.224 21.086  7.436 21.139  7.224 21.192 /
%
%
%
\linethickness= 0.500pt
\setplotsymbol ({\thinlinefont .})
\plot  3.505 21.129  5.478 23.732 /
%
%
\plot  5.392 23.531  5.478 23.732  5.308 23.595 /
%
%
%
\put{$g$} [lB] at  3.863 22.6
%
%
\put{$h$} [lB] at  6.615 22.6
%
%
\put{$k$} [lB] at  5.4 20.5
%
%
\put{ $c ~ \circlearrowleft$} [lB] at  4.7  22

%
\linethickness= 0.500pt
\setplotsymbol ({\thinlinefont .})
\plot  5.503 23.732  7.436 21.220 /
%
%
\plot  7.265 21.355  7.436 21.220  7.349 21.420 /
%
%
%
\put{0} [lB] at  3.1    21.
%
%

%
\put{1} [lB] at  5.35 24.
%
%
\put{2} [lB] at  7.6 21.
\linethickness=0pt
\putrectangle corners at  3.258 24.115 and  8.105 20.752
\endpicture}
\end{center}\end{figure}
we get the boundary condition $\partial c = kh^{-1}g^{-1}$, which was read off starting at vertex 0: first $k$, back along $h$ giving $h^{-1}$, then the same for $g$ giving $g^{-1}$.  The element $c$ is assigned to this 2-simplex with this ordering / orientation, but if we tried to read off the boundary starting at vertex 1, we would get $g^{-1}kh^{-1}$, which is not $\partial c$, but is $\partial ({}^{g^{-1}}c)$. We thus have that the $P$-action on $C$ is precisely encoding the change of starting vertex.

\medskip										

\textbf{Remark:}

Our simplices will have a marked vertex to enable the boundary condition, and later on a cocycle condition, to be read off unambiguously.  We could equally well work with a pair of marked vertices corresponding to `start' and `finish' or `source' and `target'.  For triangles this would give, for instance, the above with start at 0 and finish at 2, and would give a boundary condition read off as $k = \partial c \cdot gh$. This can lead to a 2-categorical formulation of formal $\mathcal{C}$-maps, which is connected with the way in which a crossed module $\mathcal{C}$ is equivalent to a strict 2-group.  This latter approach was used to develop some of the theory outlined below, and may be useful for future development.  For the cellular version, this leads to globular diagrams 
\begin{figure}[h]\begin{center}
\font\thinlinefont=cmr5
\begingroup\makeatletter\ifx\SetFigFont\undefined
\def\x#1#2#3#4#5#6#7\relax{\def\x{#1#2#3#4#5#6}}%
\expandafter\x\fmtname xxxxxx\relax \def\y{splain}%
\ifx\x\y   
\gdef\SetFigFont#1#2#3{%
  \ifnum #1<17\tiny\else \ifnum #1<20\small\else
  \ifnum #1<24\normalsize\else \ifnum #1<29\large\else
  \ifnum #1<34\Large\else \ifnum #1<41\LARGE\else
     \huge\fi\fi\fi\fi\fi\fi
  \csname #3\endcsname}%
\else
\gdef\SetFigFont#1#2#3{\begingroup
  \count@#1\relax \ifnum 25<\count@\count@25\fi
  \def\x{\endgroup\@setsize\SetFigFont{#2pt}}%
  \expandafter\x
    \csname \romannumeral\the\count@ pt\expandafter\endcsname
    \csname @\romannumeral\the\count@ pt\endcsname
  \csname #3\endcsname}%
\fi
\fi\endgroup
\mbox{\beginpicture
\setcoordinatesystem units <0.75000cm,0.75000cm>
\unitlength=0.75000cm
\linethickness=1pt
\setplotsymbol ({\makebox(0,0)[l]{\tencirc\symbol{'160}}})
\setshadesymbol ({\thinlinefont .})
\setlinear
%
%
\linethickness= 0.500pt
\setplotsymbol ({\thinlinefont .})
\ellipticalarc axes ratio  0.068:0.068  360 degrees 
	from  3.535 20.637 center at  3.467 20.637
%
%
\linethickness= 0.500pt
\setplotsymbol ({\thinlinefont .})
\plot  3.467 20.637  3.469 20.640 /
\plot  3.469 20.640  3.473 20.644 /
\plot  3.473 20.644  3.480 20.652 /
\plot  3.480 20.652  3.490 20.665 /
\plot  3.490 20.665  3.505 20.682 /
\plot  3.505 20.682  3.526 20.705 /
\plot  3.526 20.705  3.552 20.733 /
\plot  3.552 20.733  3.581 20.767 /
\plot  3.581 20.767  3.617 20.805 /
\plot  3.617 20.805  3.655 20.845 /
\plot  3.655 20.845  3.700 20.889 /
\plot  3.700 20.889  3.747 20.936 /
\plot  3.747 20.936  3.797 20.983 /
\plot  3.797 20.983  3.850 21.031 /
\plot  3.850 21.031  3.907 21.082 /
\plot  3.907 21.082  3.967 21.131 /
\plot  3.967 21.131  4.030 21.181 /
\plot  4.030 21.181  4.098 21.230 /
\plot  4.098 21.230  4.170 21.281 /
\plot  4.170 21.281  4.248 21.330 /
\plot  4.248 21.330  4.331 21.378 /
\plot  4.331 21.378  4.420 21.425 /
\plot  4.420 21.425  4.517 21.471 /
\plot  4.517 21.471  4.621 21.518 /
\plot  4.621 21.518  4.733 21.562 /
\plot  4.733 21.562  4.851 21.605 /
\plot  4.851 21.605  4.974 21.643 /
\plot  4.974 21.643  5.093 21.675 /
\plot  5.093 21.675  5.207 21.700 /
\plot  5.207 21.700  5.317 21.723 /
\plot  5.317 21.723  5.419 21.742 /
\plot  5.419 21.742  5.512 21.759 /
\plot  5.512 21.759  5.594 21.772 /
\plot  5.594 21.772  5.668 21.781 /
\plot  5.668 21.781  5.732 21.789 /
\plot  5.732 21.789  5.785 21.795 /
\plot  5.785 21.795  5.831 21.800 /
\plot  5.831 21.800  5.870 21.804 /
\plot  5.870 21.804  5.903 21.806 /
\plot  5.903 21.806  5.933 21.808 /
\putrule from  5.933 21.808 to  5.958 21.808
\putrule from  5.958 21.808 to  5.984 21.808
\putrule from  5.984 21.808 to  6.009 21.808
\plot  6.009 21.808  6.035 21.806 /
\plot  6.035 21.806  6.064 21.804 /
\plot  6.064 21.804  6.098 21.802 /
\plot  6.098 21.802  6.136 21.797 /
\plot  6.136 21.797  6.183 21.793 /
\plot  6.183 21.793  6.236 21.787 /
\plot  6.236 21.787  6.299 21.778 /
\plot  6.299 21.778  6.371 21.768 /
\plot  6.371 21.768  6.452 21.755 /
\plot  6.452 21.755  6.545 21.740 /
\plot  6.545 21.740  6.646 21.721 /
\plot  6.646 21.721  6.754 21.700 /
\plot  6.754 21.700  6.869 21.673 /
\plot  6.869 21.673  6.985 21.643 /
\plot  6.985 21.643  7.108 21.605 /
\plot  7.108 21.605  7.224 21.565 /
\plot  7.224 21.565  7.332 21.522 /
\plot  7.332 21.522  7.434 21.476 /
\plot  7.434 21.476  7.529 21.431 /
\plot  7.529 21.431  7.616 21.385 /
\plot  7.616 21.385  7.696 21.336 /
\plot  7.696 21.336  7.770 21.289 /
\plot  7.770 21.289  7.840 21.241 /
\plot  7.840 21.241  7.906 21.192 /
\plot  7.906 21.192  7.967 21.143 /
\plot  7.967 21.143  8.024 21.097 /
\plot  8.024 21.097  8.077 21.048 /
\plot  8.077 21.048  8.130 20.999 /
\plot  8.130 20.999  8.177 20.953 /
\plot  8.177 20.953  8.221 20.908 /
\plot  8.221 20.908  8.263 20.866 /
\plot  8.263 20.866  8.299 20.826 /
\plot  8.299 20.826  8.333 20.790 /
\plot  8.333 20.790  8.361 20.756 /
\plot  8.361 20.756  8.386 20.729 /
\plot  8.386 20.729  8.405 20.707 /
\plot  8.405 20.707  8.420 20.690 /
\plot  8.420 20.690  8.441 20.663 /
%
%
\plot  8.270 20.798  8.441 20.663  8.354 20.863 /
%
%
%
\linethickness= 0.500pt
\setplotsymbol ({\thinlinefont .})
\plot  3.467 20.663  3.469 20.661 /
\plot  3.469 20.661  3.473 20.657 /
\plot  3.473 20.657  3.480 20.648 /
\plot  3.480 20.648  3.490 20.635 /
\plot  3.490 20.635  3.505 20.618 /
\plot  3.505 20.618  3.526 20.595 /
\plot  3.526 20.595  3.552 20.568 /
\plot  3.552 20.568  3.581 20.534 /
\plot  3.581 20.534  3.617 20.496 /
\plot  3.617 20.496  3.655 20.455 /
\plot  3.655 20.455  3.700 20.411 /
\plot  3.700 20.411  3.747 20.364 /
\plot  3.747 20.364  3.797 20.318 /
\plot  3.797 20.318  3.850 20.269 /
\plot  3.850 20.269  3.907 20.218 /
\plot  3.907 20.218  3.967 20.170 /
\plot  3.967 20.170  4.030 20.119 /
\plot  4.030 20.119  4.098 20.070 /
\plot  4.098 20.070  4.170 20.019 /
\plot  4.170 20.019  4.248 19.971 /
\plot  4.248 19.971  4.331 19.922 /
\plot  4.331 19.922  4.420 19.876 /
\plot  4.420 19.876  4.517 19.829 /
\plot  4.517 19.829  4.621 19.782 /
\plot  4.621 19.782  4.733 19.738 /
\plot  4.733 19.738  4.851 19.696 /
\plot  4.851 19.696  4.974 19.657 /
\plot  4.974 19.657  5.093 19.626 /
\plot  5.093 19.626  5.207 19.600 /
\plot  5.207 19.600  5.317 19.577 /
\plot  5.317 19.577  5.419 19.558 /
\plot  5.419 19.558  5.512 19.541 /
\plot  5.512 19.541  5.594 19.528 /
\plot  5.594 19.528  5.668 19.520 /
\plot  5.668 19.520  5.732 19.511 /
\plot  5.732 19.511  5.785 19.505 /
\plot  5.785 19.505  5.831 19.501 /
\plot  5.831 19.501  5.870 19.497 /
\plot  5.870 19.497  5.903 19.494 /
\plot  5.903 19.494  5.933 19.492 /
\putrule from  5.933 19.492 to  5.958 19.492
\putrule from  5.958 19.492 to  5.984 19.492
\putrule from  5.984 19.492 to  6.009 19.492
\plot  6.009 19.492  6.035 19.494 /
\plot  6.035 19.494  6.064 19.497 /
\plot  6.064 19.497  6.098 19.499 /
\plot  6.098 19.499  6.136 19.503 /
\plot  6.136 19.503  6.183 19.507 /
\plot  6.183 19.507  6.236 19.514 /
\plot  6.236 19.514  6.299 19.522 /
\plot  6.299 19.522  6.371 19.533 /
\plot  6.371 19.533  6.452 19.545 /
\plot  6.452 19.545  6.545 19.560 /
\plot  6.545 19.560  6.646 19.579 /
\plot  6.646 19.579  6.754 19.600 /
\plot  6.754 19.600  6.869 19.628 /
\plot  6.869 19.628  6.985 19.657 /
\plot  6.985 19.657  7.108 19.696 /
\plot  7.108 19.696  7.224 19.736 /
\plot  7.224 19.736  7.332 19.778 /
\plot  7.332 19.778  7.434 19.825 /
\plot  7.434 19.825  7.529 19.869 /
\plot  7.529 19.869  7.616 19.916 /
\plot  7.616 19.916  7.696 19.964 /
\plot  7.696 19.964  7.770 20.011 /
\plot  7.770 20.011  7.840 20.060 /
\plot  7.840 20.060  7.906 20.108 /
\plot  7.906 20.108  7.967 20.157 /
\plot  7.967 20.157  8.024 20.204 /
\plot  8.024 20.204  8.077 20.252 /
\plot  8.077 20.252  8.130 20.301 /
\plot  8.130 20.301  8.177 20.348 /
\plot  8.177 20.348  8.221 20.392 /
\plot  8.221 20.392  8.263 20.434 /
\plot  8.263 20.434  8.299 20.475 /
\plot  8.299 20.475  8.333 20.510 /
\plot  8.333 20.510  8.361 20.544 /
\plot  8.361 20.544  8.386 20.572 /
\plot  8.386 20.572  8.405 20.593 /
\plot  8.405 20.593  8.420 20.610 /
\plot  8.420 20.610  8.441 20.637 /
%
%
\plot  8.354 20.437  8.441 20.637  8.270 20.502 /
%
%
%
\put{$g$} [lB] at 6 22.119
%
%
\put{$h$} [lB] at  6 19.022
%
%
\put{$c  ~\circlearrowleft$} [lB] at  5.5 20.6
%
%
\linethickness= 0.500pt
\setplotsymbol ({\thinlinefont .})
\ellipticalarc axes ratio  0.068:0.068  360 degrees 
	from  8.562 20.585 center at  8.494 20.585
\linethickness=0pt
\putrectangle corners at  3.382 22.373 and  8.579 18.927
\endpicture}
\end{center}\end{figure}
and a boundary condition $h = \partial c \cdot g$, which can be viewed as a `2-cell' from $g$ to $h$, labelled $(c,g)$:
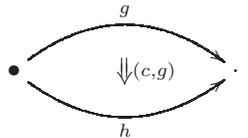
$$\xymatrix@+10pt{\bullet\rrtwocell<7>^g_h{~~~(c,g)} &&. }.$$
The use of marked vertices is, in fact unnecessary.  It can be avoided, but at the cost of repeating information, or of introducing a moderate amount of theory.  If we included expressions for each possible `marking', then any one of them could be deduced, by change of base point' from any other.  We would more naturally then use a groupoid based intuition.  The lack of `naturality' is the price we pay for sticking with a more group based system.

\medskip

As our main initial use of formal $\mathcal{C}$-maps will be in low dimensions, we will first describe them for closed 1-manifolds, then for surfaces, etc.

\medskip

Let $C_n$ denote an oriented $n$-circuit, that is, a triangulated oriented circle with $n$-edges and a choice of start-vertex.  A \emph{formal $\mathcal{C}$-map} on $C_n$ is a sequence of elements of $P$, $\mathbf{g} = (g_1,\ldots, g_n)$, thought of as labelling the edges in turn.  We will also call this a \emph{formal $\mathcal{C}$-circuit}. Two formal $\mathcal{C}$-circuits will be isomorphic if there is a simplicial isomorphism between the underlying circuits preserving the orientation and labelling.

If $S$ is a closed 1-manifold, it will be a $k$-fold disjoint union of circles and an oriented triangulation of $S$ gives a family of $n$-circuits for varying $n$.  A \emph{formal $\mathcal{C}$-map} on $S$ will be a family of formal $\mathcal{C}$-maps on the various $C_n$s.  (This includes the empty family as an instance where $S$ is the empty 1-manifold.)\label{1Cmap} 

It will be technically useful to have chosen an ordering of the vertices in any 1-manifold or, later, cobordism / triangulated surface.  This ordering may be a total order, in which case it can be used to replace the orientation, but a partial order in which the vertices of each simplex  and equally the base points of components, are totally ordered, will suffice.  The main initial reason for this imposition of an order is that it allows us to handle disjoint unions of circuits, etc., in an unambiguous way in our notation, but for most of the time it is merely for convenience.  

With such an order on the vertices of a 1-manifold, we have that a formal $\mathcal{C}$-map on it is able to be written as an \emph{ordered} family of formal $\mathcal{C}$-circuits, that is, a list of lists of elements of $P$.  Of course, the end result depends  on that order and care must be taken with this, just as care needs to be taken with the order of the constituent spaces in a vector product decomposition - and for the same reasons.

Given two formal $\mathcal{C}$-maps $\mathbf{g}$ on $S_1$, $\mathbf{h}$ on $S_2$, we can take their disjoint union to obtain a $\mathcal{C}$-map $\mathbf{g}\sqcup \mathbf{h}$ on $S_1 \sqcup S_2$. We note that $\mathbf{g}\sqcup \mathbf{h}$ and $\mathbf{h}\sqcup \mathbf{g}$ are not identical, merely `isomorphic', via an action of the symmetric group of suitable order, but, of course, this can be handled in the usual ways, depending to some extent on taste, for instance via the standard technical machinery of symmetric monoidal categories.  In this paper however we will tend to avoid the detailed technicalities where they are inessential to our aim of building the intuition of what is going on.

If we have a closed oriented triangulated 1-manifold, we mark each initial vertex of each edge as such.  If we reverse the orientation on the 1-manifold, we reverse the order of the elements in the sequence, and invert each in turn. If we change the start vertex, we merely cyclically permute the sequence in the obvious way.  

Note that the `top group' $C$ of $\mathcal{C}$ plays no role in this dimension.

Now let $M$ be an oriented (triangulated) cobordism between two such 1-manifolds $S_0$ and $S_1$, and suppose given formal $\mathcal{C}$-maps, $\mathbf{g}_0$, $\mathbf{g}_1$, on $S_0$ and $S_1$ respectively. A \emph{formal $\mathcal{C}$-map}, $\mathbf{F}$, on $M$ consists of a family of elements $\{c_t\}$ of $C$, indexed by the triangles $t$ of $M$, a family, $\{p_e\}$ of elements of $P$ indexed by the edges of $M$ and for each $t$, a choice of base vertex, $b(t)$, such that the boundary condition below is satisfied:\\
in any triangle $t$,
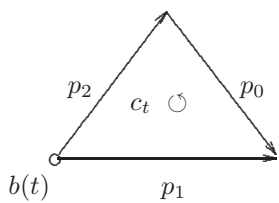
\begin{center}
\font\thinlinefont=cmr5
\begingroup\makeatletter\ifx\SetFigFont\undefined
\def\x#1#2#3#4#5#6#7\relax{\def\x{#1#2#3#4#5#6}}%
\expandafter\x\fmtname xxxxxx\relax \def\y{splain}%
\ifx\x\y   
\gdef\SetFigFont#1#2#3{%
  \ifnum #1<17\tiny\else \ifnum #1<20\small\else
  \ifnum #1<24\normalsize\else \ifnum #1<29\large\else
  \ifnum #1<34\Large\else \ifnum #1<41\LARGE\else
     \huge\fi\fi\fi\fi\fi\fi
  \csname #3\endcsname}%
\else
\gdef\SetFigFont#1#2#3{\begingroup
  \count@#1\relax \ifnum 25<\count@\count@25\fi
  \def\x{\endgroup\@setsize\SetFigFont{#2pt}}%
  \expandafter\x
    \csname \romannumeral\the\count@ pt\expandafter\endcsname
    \csname @\romannumeral\the\count@ pt\endcsname
  \csname #3\endcsname}%
\fi
\fi\endgroup
\mbox{\beginpicture
\setcoordinatesystem units <0.750000cm,0.750000cm>
\unitlength=0.750000cm
\linethickness=1pt
\setplotsymbol ({\makebox(0,0)[l]{\tencirc\symbol{'160}}})
\setshadesymbol ({\thinlinefont .})
\setlinear
%
%
\linethickness= 0.500pt
\setplotsymbol ({\thinlinefont .})
\plot  3.6 21.220  5.5 23.732 /
%
%
\plot  5.394 23.531  5.5 23.732  5.309 23.594 /
%
%
%
\linethickness= 0.500pt
\setplotsymbol ({\thinlinefont .})
\plot  5.50 23.732  7.436 21.220 /
%
%
\plot  7.265 21.355  7.436 21.220  7.349 21.420 /
%
%
%
\linethickness= 0.500pt
\setplotsymbol ({\thinlinefont .})
\putrule from  3.573 21.139 to  7.436 21.139
%
%
\plot  7.224 21.086  7.436 21.139  7.224 21.192 /
%
%
%
\put{$p_2$} [lB] at  3.75 22.3
%
%
\put{$p_0$} [lB] at  6.8 22.3
%
%
\put{$p_1$} [lB] at  5.4 20.5
%
%
\put{ $c_t ~ \circlearrowleft$} [lB] at  4.7  22
%
%
\linethickness= 5.500pt
\setplotsymbol ({\thinlinefont .})
\ellipticalarc axes ratio  0.093:0.093  360 degrees 
	from  3.613 21.114 center at  3.520 21.114
%
%
\put{$b(t)$} [lB] at  2.7 20.5

\linethickness=0pt
\putrectangle corners at  3.410 23.757 and  7.461 20.701
\endpicture}

\end{center}
we have $$\partial c_t = p_1 p_0^{-1}p_2^{-1}.$$

We call such a formal $\mathcal{C}$-map  on $M$ a \emph{formal $\mathcal{C}$-cobordism} from $(S_0,\mathbf{g}_0)$  to $(S_1,\mathbf{g}_1)$  if it restricts to these formal $\mathcal{C}$-maps on the boundary 1-manifolds.  We will denote it $(M, \mathbf{F})$.

To be able to handle manipulation of formal $\mathcal{C}$-cobordisms `up to equivalence', so as to be able to absorb choices of triangulation, base vertices, etc. and eventually to pass to regular cellular decompositions, we need to consider triangulations of 3-dimensional simplicial complexes and formal $\mathcal{C}$-maps on these.   We, in fact, can use a common generalisation to all simplicial complexes.

\

\textbf{Definition}

Let $K$ be a simplicial complex.  A \emph{(simplicial) formal $\mathcal{C}$-map}, $\lambda$, on $K$ consists of families of elements\\
(i) $\{c_t\}$ of $C$, indexed by the set, $K_2$, of 2-simplices of $K$, 
\\(ii) $\{p_e\}$ of $P$, indexed by the set of 1-simplices, $K_1$, of $K$\\
 and a partial order on the vertices of $K$, so that each simplex is totally ordered (this replaces the orientation and gives start vertices to all edges and triangles without problem).  The assignments of $c_t$ and $p_e$, etc. are to satisfy\\
(a) the boundary condition  $$\partial c_t = p_1 p_0^{-1}p_2^{-1},$$
where the vertices of $t$, labelled $v_0$, $v_1$, $v_2$ in order, determine the numbering of the opposite edges, e.g., $e_0$ is between $v_1$ and $v_2$, and $p_{e_i}$ is abbreviated to $p_i$;\\
and \\
(b) the cocycle condition:\\
in a tetrahedron yielding two composite faces
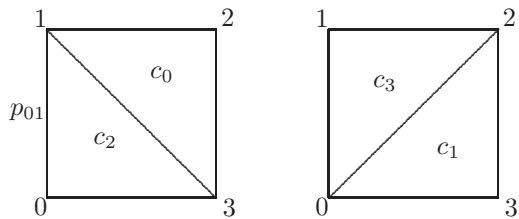
\begin{figure}[h]\begin{center}$\begin{array}{cc}
\font\thinlinefont=cmr5
\begingroup\makeatletter\ifx\SetFigFont\undefined
\def\x#1#2#3#4#5#6#7\relax{\def\x{#1#2#3#4#5#6}}%
\expandafter\x\fmtname xxxxxx\relax \def\y{splain}%
\ifx\x\y   
\gdef\SetFigFont#1#2#3{%
  \ifnum #1<17\tiny\else \ifnum #1<20\small\else
  \ifnum #1<24\normalsize\else \ifnum #1<29\large\else
  \ifnum #1<34\Large\else \ifnum #1<41\LARGE\else
     \huge\fi\fi\fi\fi\fi\fi
  \csname #3\endcsname}%
\else
\gdef\SetFigFont#1#2#3{\begingroup
  \count@#1\relax \ifnum 25<\count@\count@25\fi
  \def\x{\endgroup\@setsize\SetFigFont{#2pt}}%
  \expandafter\x
    \csname \romannumeral\the\count@ pt\expandafter\endcsname
    \csname @\romannumeral\the\count@ pt\endcsname
  \csname #3\endcsname}%
\fi
\fi\endgroup
\mbox{\beginpicture
\setcoordinatesystem units <0.5000cm,0.5000cm>
\unitlength=0.5000cm
\linethickness=1pt
\setplotsymbol ({\makebox(0,0)[l]{\tencirc\symbol{'160}}})
\setshadesymbol ({\thinlinefont .})
\setlinear
%
%
\put{2} [lB] at  9.129 22.754
%
%
\linethickness= 0.500pt
\setplotsymbol ({\thinlinefont .})
\plot  4.498 22.648  8.970 18.203 /
%
%
\put{0} [lB] at  4.155 17.727
%
%
\put{1} [lB] at  4.155 22.727
%
%
\linethickness= 0.500pt
\setplotsymbol ({\thinlinefont .})
\putrule from  4.498 18.176 to  4.498 22.648
\putrule from  4.498 22.648 to  8.996 22.648
\putrule from  8.996 22.648 to  8.996 18.176
\putrule from  8.996 18.176 to  4.498 18.176
%
%
\put{3} [lB] at  9.182 17.700
%
%
\put{$p_{01}$} [lB] at  3.520 20.373
%
%
\put{$c_2$} [lB] at  5.743 19.605
%
%
\put{$c_0$} [lB] at  7.250 21.351
\linethickness=0pt
\putrectangle corners at  3.520 22.754 and  9.390 17.700
\endpicture}
\quad&\quad
\font\thinlinefont=cmr5
\begingroup\makeatletter\ifx\SetFigFont\undefined
\def\x#1#2#3#4#5#6#7\relax{\def\x{#1#2#3#4#5#6}}%
\expandafter\x\fmtname xxxxxx\relax \def\y{splain}%
\ifx\x\y   
\gdef\SetFigFont#1#2#3{%
  \ifnum #1<17\tiny\else \ifnum #1<20\small\else
  \ifnum #1<24\normalsize\else \ifnum #1<29\large\else
  \ifnum #1<34\Large\else \ifnum #1<41\LARGE\else
     \huge\fi\fi\fi\fi\fi\fi
  \csname #3\endcsname}%
\else
\gdef\SetFigFont#1#2#3{\begingroup
  \count@#1\relax \ifnum 25<\count@\count@25\fi
  \def\x{\endgroup\@setsize\SetFigFont{#2pt}}%
  \expandafter\x
    \csname \romannumeral\the\count@ pt\expandafter\endcsname
    \csname @\romannumeral\the\count@ pt\endcsname
  \csname #3\endcsname}%
\fi
\fi\endgroup
\mbox{\beginpicture
\setcoordinatesystem units <0.5000cm,0.5000cm>
\unitlength=0.5000cm
\linethickness=1pt
\setplotsymbol ({\makebox(0,0)[l]{\tencirc\symbol{'160}}})
\setshadesymbol ({\thinlinefont .})
\setlinear
%
%
\linethickness= 0.500pt
\setplotsymbol ({\thinlinefont .})
\plot  4.551 18.212  9.023 22.657 /
%
%
\put{0} [lB] at  4.155 17.727
%
%
\put{1} [lB] at  4.155 22.727
%
%
\put{2} [lB] at  9.129 22.754
%
%
\put{3} [lB] at  9.182 17.700
%
%
\put{$c_3$} [lB] at  5.690 21.114
%
%
\put{$c_1$} [lB] at  7.383 19.340
%
%
\linethickness= 0.500pt
\setplotsymbol ({\thinlinefont .})
\putrule from  4.498 18.176 to  4.498 22.648
\putrule from  4.498 22.648 to  8.996 22.648
\putrule from  8.996 22.648 to  8.996 18.176
\putrule from  8.996 18.176 to  4.498 18.176
\linethickness=0pt
\putrectangle corners at  4.155 23.072 and  9.390 17.700
\endpicture}
\end{array}$
\end{center}\end{figure}\\
we have $$c_2{}^{p_{01}}c_0 = c_1c_3.$$

\medskip

\textbf{ Explanation of the cocycle condition.}

The left hand and right hand sides of the cocycle condition have the same boundary, namely the boundary of the square, so $c_2{}^{p_{01}}c_0(c_1c_3)^{-1}$ is a cycle.   
A crossed module has elements in dimensions 1 and 2, but nothing in dimension 3, therefore just as the case where $B = BG$ for $G$ a group led to a cocycle condition in dimension 2, so when labelling with elements of a crossed module, we should expect the cocycle condition to be a `tetrahedral equation', hence in dimension 3. In future developments, it may be useful to replace a crossed module $\partial :C\to P$ by a longer `crossed complex', $C_3\to C_2\to C_1$, and then we would expect to have a slightly more complex labelling and a correspondingly adjusted cocycle condition.

When `integrating ' a labelling over a surface corresponding to three faces of a tetrahedron, the composite label is on the remaining face, so given a formal $\mathcal{C}$-map on the tetrahedron, and a specification of $p_{01}$, any one of $c_0, \ldots, c_3$ is  determined by the others. (For example if all but $c_0$ are given, then 
$${}^{p_{01}}c_0 = c_2^{-1}c_1c_3,$$and acting throughout with $p_{01}^{-1}$ yields $c_0$.)

A third related view is that coming from the homotopy addition lemma, \cite{RBRS:2005}, which loosely says that any one face of an $n$-simplex is a (suitably defined) composite of the others.

 \medskip

The way the cocycle condition will be used is to show that the manner in which the interior of a polyhedral disc is triangulated in a formal $\mathcal{C}$-map yields a single label on that polyhedral 2-cell that is independent of the actual decomposition used, although dependent on labellings up to a notion of equivalence to be given shortly.  It replaces the use of `moves' on the triangulation in this respect. This will allow us to simplify formal $\mathcal{C}$-maps from the above simplicial form to a neat cellular form, see later.  The reader may already see the basic idea of attaching elements of $P$ to edges of a cellular decomposition and elements of $C$ to the 2-cells satisfying a boundary condition.  The one more subtle but important point is however the handling of the cocycle condition, which takes a bit of more care.

\medskip

We will restrict attention to 1+1 HQFTs and to formal $\mathcal{C}$-maps on 1-manifolds, surfaces and 3-manifolds.  If a higher dimensional theory was being considered based on $B$-manifolds of dimension $d$, the cocycle condition would naturally occur in dimension $d + 2$.  In that case, the natural coefficients would be in one of the higher dimensional analogues of a crossed module such as crossed complexes, or truncated hypercrossed complexes (or, equivalently, simplicial groups). References to these notions can be found in  Baues, \cite{HJB:AH,HJB:4D}, Brown, Higgins  and Sivera,  \cite{RBRS:2005}, Carrasco and Cegarra, \cite{carrasco-cegarra}, Porter, \cite{TP:n-types,tqft2,cubo}, etc, ..., depending on the level of generality desired. This will be explored more fully in \cite{TP:fhqft2}.
\medskip

\textbf{Equivalence of formal $\mathcal{C}$-maps}

Suppose $X$ is a polyhedron with a given non-empty family of base points $\mathbf{m} = \{m_i\}$, and $K_0$, $K_1$ two triangulations of $X$, i.e., $K_0$ and $K_1$ are simplicial complexes with geometric realisations homeomorphic to $X$ (by specified homeomorphisms) with the given base points among the vertices of the triangulation. 
\medskip

\textbf{Definition}

 Given two formal $\mathcal{C}$-maps $(K_0,\lambda_0),$ $(K_1,\lambda_1)$, then we say they are \emph{equivalent}  if there is a triangulation, $T$, of $X\times I$ extending $K_0$ and $K_1$ on $X \times \{0\}$ and $X\times \{1\}$ respectively, and a formal $\mathcal{C}$-map, $\Lambda$, on $T$ extending the given ones on the two ends and respecting the base points, in the sense that $T$ contains a subdivided $\{m_i\}\times I$ for each basepoint $m_i$ and $\Lambda $ assigns the identity element $1_P$ of $P$ to each 1-simplex of $\{m_i\}\times I$.

\,
	We will use the term `ordered simplicial complex' for a simplicial complex, $K$, together with a partial order on its set of vertices such that the vertices in any simplex of $K$ form a totally ordered set.  If we give the unit interval, $I$, the standard structure of an ordered simplicial complex with $0 < 1$, then the cylinder $|K|\times I$ has a canonical triangulation as an ordered simplicial complex and we will write $K\times I$ for this.  We will assume some base points are given.  We can, for instance, consider all vertices as base points.

If we are given two formal $\mathcal{C}$-maps defined on the same ordered $K$, $(K,\lambda_0),$ and $(K,\lambda_1)$, we say they are \emph{simplicially homotopic}  as formal maps, if there is a formal $\mathcal{C}$-map defined on the ordered simplicial complex $K\times I$ extending them both and respecting base points.
\begin{lemma}\hspace*{1mm}\\ \label{refl-trans}
Equivalence is an equivalence relation.
\end{lemma}

\textbf{Proof}

This is mostly routine. Transitivity and symmetricity are  easy, whilst reflexivity merely requires the construction of the standard triangulation of $X\times I$, followed by the obvious construction of a formal map on it. The details are omitted. \hfill $\square$
\medskip

  Equivalence combines  the intuition of the geometry of triangulating a (topological) homotopy,  where the triangulations of the two ends may differ,  with some idea of a  combinatorially defined simplicial homotopy of formal maps.

\begin{lemma}~\\
If   $(K,\lambda_0),$ and $(K,\lambda_1)$ are two formal  $\mathcal{C}$-maps, which are simplicially homotopic, then they are equivalent. \hspace*{1cm}\hfill $\square$
\end{lemma}

The proof is immediate from the definition and is omitted.

There are several possible proofs of the following result. We give one that is amongst the longer ones as it illustrates more clearly the processes of combination of labellings of simplices given by a formal $\mathcal{C}$-map by explicitly constructing the required extension.
\begin{proposition}~\\
Given a simplicial complex, $K$, with geometric realisation $X = |K|$, and a subdivision $K^\prime$ of $K$.\\
(a)  Suppose $\lambda$ is a formal $\mathcal{C}$-map on $K$, then there is a formal $\mathcal{C}$-map, $\lambda^\prime$ on $K^\prime$ equivalent to $\lambda$.\\
(b) Suppose $\lambda^\prime$ is a formal $\mathcal{C}$-map on $K^\prime$, then there is a formal $\mathcal{C}$-map, $\lambda$ on $K$ equivalent to $\lambda^\prime$.
\end{proposition}

\textbf{Proof}

(The proof that follows is moderately `technical', so if the reader is willing to accept the results as `clear', it can safely be omitted or `skimmed' at first reading . This is also the case for several other proofs in the following pages.The method of proof is clear however. )

We first need to construct a good triangulation, $T$, of the cylinder $X\times I$ extending $K$ on $X \times \{0\}$ and $K^\prime$ on $X \times \{1\}$, then given a formal $\mathcal{C}$-map defined on one end extend it to one on the whole triangulated cylinder so that restricting to the opposite end gives the required equivalent formal $\mathcal{C}$-map.  We can assume that both $K$ and $K^\prime$ are ordered and, for convenience, will assume that the vertices in $K^\prime$  that are also in $K$ have the same order as there, whilst those new vertices in $K^\prime$ are ordered after those in $K$. (The other possibilities can be reduced to this using Proposition 2 and Lemma 1 if need be.)

We first consider the trivial situation in which $K$= $K^\prime$, so we get an obvious triangulation of $X\times I$.  If $\sigma=\langle v_0, \ldots , v_n\rangle$ is a simplex in $K$, then $$\langle (v_0,0), \ldots , (v_k,0),(v_k,1),\ldots , (v_n,1) \rangle$$ gives a simplex in this triangulation of the cylinder and simplices of this general form generate that triangulation.  Such a simplex is given by an initial segment $\sigma_k =  \langle v_0, \ldots , v_k\rangle$ of the ordered set $\{v_0< \ldots < v_n\}$, which then determines $\sigma^c_{n-k} = \langle v_k, \ldots , v_n\rangle$, such that the corresponding $(n+1)$-simplex of $X \times I$ is the join of $\sigma_k \times \{0\}$ and $\sigma^c_{n-k}\times \{1\}$. We will call these simplices `large' simplices.  The triangulation we need in general will be obtained by subdividing this simple one.  

Within $X\times I$, we have $|\sigma|\times I$, which is a prism with $\sigma$ on its base and a possibly subdivided version of $\sigma$ on its top. If  $\sigma$  is not subdivided within $K^\prime$, then we use the simple case discussed in the previous paragraph to triangulate $|\sigma|\times I$. If it is subdivided, then we look at the joins, $\sigma_k \times \{0\} * \sigma^c_{n-k}\times \{1\}$.  The general picture is that $\sigma^c_{n-k}$ may be subdivided to obtain $\sigma^{c,\prime}_{n-k}$, say, and so we triangulate $\sigma_k \times \{0\} * \sigma^c_{n-k}\times \{1\}$ using the join
 $\sigma_k \times \{0\} * \sigma^{c,\prime}_{n-k}\times \{1\}$.  The neat way to do this is by induction up the skeleton of $K$ so handling induction on $n$ first and then considering $k = n, n-1, ..., 0$ in turn so that each subcomplex is built neatly on ones that have been previously constructed.  The essential step, however, is always the same : the triangulation subdivides the top part of the join $\sigma_k \times \{0\} * \sigma^c_{n-k}\times \{1\}$ when necessary. (We leave the detailed induction to the reader.  In the cases that we will need here, $n$ is small as it is never bigger than 4 in any argument used in this paper, so a detailed induction seems `overkill', but the result is true without restriction on the dimensions.)

We now start building a formal $\mathcal{C}$-map, $\Lambda$, on $T$ extending $\lambda$ on $K$ (considered as ``$K \times \{0\}$'' within $T$).  If an edge, $\langle u,u^\prime\rangle$, of $K$ does not get subdivided in $K^\prime$, then in the square with base $\langle (u,0),(u^\prime,0)\rangle$, we label the vertical edges with $1_P$ and the 2-simplices with $1_C$.  The boundary rule then determines that the edge $\langle (u,0),(u^\prime,1)\rangle$ is labelled by the same element as the base.  A repeat use of this argument then shows that $\langle (u,1),(u^\prime,1)\rangle$ is again labelled by that same element.    

A similar thing happens over a 2-simplex of $K$ not involving any new vertex $v$.  We can consider the edges of this simplex as  having already been handled, so if the simplex is $\langle u_0,u_1,u_2\rangle$, labelled $c\in C$, we look at the tetrahedron $\langle (u_0,0), (u_1,0), (u_2,0), (u_2,1)\rangle$ and check its faces:
\begin{itemize}
\item $\langle (u_0,0), (u_1,0), (u_2,0)\rangle$ is in the base, so is handled by $\lambda$, and is labelled $c$;
\item $(*) = \langle (u_0,0), (u_1,0), (u_2,1)\rangle$, we do not know yet as it is a `free face';
\item $\langle (u_0,0),  (u_2,0), (u_2,1)\rangle$ has already been handled, as it is part of a face over $\langle u_0, u_2\rangle$ and has been labelled $1_C$ by our previous step;
\item $\langle (u_1,0), (u_2,0),(u_2,1)\rangle$ has likewise been labelled $1_C$.
\end{itemize}
Now the cocycle condition implies that $(*)$ must also be labelled $c$.

Examination of the 3-simplex $\langle (u_0,0), (u_1,0), (u_1,1), (u_2,1)\rangle$  next shows that it has just one face `free' and there is a unique value, $c$ again, with which it can be  labelled consistently with the cocycle condition.  The final simplex $$\langle (u_0,0), (u_0,1), (u_1,1), (u_2,1)\rangle$$  is similar and causes no problem.

This case was, of course, easy, but it suggests the general process. In general, there will be no bound on the dimension of the base (labelled) simplex and no difficulty in creating the formal $\mathcal{C}$-map, $\Lambda$ on the corresponding prism as above dimension 2, there are no more labels.  The new $\lambda^\prime$, i.e., $\Lambda$ restricted to the other end, will agree with $\lambda$ on these simplices.  Note that the cocycle condition was all that was needed for this, given the fact that we used $1_C$ to label the vertical faces of the prism.
 
We next turn to those $\sigma$ in $K$ which will be subdivided in $K^\prime$. We start with $\sigma$ of dimension 1, so $\sigma = \langle u_0,u_1\rangle$ and we have new vertices $w_1,w_2, \ldots, w_l$, which we will assume ordered as given, in the subdivided edge $\sigma^\prime$ of $K^\prime$.  (For simplicity we will assume the vertices occur in order along the edge.  A reversal of order just corresponds to replacing a label $p$ in $P$ by its inverse $p^{-1}$, and so changes nothing of importance.)  In the face $|\sigma|\times I$, we have two `large' simplices: $\langle (u_0,0),(u_1,0),(u_1,1)\rangle,  $ which is not subdivided, and  a second one $\langle (u_0,0),(u_0,1),(u_1,1)\rangle$, which is subdivided into $l+1$ parts: $\langle (u_0,0),(u_0,1),(w_1,1)\rangle$, $\langle (u_0,0),(w_1,1),(w_2,1)\rangle$, ..., $\langle (u_0,0),(w_{l-1},1),(w_l,1)\rangle$ and finally  $\langle (u_0,0),(u_1,1),(w_l,1)\rangle$.

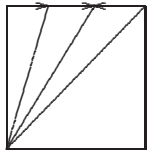
\begin{center}
\font\thinlinefont=cmr5
\begingroup\makeatletter\ifx\SetFigFont\undefined
\def\x#1#2#3#4#5#6#7\relax{\def\x{#1#2#3#4#5#6}}%
\expandafter\x\fmtname xxxxxx\relax \def\y{splain}%
\ifx\x\y   
\gdef\SetFigFont#1#2#3{%
  \ifnum #1<17\tiny\else \ifnum #1<20\small\else
  \ifnum #1<24\normalsize\else \ifnum #1<29\large\else
  \ifnum #1<34\Large\else \ifnum #1<41\LARGE\else
     \huge\fi\fi\fi\fi\fi\fi
  \csname #3\endcsname}%
\else
\gdef\SetFigFont#1#2#3{\begingroup
  \count@#1\relax \ifnum 25<\count@\count@25\fi
  \def\x{\endgroup\@setsize\SetFigFont{#2pt}}%
  \expandafter\x
    \csname \romannumeral\the\count@ pt\expandafter\endcsname
    \csname @\romannumeral\the\count@ pt\endcsname
  \csname #3\endcsname}%
\fi
\fi\endgroup
\mbox{\beginpicture
\setcoordinatesystem units <0.750000cm,0.750000cm>
\unitlength=0.750000cm
\linethickness=1pt
\setplotsymbol ({\makebox(0,0)[l]{\tencirc\symbol{'160}}})
\setshadesymbol ({\thinlinefont .})
\setlinear
%
%
\linethickness= 0.500pt
\setplotsymbol ({\thinlinefont .})
\putrule from  4.551 23.150 to  5.292 23.150
%
%
\plot  5.080 23.097  5.292 23.150  5.080 23.203 /
%
%
%
\linethickness= 0.500pt
\setplotsymbol ({\thinlinefont .})
\putrule from  4.551 20.610 to  7.013 20.610
\putrule from  7.013 20.610 to  7.013 23.150
%
%
\linethickness= 0.500pt
\setplotsymbol ({\thinlinefont .})
\plot  6.085 23.097  6.085 23.097 /
%
%
\linethickness= 0.500pt
\setplotsymbol ({\thinlinefont .})
\putrule from  7.013 23.150 to  6.085 23.150
%
%
\plot  6.297 23.203  6.085 23.150  6.297 23.097 /
%
%
%
\linethickness= 0.500pt
\setplotsymbol ({\thinlinefont .})
\putrule from  5.266 23.150 to  6.113 23.150
%
%
\plot  5.901 23.097  6.113 23.150  5.901 23.203 /
%
%
%
\linethickness= 0.500pt
\setplotsymbol ({\thinlinefont .})
\putrule from  4.551 20.610 to  4.551 23.150
%
%
\linethickness= 0.500pt
\setplotsymbol ({\thinlinefont .})
\plot  4.551 20.61  7.013 23.150 /
%
%
\linethickness= 0.500pt
\setplotsymbol ({\thinlinefont .})
\plot  4.551 20.61  6.113 23.150 /
%
%
\linethickness= 0.500pt
\setplotsymbol ({\thinlinefont .})
\plot  4.551 20.61  5.292 23.15 /
%
%
\put{case  $l = 2$} [lB] at  9.2 21.696
\linethickness=0pt
\putrectangle corners at  4.473 23.228 and  9.432 20.559
\endpicture}
\end{center}
As previously, we label any edge of form $\langle(v,0),(v,1)\rangle$ with $1_P \in P$. Assuming $\sigma = \langle u_0,u_1\rangle$ in $K$ was labelled $p$, the obvious thing to do is to use $1_C$ as the $C$-part of all triangles within the face and then there is a choice of labels for all but the last edge labelled in the top, in fact that edge will be $\langle (u_0,1),(w_1,1)$.  For instance, we clearly must label the edge  $\langle (u_0,0),(u_1,1)\rangle$ by $p$, then make an arbitrary choice $p_l\in P$ to label $(u_1,1),(w_l,1)\rangle$.  The cocycle condition then says $\langle(u_0,0),(w_l,1)\rangle$ is labelled $pp_l$.  Applying the same process to the next triangle along, pick some $p_{l-1}$ as a label for the top edge $\langle (w_{l-1},1),(w_l,1)\rangle$ and then the cocycle condition will give a label $x$ satisfying $xp_{l-1} = pp_l$, and so on. As the proof only needs the existence of a formal $\mathcal{C}$-map with the required properties on $X\times I$, we can simplify things and always choose $1_P$ as the label for the top edge, \emph{but} for a fuller picture of what is going on, it is important to realise that the choices could be made otherwise. Assuming the simple choice is made, each edge $\langle (u_0,0),(w_k,1)\rangle$ is labelled $p$, and finally 
$\langle (u_0,1),(w_k,1)\rangle$ is as well. (For  any choice the product across the top should give $p$ as the face is labelled $1_C$ and the vertical edges $1_P$.)

The higher dimensional cases have a similar pattern, but, of course, are more complex to describe.  We will look at $\sigma = \langle u_0,u_1,u_2\rangle$ in detail, but in higher dimensions, the fact that the cocycle condition applies as if there was a labelling by identity elements (together with a boundary condition) makes the extension to those dimensions more or less trivial.  

The possible subdivisions of a 2-simplex are much more complex than those for an edge; there can be new `interior' vertices in an old simplex, new edges that cross the simplex from vertices in the boundary and so on.  We assume $\sigma$ is labelled with $(c,p_0,p_1,p_2)$ with the obvious convention, $c\in C$, etc., satisfying the boundary condition $\partial c = p_1p_0^{-1}p_2^{-1}$, and we proceed to label the prism with base $\sigma$.  The vertical edges are labelled $1_P$ as before and the vertical faces by the previous step. The triangulation using joins, as in our previous discussion, gives a first tetrahedron $\langle (u_0,0),(u_1,0),(u_2,0),(u_2,1)\rangle$, which, as is easily seen, gives, by the cocycle condition, a labelling $(c,p_0,p_1,p_2)$ on the face $\langle (u_0,0),(u_1,0),(u_2,1)\rangle$.  The next `large' simplex is that given by $\langle (u_0,0),(u_1,0),(u_1,1),(u_2,1)\rangle$. We already have a labelling on the vertical face over $\langle (u_0,0),(u_2,0))\rangle$ with any subdivision of $\langle (u_1,1),(u_2,1)\rangle$ already used.  We also have a labelling of $\langle (u_0,0),(u_1,0),(u_1,1)\rangle$ as it also is a vertical face.  Suppose $\langle u_1,u_2\rangle$ is subdivided in $K^\prime$ and we adopt the same notation, $w_1,   \ldots $, as before. The 3-simplex  
 $\langle (u_0,0),(u_1,0),(u_1,1),(w_l,1)\rangle$ can be labelled with a $1_C$ on  $\langle (u_0,0),(u_1,1),(w_l,1)\rangle$, giving once again $(c,p_0,p_1,p_2)$ labelling  $\langle (u_0,0),(u_1,0),(w_l,1)\rangle$.  We have pushed that labelling up from the base and can now push it along the top until we get to the last 3-simplex of our subdivided `large' one.  Here we already have 3 of the 4 faces predetermined, so can use the cocycle condition to solve for the last one.

This leaves us with the triangulated version of our large simplex $$\langle (u_0,0),(u_0,1),(u_1,1),(u_2,1)\rangle.$$ If the subdivided top face has any interior vertices $(v,1)$, then label $\langle (u_0,0),(v,1)\rangle$ with $1_P$. If it has any extra edges, then $\langle(u_0,0),(v,1),(w,1)\rangle$ will be a 2-simplex and we already have labelled two of its edges, so using a value of $1_C$ on the 2-simplex yields the value on the edge $\langle(v,1),(w,1)\rangle$.  This just leaves any 2-simplices in this subdivided $\sigma$.  If $\langle v_0,v_1,v_2\rangle$ is any such, $\langle(u_0,0), (v_0,1),(v_1,1),(v_2,1)\rangle$ is a 3-simplex on which we know the labelling on all but one of the faces, so using the cocycle condition we obtain the final face and, by repeating for all such, our extended formal $\mathcal{C}$-map.  Restricting to the top face, $K^\prime$, we obtained our required $\lambda^\prime$ equivalent to $\lambda$ thus proving a).

We note that the extension can be reversed with minor alterations to prove b).   Now the top edges and faces are already labelled; we can label vertical edges as before with $1_P$ and vertical faces with $1_C$.  The cocycle condition then gives us the labels on diagonal edges and faces. The reversal of the `algorithm' proves b).
\hfill $\square$

\medskip

\textbf{Remarks}

(i) If we look at this proof in detail, we can see it as a series of nested inductions `up the skeleton' of various parts of the structure. To handle higher dimensions, we continue that process only handling $\sigma \in K_n$ when all its faces have been done, then using inverse induction and the join formulation of the triangulation as above for the case $n = 2$.

(ii)  There is a simplicial set formulation of the above in terms of the Kan complex condition on the simplicial nerve of $\mathcal{C}$.  This is useful for the extension of this theory to higher dimensions, but we have avoided its use here as  the extra technical machinery required for its development might tend to obscure the basic simplicity of the extension of the theory of \cite{turaev:hqft1} from handling a $K(G,1)$ to handling a general 2-type, $B\mathcal{C}$. We explore this more fully in the second paper, \cite{TP:fhqft2}, of this series. Similar methods were used in \cite{tqft1}.

(iii) The idea of a formal $\mathcal{C}$-map is to represent, combinatorially, the characteristic map of a $B$-manifold  or $B$-cobordism, and from this perspective, equivalent formal maps will correspond to homotopic characteristic maps.
\begin{proposition}\label{reorder}~\\
A change of partial order on the vertices of $K$, or a change in choice of start vertices for simplices, generates an equivalent formal $\mathcal{C}$-map.
\end{proposition}
\textbf{Proof}

More formally, let $K_0$ be $K$ with the given order and $K_1$ the same simplicial complex with a new ordering.  Construct a triangulation $T$ of $|K|\times I$ having $K_0$ and $K_1$ on the two ends.  (Inductively, we can suppose just one pair of elements has been transposed in the order.)  It is now easy to adapt the method of the previous proposition to extend any given $\lambda_0$ on $K_0$ over $T$ and then to restrict to get an equivalent $\lambda_1$ on $K_1$.\hfill $\square$

\medskip

Note if $\langle v_0,v_1 \rangle$ is an ordered edge of $K_0$ and, with the reordering, $\langle v_1,v_0 \rangle$ is the corresponding one in $K_1$, then if $\lambda_0$ assigns $p$ to $\langle v_0,v_1 \rangle$, $\lambda_1$ assigns $p^{-1}$ to $\langle v_1,v_0 \rangle$ as is clear for the simplest assignment scheme:
$$\xymatrix{(v_0,1)&(v_1,1)\ar[l]_{p^{-1}}\\
(v_0,0)\ar[u]_{1_P}\ar[r]_p&(v_1,0)\ar[ul]_{p^{-1}}\ar[u]_{1_P}}
$$
(The triangulation $T$ assumes here that vertices of $|K|\times	 \{1\}$ are always listed after those of $|K|\times \{0\}$.)  A similar, but more complex, observation is valid for higher dimensional simplices.  Once the use of the boundary and cocycle conditions is understood, the choice of local ordering within the triangulation easily determines the simplest choice of extension.  That extension can be perturbed or deformed by changing the choice of fillers for the 2-simplices in the faces of the prisms however.

\subsection{An interlude on combining simplices}
We can use the cocycle condition to combine formal $\mathcal{C}$-data given locally on simplices into cellular blocks, \emph{up to equivalence}.

As homotopy of characteristic maps is mirrored combinatorially by equivalence of formal maps, we can study $B$-manifolds and the resulting HQFTs by manipulating formal maps up to equivalence. We will mainly use examples to illustrate the process.

\medskip

\textbf{Examples}

(i) Suppose we have a simplicial complex $K$ and two adjacent 2-simplices with formal map data relative to the crossed module $\mathcal{C} = (C,G, \partial)$,
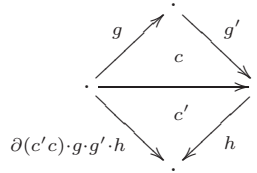
$$\xymatrix{&.\ar[dr]^{g^\prime}&\\
.\ar[ur]^{g}_{\hspace{.5cm}c}\ar[rr]\ar[dr]_{\partial(c^\prime c)\cdot g\cdot g^\prime\cdot h}^{\hspace{.5cm}c^\prime}&&.\ar[dl]^h\\
&.&
}$$
with the horizontal edge labelled $\partial c\cdot g\cdot g^\prime$.  In such a diagram, one can `compose' the two 2-simplices to get an equivalent labelling, locally, without changing the overall boundary of this subdiagram. This can be done in several ways, but notably so as to get: 
$$
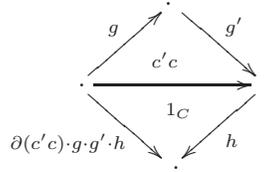
\xymatrix{&.\hspace{1mm}\ar[dr]^{g^\prime}&\\
\hspace{1mm}.\ar[ur]^{g}_{\hspace{.3cm}c^\prime c}\ar[rr]\ar[dr]_{\partial(c^\prime c)\cdot g\cdot g^\prime\cdot h}^{\hspace{.5cm}1_C}&&\hspace{1mm}.\ar[dl]^h\\
&\hspace{1mm}.&
}$$
To prove equivalence we just triangulate the square prism, label the top and bottom in the required ways and the side panels with `constant' equivalences, then there is an obvious labelling on the interior edges and faces.  Try it!

Again we could have moved   $c$ to the bottom triangle replacing it with a $1_C$. 
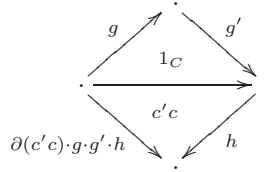
$$\xymatrix{&\hspace{1mm}.\ar[dr]^{g^\prime}&\\
\hspace{1mm}.\ar[ur]^{g}_{\hspace{.4cm}1_C}\ar[rr]\ar[dr]_{\partial(c^\prime c)\cdot g\cdot g^\prime\cdot h}^{\hspace{.3cm}c^\prime c}&&\hspace{1mm}.\ar[dl]^h\\
&\hspace{1mm}.&
}$$
 In each case, the horizontal edge gets a different labelling or `colouring' from its original one, in the first case by $\partial(c^\prime c)\cdot g\cdot g^\prime$ and in the second by $g\cdot g^\prime$.  This process thus can replace a labelling by an equivalent one in which all the non-trivial 2-dimensional `colour' is concentrated in one of the 2-simplices.   From the perspective of cellular formal maps, the obvious cellular encoding of the above combined map is
$$
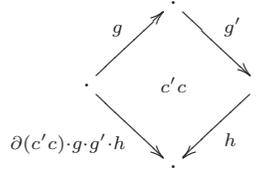
\xymatrix{&.\ar[dr]^{g^\prime}&\\
.\ar[ur]^{g}\ar@{}[rr] | {c^\prime c} \ar[dr]_{\partial(c^\prime c)\cdot g\cdot g^\prime\cdot h}&&.\ar[dl]^h\\
&.&
}$$

The geometric interpretation of this is that `integrating' around the boundary of the square picks up an element $c^\prime c$ of $C$ related to that boundary by $\partial(c^\prime c)$ being the evaluation of the labelling on the boundary path.  The analogy with integration around a curve seems important.  The integral is the non-commutative `sum' of the integrals over the parts of the subdivision, so, of course, the value of the labelling on this `horizontal' edge is of little importance.  It is added then subtracted.

Here the basic 2-cell stayed the same, but in other configurations,  it  may get conjugated or inverted or both as in the following variant. 

 We could equally well have had two different 2-simplices with  the same overall boundary.
$$
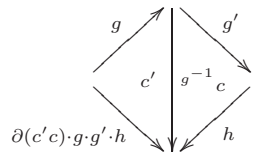
\xymatrix{&\hspace{1mm}\ar[dr]^{g^\prime}\ar[dd]^{{}^{g^{-1}}c}_{c^\prime\hspace{1mm}}&\\
\hspace{1mm}\ar[ur]^{g}\ar[dr]_{\partial(c^\prime c)\cdot g\cdot g^\prime\cdot h}&&\hspace{1mm}\ar[dl]^h\\
&\hspace{1mm}&
}$$
It is relatively simple to show, using the boundary and cocycle conditions and the extension schemes discussed in our earlier propositions, that this labelled triangulation and either of the previous ones are equivalent, likewise this second type of subdivision can be relabelled with $1_C$ in either of the two triangles with $c^\prime c$ in the other. 

The exercise is quite revealing of how thing behave in the composition process, but is left `to the reader'. We will instead examine a second example namely a subdivided annulus.

(ii) Consider an annulus with 1-skeleton labelled as follows:
$$
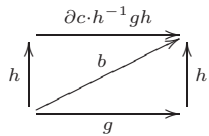
\xymatrix{\ar[rr]^{\partial c\cdot h^{-1}gh}&&\\
\ar[u]^h\ar[urr]^b\ar[rr]_g&&\ar[u]_h}$$and also an element $c\in C$, contributing to the label on the top edge.
The problem is thus to decide what $b$ is if the contribution of $c$ is concentrated in the top left triangle or in the bottom right. Likewise we could triangulate differently
$$
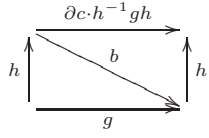
\xymatrix{\ar[rr]^{\partial c\cdot  h^{-1}gh}\ar[rrd]^b&&\\
\ar[u]^h\ar[rr]_g&&\ar[u]_h}$$
In the first case we can put $b = gh$, using an identity 2-cell in the bottom right and $~^hc^{-1}$ in the top left.  We could equally well set $b = h\cdot\partial c\cdot h^{-1}g\cdot h$ with the identity 2-label $1_P$ in the top left and $~^hc^{-1}\cdot g\cdot h$ in the bottom right.   We leave the other subdivision to the reader.  It is of note that here the use of `moves' as in many treatments in this area, is replaced by a geometric notion of equivalence, which is dominated by the cocycle condition.  This gives `for free' the independence of the end result on the order of combination of the local values, since any two such combinations will be the ends of an equivalence with the original in the middle!

\subsection{Cellular formal $\mathcal{C}$-maps} 
Combining simplices thus provides a simplification process which allows us to replace triangulated manifolds by manifolds with a given regular cellular decomposition.  These are much easier to handle.  We still will need base points in each 1-manifold and start vertices in each cell.  

Assume given a regular CW-complex $X$ having, for each cell, a specified `start 0-cell' among which are a set of distinguished base points.  Assume further that each cell has a specified orientation (so as to ensure that the boundary formulae make sense unambiguously).

\medskip

\textbf{Definition}

A \emph{cellular formal $\mathcal{C}$-map} $\lambda$ on $X$ consists of families of elements \\
(i) $\{c_f\}$ of $C$ indexed by the 2-cells, $f$, of $X$, and \\
(ii) $\{p_e\}$ of $P$ indexed by the 1-cells, $e$, of $X$ such that 

a) the boundary condition\\
\centerline{$\partial c_f$ = the ordered product of the edge labels of $f$}
is satisfied;\\
and

b) the cocycle condition is satisfied for each 3-cell. 

(In words b) gives, for each 3-cell $\sigma$, that the product of the labels on the boundary cells of $\sigma$ is trivial.)

\medskip

For a connected 1-manifold, $S$, decomposed as a CW-complex, (thus a subdivided circle),  there is no difference from the simplicial description we had before.  We have notions of formal $\mathcal{C}$-circuit given by a sequence of elements of $P$ and, more generally, if $S$ is not connected, we have a list of such formal $\mathcal{C}$-circuits.  

A \emph{cellular formal  $\mathcal{C}$-cobordism} between cellular formal $\mathcal{C}$-maps is the obvious thing. It is a cellular cobordism between the underlying 1-manifolds endowed with a formal $\mathcal{C}$-map that agrees with the two given $\mathcal{C}$-maps on the two ends of the cobordism. Here again the important ingredient is the cocycle condition and before going further we must say something more about both this and the boundary condition.

\medskip

The algebraic-combinatorial description of the cellular version formal $\mathcal{C}$-map is less explicitly  given above than for the simplicial version as a full description will require the introduction of some additional machinery, but this is not essential for the \emph{intuitive} development of the ideas.  We will, however, briefly sketch this extra theory in order to point the reader to sources which provide enough to construct a full development of the cellular theory. A more detailed treatment of this point will be given in \cite{TP:fhqft2}.

A few extra concepts are needed:
\begin{itemize}
\item \emph{Crossed complex}:  The basic idea is that of a chain complex of groups $(C_n,\partial)$, which are Abelian for $n\geq 3$, but with $C_2\to C_1$ being a crossed module.  The main example for us is the crossed complex of $X$, a CW-complex as above.  This has $C_n = \pi_n(X_n,X_{n-1},\mathbf{x})$, $X_n$, being as usual, the $n$-skeleton of $X$ and with $\partial$ the usual boundary map. Here we really need a many-object /groupoid version working with the multiple base points $\mathbf{x}$, but we will omit the detailed changes to the basic idea.  We write $\mathbf{\pi}(\mathbf{X})$ for this crossed complex. 
\item \emph{Free crossed module}:  The case of a 2-dimensional CW-complex $X$ is of some importance for our theory as the $B$-cobordisms will be surfaces and hence 2-dimensional regular CW-complexes once a decomposition is given. Any such 2-dimensional CW-complex yields a \emph{free} crossed module
$$\pi_2(X_2,X_1,X_0)\to \pi_1(X_1,X_0)$$
with $\pi_1(X_1,X_0)$, the fundamental groupoid of the 1-skeleton $X_1$ of $X$ based at the set of vertices $X_0$ of $X$.  Each 2-cell of $X$ gives a generating element in $\pi_2(X_2,X_1,X_0)$ and the assignment of the data for a cellular formal $\mathcal{C}$-map satisfying the boundary condition, is equivalent to specifying a morphism, $\lambda$, of crossed modules
$$\xymatrix{\pi_2(X_2,X_1,X_0)\ar[d]_{\lambda_2}\ar[r]^\partial &\pi_1(X_1,X_0)\ar[d]^{\lambda_1}\\
C\ar[r]_\partial & P}.$$
The boundary condition just states $\lambda_1 \partial = \partial \lambda_2$.
\item \emph{Free crossed complex}: The idea of free crossed complex is an extension of the above and $\pi(\mathbf{X})$ is free on the cells of $X$.  (In particular, $C_3 = \pi_3(X_3,X_2,\mathbf{x})$ is a collection of free $\pi_1(X)$-modules over the various basepoints. The generating set  is the set of 3-cells of $X$.) \end{itemize}

A formal $\mathcal{C}$-map, $\lambda$, is equivalent to a morphism of crossed complexes
$$\lambda : \pi(\mathbf{X})\to \mathcal{C},$$
or, expanding this, to 
$$\xymatrix{\ar[r] & \pi_3(X_3,X_2,X_0)\ar[d]_{\lambda_3}\ar[r]^\partial &\pi_2(X_2,X_1,X_0)\ar[d]_{\lambda_2}\ar[r]^\partial &\pi_1(X_1,X_0)\ar[d]_{\lambda_1}\\
\ar[r] & 1 \ar[r]_\partial &C\ar[r]_\partial & P}.$$
Each 3-cell gives an element of $\pi_3(X_3,X_2,X_0)$.  More exactly, if $\sigma$ is a 3-cell of $X$, then it can be specified by a characteristic map $\phi_\sigma : (B^3,S^2,\mathbf{s})\to (X_3,X_2,X_0)$ and thus we get an induced crossed complex morphism, which in the crucial dimensions gives
$$\xymatrix{\ar[r] & \pi_3(B^3,S^2,\mathbf{s})\ar[d]_{\phi_{\sigma,3}}\ar[r]^\partial &\pi_2(S^2,\phi_\sigma^{-1}(X_1),\mathbf{s})\ar[d]_{\phi_{\sigma,2}}\ar[r]^\partial &\pi_1(\phi_\sigma^{-1}(X_1),\mathbf{s})\ar[d]_{\phi_{\sigma,1}}\\
\ar[r] & \pi_3(X_3,X_2,X_0)\ar[d]_{\lambda_3}\ar[r]^\partial &\pi_2(X_2,X_1,X_0)\ar[d]_{\lambda_2}\ar[r]^\partial &\pi_1(X_1,X_0)\ar[d]_{\lambda_1}\\ 
\ar[r] & 1 \ar[r]_\partial &C\ar[r]_\partial & P}.$$
We have $\pi_3(B^3,S^2,\mathbf{s})$ is generated by the class of the 3-cell, $\langle e^3\rangle$ and $\phi_{\sigma,3}(\langle e^3\rangle) = \langle \sigma \rangle$.  The cocycle condition is then explicitly given by $\lambda_2\partial \langle\sigma\rangle = 1$.

The  explicit \emph{combinatorial} form of the cocycle condition for $\sigma$ will depend on the decomposition of the boundary $S^2$ given by   $\phi_\sigma^{-1}(X_1)$.  (This type of argument was first introduced in the original paper by J. H. C. Whitehead, \cite{JHCW:CH2}.  It can also be found in the forthcoming book by Brown, Higgins and Sivera, \cite{RBRS:2005},  work by Brown and Higgins, \cite{B&H1987,B&H1991} and by Baues, \cite{HJB:AH,HJB:4D}, where, however, crossed complexes are called \emph{crossed chain complexes}.) Our use of this cocycle condition does not require such a detailed description so we will not attempt to give one here.

	The next ingredient is to cellularise `equivalence'.  We can do this for arbitrary formal $\mathcal{C}$-maps specialising to 1- or 2-dimensions (cobordisms) afterwards. We use a regular cellular decomposition of the space $X\times I$, with possibly different regular CW-complex decompositions on the two ends, but with the base points `fixed' so that $\mathbf{x}\times I$ is a subcomplex of $X\times I$. 

\

\textbf{Definition}
 
Given cellular formal $\mathcal{C}$-maps $\lambda_i$ on $X_i = X\times \{i\}$, for $i = 0,1$, they will be \emph{equivalent} if there is a cellular formal $\mathcal{C}$-map $\Lambda$ on a cellular decomposition of $X\times I$ extending $\lambda_0$ and $\lambda_1$ and assigning $1_p$ to each edge in $\mathbf{x}\times I$.  

Again equivalence is an equivalence relation.  It allows the combination and collection processes examined in the previous subsection to be made precise.  In other words:\emph{
\begin{itemize}
\item if we triangulate each cell of a CW-complex $X$ in such a way that the result gives a triangulation $K$ of the space, then a formal $\mathcal{C}$-map, $\lambda$, on $K$ determines a cellular formal $\mathcal{C}$-map on $X$;
\item equivalent simplicial formal $\mathcal{C}$-maps on (possibly different) such triangulations  yield equivalent formal $\mathcal{C}$-maps on $X$;
\item given any cellular formal $\mathcal{C}$-map, $\mu$, on $X$ and a triangulation, $K$, of $X$ subdividing the cells of $X$, there is a simplicial formal $\mathcal{C}$-map on $K$ that combines to give $\mu$.
\item Any two different ways of combining a formal $\mathcal{C}$-map into a cellular one on $X$ will be equivalent.  (In other words, the order of combination and the choices made make no difference \emph{up to equivalence}.)
\end{itemize}}
\textbf{Remarks:}

(i) Full proofs of these would use cellular and simplicial decompositions of $X\times I$, but would also need the introduction of far more of the theory of crossed modules, crossed complexes and their classifying spaces than we have available here. Because of that, the proofs are omitted here in order to make this introduction to formal $\mathcal{C}$-maps easier to approach.

(ii) Any simplicial formal $\mathcal{C}$-map on $K$ is, of course, a cellular one for the obvious regular CW-structure on $|K|$.

\medskip

The notion of equivalent cellular formal $\mathcal{C}$-cobordisms can now be formulated. Given the obvious set-up with $\mathbf{F}$ and $\mathbf{G}$, two such cobordisms between $\mathbf{g_1}$ and $\mathbf{g_2}$, they will be equivalent if they are equivalent as formal $\mathcal{C}$-maps by an equivalence that is constant on the two `ends'.

\subsection{2-dimensional formal $\mathcal{C}$-maps}\label{building blocks}

It is now easy to describe a set of `building blocks' for all cellular formal $\mathcal{C}$-maps on orientable surfaces and thus all cobordisms between 1-dimensional formal $\mathcal{C}$-maps.  Again we want to emphasise the fact that  these models provide \emph{formal combinatorial models} for the characteristic maps with target a 2-type.

We will shortly introduce the formal version of 1+1 HQFTs with a `background' crossed module, $\mathcal{C}$, which is a model for a 2-type $B$, represented by that crossed module.  As the basic manifolds  are 1-dimensional, they are just disjoint unions of pointed oriented circles, and  so a formal $\mathcal{C}$-map on a 1-manifold, as we saw earlier  (page \pageref{1Cmap}), is specified by a list of lists of  elements in $P$, one list for each connected component. Cellularly we can assume that the lists have just one element in them, obtained from the simplicial case by multiplying the elements in the list together in order. The corresponding cellular cobordisms are then compact oriented surfaces $W$ with pointed oriented boundary endowed with a formal $\mathcal{C}$-map $\Lambda$ as above. Since such surfaces can be built up from three basic models, the disc, annulus and disc with two holes (pair of trousers), we need only examine what formal $\mathcal{C}$-maps  look like on these basic example spaces and how they compose and combine, as any formal 1 + 1 `$\mathcal{C}$-HQFT' will be determined completely by its behaviour on the formal maps on these basic surfaces.

\vspace{5mm}

\textbf{Formal $\mathcal{C}$-Discs.}

The only formal $\mathcal{C}$-maps that makes sense on the disc must have an element $c\in C$ assigned to the interior 2-cell with the boundary $\partial c$ assigned to the single 1-cell, i.e.

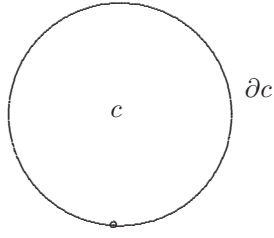
\begin{center}
\font\thinlinefont=cmr5
\begingroup\makeatletter\ifx\SetFigFont\undefined
\def\x#1#2#3#4#5#6#7\relax{\def\x{#1#2#3#4#5#6}}%
\expandafter\x\fmtname xxxxxx\relax \def\y{splain}%
\ifx\x\y   
\gdef\SetFigFont#1#2#3{%
  \ifnum #1<17\tiny\else \ifnum #1<20\small\else
  \ifnum #1<24\normalsize\else \ifnum #1<29\large\else
  \ifnum #1<34\Large\else \ifnum #1<41\LARGE\else
     \huge\fi\fi\fi\fi\fi\fi
  \csname #3\endcsname}%
\else
\gdef\SetFigFont#1#2#3{\begingroup
  \count@#1\relax \ifnum 25<\count@\count@25\fi
  \def\x{\endgroup\@setsize\SetFigFont{#2pt}}%
  \expandafter\x
    \csname \romannumeral\the\count@ pt\expandafter\endcsname
    \csname @\romannumeral\the\count@ pt\endcsname
  \csname #3\endcsname}%
\fi
\fi\endgroup
\mbox{\beginpicture
\setcoordinatesystem units <0.50000cm,0.50000cm>
\unitlength=0.50000cm
\linethickness=1pt
\setplotsymbol ({\makebox(0,0)[l]{\tencirc\symbol{'160}}})
\setshadesymbol ({\thinlinefont .})
\setlinear
%
%
\linethickness= 0.500pt
\setplotsymbol ({\thinlinefont .})
\ellipticalarc axes ratio  0.078:0.078  360 degrees 
	from  5.793 18.944 center at  5.715 18.944
%
%
\put{$\partial c$} [lB] at  9.207 22.331
%
%
\put{$c$} [lB] at  5.637 21.827
%
%
\linethickness= 0.500pt
\setplotsymbol ({\thinlinefont .})
\ellipticalarc axes ratio  2.963:2.963  360 degrees 
	from  8.860 21.878 center at  5.897 21.878
\linethickness=0pt
\putrectangle corners at  2.917 24.856 and  9.207 18.849
\endpicture}
$$Disc(c) : \emptyset \to \mathbf{\partial c}$$
\end{center}
(Remember that here $\empty$ is the notation for the empty 1-manifold with the empty map as characteristic map.)
Later we will see that these give the crucial difference between the formal $\mathcal{C}$-theory and the standard form of \cite{turaev:hqft1}.

\vspace{5mm}

\textbf{Formal $\mathcal{C}$-Annuli.}

Let $Cyl$ denote the cylinder/annulus, $S^1 \times [0,1]$.  We fix an orientation of $Cyl$ once and for all, and set $Cyl^0 = S^1\times (0) \subset \partial Cyl$ and $Cyl^1 = S^1\times (1) \subset \partial Cyl$.  We provide $Cyl^0$ and $Cyl^1$ with base points $z^0 = (s,0)$, $z^1=(s,1)$, respectively, where $s \in S^1$.  As in \cite{turaev:hqft1}, let $\varepsilon, \mu =\pm $, and denote by $Cyl_{\varepsilon,\mu}$ the triple $(Cyl,Cyl^0_\varepsilon,Cyl^1_\mu)$. This is an annulus with oriented pointed boundary,
$$\partial Cyl_{\varepsilon,\mu} = (\varepsilon Cyl^0_\varepsilon)\cup(\mu Cyl^1_\mu),$$
where by $-X$ we mean $X$ with opposite orientation.
A formal $\mathcal{C}$-map, $\Lambda$, on $Cyl$ may be drawn diagrammatically as:
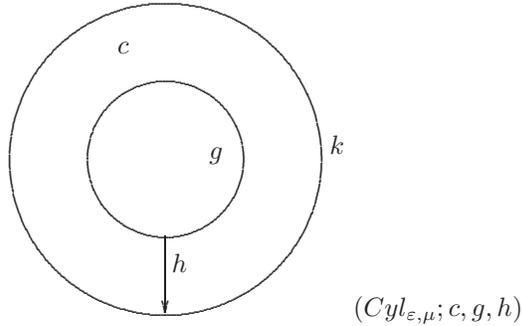
\begin{center}\hspace*{2cm}
\font\thinlinefont=cmr5
\begingroup\makeatletter\ifx\SetFigFont\undefined
\def\x#1#2#3#4#5#6#7\relax{\def\x{#1#2#3#4#5#6}}%
\expandafter\x\fmtname xxxxxx\relax \def\y{splain}%
\ifx\x\y   
\gdef\SetFigFont#1#2#3{%
  \ifnum #1<17\tiny\else \ifnum #1<20\small\else
  \ifnum #1<24\normalsize\else \ifnum #1<29\large\else
  \ifnum #1<34\Large\else \ifnum #1<41\LARGE\else
     \huge\fi\fi\fi\fi\fi\fi
  \csname #3\endcsname}%
\else
\gdef\SetFigFont#1#2#3{\begingroup
  \count@#1\relax \ifnum 25<\count@\count@25\fi
  \def\x{\endgroup\@setsize\SetFigFont{#2pt}}%
  \expandafter\x
    \csname \romannumeral\the\count@ pt\expandafter\endcsname
    \csname @\romannumeral\the\count@ pt\endcsname
  \csname #3\endcsname}%
\fi
\fi\endgroup
\mbox{\beginpicture
\setcoordinatesystem units <0.70000cm,0.70000cm>
\unitlength=0.50000cm
\linethickness=1pt
\setplotsymbol ({\makebox(0,0)[l]{\tencirc\symbol{'160}}})
\setshadesymbol ({\thinlinefont .})
\setlinear
%
%
\linethickness= 0.500pt
\setplotsymbol ({\thinlinefont .})
\ellipticalarc axes ratio  2.963:2.963  360 degrees 
	from  8.970 21.643 center at  6.007 21.643
%
%
\linethickness= 0.500pt
\setplotsymbol ({\thinlinefont .})
\putrule from  5.997 20.180 to  5.997 18.735
%
%
\plot  5.944 18.946  5.997 18.735  6.049 18.946 /
%
%
%
\put{$c$} [lB] at  5.093 23.664
%
%
\put{$g$} [lB] at  6.835 21.660
%
%
\put{$h$} [lB] at  6.128 19.490
%
%
\put{$k$} [lB] at  9.129 21.749
%
%
\linethickness= 0.500pt
\setplotsymbol ({\thinlinefont .})
\ellipticalarc axes ratio  1.482:1.482  360 degrees 
	from  7.489 21.643 center at  6.007 21.643
\linethickness=0pt
\putrectangle corners at  3.027 24.621 and  9.129 18.665
\endpicture} $(Cyl_{\varepsilon,\mu};c,g,h)$\label{CYL}
\end{center}
with initial vertex, $s$, for the 2-cell at the head of $h$, i.e. on the outer circle.
This diagram will represent the cobordism that we will denote $(Cyl_{\varepsilon,\mu};c,g,h)$.  Similar notation may be used in other contexts without further comment. 

We omit the orientations on the boundary circles so as to avoid the need to repeat more or less the same diagram several times.  The exact expression for $k$ will change depending on the orientations and which vertex is used as the `start' of the 2-cell.  Reading off clockwise $\partial c = kh^{-1}g^{-1}h$, so $k = \partial c.h^{-1}gh$ if we assume both boundaries are clockwise oriented.  If we change the start vertex to the inner circle we need to act on $c$ with $h$ to keep the same element $k$ on the outer circle.  With the same labelling on the edges, that change of base point changes $c$ to ${}^hc$.

The loop, $\Lambda|_{Cyl^1_\mu}$, i.e. $k$, represents $(\partial c\cdot h^{-1}gh)$ or its inverse depending on the sign of $\mu$. There are two special cases that generate all the others: (i) $c = 1$, which corresponds to the case already handled in \cite{turaev:hqft1}, and (ii) $h=1$, where the base point $s$ does not move during the cobordism.  The general case, illustrated in the figure, is the composite of particular instances of the two cases.
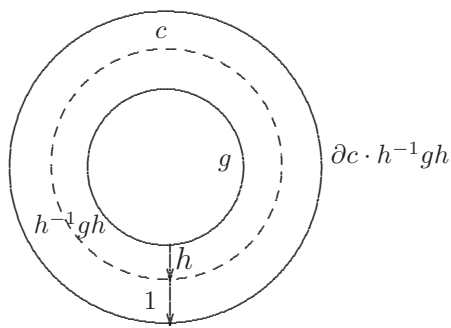
\begin{center}
\font\thinlinefont=cmr5
\begingroup\makeatletter\ifx\SetFigFont\undefined
\def\x#1#2#3#4#5#6#7\relax{\def\x{#1#2#3#4#5#6}}%
\expandafter\x\fmtname xxxxxx\relax \def\y{splain}%
\ifx\x\y   
\gdef\SetFigFont#1#2#3{%
  \ifnum #1<17\tiny\else \ifnum #1<20\small\else
  \ifnum #1<24\normalsize\else \ifnum #1<29\large\else
  \ifnum #1<34\Large\else \ifnum #1<41\LARGE\else
     \huge\fi\fi\fi\fi\fi\fi
  \csname #3\endcsname}%
\else
\gdef\SetFigFont#1#2#3{\begingroup
  \count@#1\relax \ifnum 25<\count@\count@25\fi
  \def\x{\endgroup\@setsize\SetFigFont{#2pt}}%
  \expandafter\x
    \csname \romannumeral\the\count@ pt\expandafter\endcsname
    \csname @\romannumeral\the\count@ pt\endcsname
  \csname #3\endcsname}%
\fi
\fi\endgroup
\mbox{\beginpicture
\setcoordinatesystem units <0.70000cm,0.70000cm>
\unitlength=1.00000cm
\linethickness=1pt
\setplotsymbol ({\makebox(0,0)[l]{\tencirc\symbol{'160}}})
\setshadesymbol ({\thinlinefont .})
\setlinear
%
%
\put{\SetFigFont{12}{14.4}{rm}$h$} [lB] at  6.191 19.717
%
%
\linethickness= 0.500pt
\setplotsymbol ({\thinlinefont .})
\ellipticalarc axes ratio  2.963:2.963  360 degrees 
	from  8.970 21.643 center at  6.007 21.643
%
%
\linethickness= 0.500pt
\setplotsymbol ({\thinlinefont .})
\setdashes < 0.1270cm>
\ellipticalarc axes ratio  2.191:2.191  360 degrees 
	from  8.217 21.700 center at  6.026 21.700
%
%
\linethickness= 0.500pt
\setplotsymbol ({\thinlinefont .})
\setsolid
\plot  6.094 20.18  6.094 18.639 /
%
%
%
\plot  6.042 19.7  6.094 19.5  6.148 19.7 /

\plot  6.042 18.851  6.094 18.639  6.148 18.851 /

%
%
\put{ $g$} [lB] at  6.835 21.660
%
%
\linethickness= 0.500pt
\setplotsymbol ({\thinlinefont .})
\ellipticalarc axes ratio  1.482:1.482  360 degrees 
	from  7.489 21.643 center at  6.007 21.643
%
%
\put{ $1$} [lB] at  5.429 18.955
%
%
\put{$\partial c\cdot h^{-1}gh$} [lB] at  9.144 21.717
%
%
\put{$h^{-1}gh$} [lB] at  3.5 20.384
%
%
\put{$c$} [lB] at  5.810 24.098
\linethickness=0pt
\putrectangle corners at  3.027 24.621 and  9.144 18.614
\endpicture}
\end{center}
A remark should be made here about the combination of cobordisms, although we will handle this in some more detail later.  The rule is more or less the obvious one.  In fact it is always possible to triangulate the cellular map and to combine the composing cobordisms followed by recombination to get a cellular map.  The result does not depend on the  triangulation, again \emph{up to equivalence}.  Of course two different choices of start vertex for the 2-cell of the combined cobordism, will give different labellings, but this is easily rectified if it occurs.  The use of globular 2-cell notation, labelling start and end vertices of each 2-cell (in a 2-categorical fashion) can help here as it combines a label, $c$, from $C$ on the 2-cell with the initial 1-cell from the start to the finish vertices $g\in P$, say, to get $(c,g)$, an element in the (group) semidirect product, $C\rtimes P$. Some of the combinations of labellings that arise in calculations then correspond in part to the semidirect product formula for multiplication (see later).

\medskip

\textbf{Formal $\mathcal{C}$-Disc with 2 holes}

Let $D$ be an oriented 2-disc with two holes.  We will denote the boundary components of $D$ for convenience by $Y$, $Z$, and $T$ and provide them with base points $y$, $z$ and $t$ respectively.  For any choice of signs $\varepsilon$, $\mu$,  $\nu = \pm$, we denote by $D_{\varepsilon, \mu,\nu}$ the tuple $(D,Y_\varepsilon,Z_\mu,T_\nu)$.  This is a 2-disc with two holes with oriented pointed boundary. By definition,
$$\partial D_{\varepsilon, \mu,\nu} = (\varepsilon Y_\varepsilon)\cup (\mu Z_\mu) \cup (\nu T_\nu).$$
Finally we fix two proper embedded arcs $yt$ and $zt$ in $D$ leading from $y$ and $z$ to $t$.  A formal $\mathcal{C}$-map $\lambda$ on $D_{\varepsilon, \mu,\nu}$  will, in general, assign elements of $P$ to each boundary component and to each arc.  As for the annulus we may assume that the formal map assigns $1_P$ to both $yt$ and $zt$, as the general case can be generated by this one together with cylinders. In addition the single 2-cell will be assigned an element  $c$ of $C$.  (As usual a start vertex for each 2-cell is used - but is not always made explicit.)

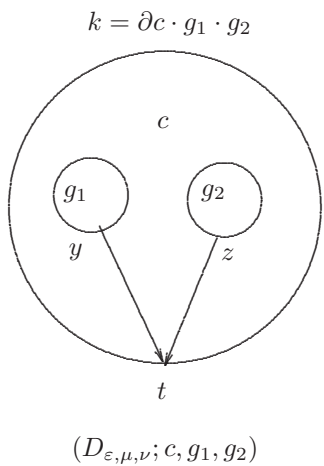
\begin{center}

\font\thinlinefont=cmr5
\begingroup\makeatletter\ifx\SetFigFont\undefined
\def\x#1#2#3#4#5#6#7\relax{\def\x{#1#2#3#4#5#6}}%
\expandafter\x\fmtname xxxxxx\relax \def\y{splain}%
\ifx\x\y   
\gdef\SetFigFont#1#2#3{%
  \ifnum #1<17\tiny\else \ifnum #1<20\small\else
  \ifnum #1<24\normalsize\else \ifnum #1<29\large\else
  \ifnum #1<34\Large\else \ifnum #1<41\LARGE\else
     \huge\fi\fi\fi\fi\fi\fi
  \csname #3\endcsname}%
\else
\gdef\SetFigFont#1#2#3{\begingroup
  \count@#1\relax \ifnum 25<\count@\count@25\fi
  \def\x{\endgroup\@setsize\SetFigFont{#2pt}}%
  \expandafter\x
    \csname \romannumeral\the\count@ pt\expandafter\endcsname
    \csname @\romannumeral\the\count@ pt\endcsname
  \csname #3\endcsname}%
\fi
\fi\endgroup
\mbox{\beginpicture
\setcoordinatesystem units <0.70000cm,0.70000cm>
\unitlength=0.50000cm
\linethickness=1pt
\setplotsymbol ({\makebox(0,0)[l]{\tencirc\symbol{'160}}})
\setshadesymbol ({\thinlinefont .})
\setlinear
%
%
\put{$y$} [lB] at  4.843 20.373
%
%
\linethickness= 0.500pt
\setplotsymbol ({\thinlinefont .})
\ellipticalarc axes ratio  0.705:0.705  360 degrees 
	from  5.971 21.522 center at  5.266 21.522
%
%
\linethickness= 0.500pt
\setplotsymbol ({\thinlinefont .})
\ellipticalarc axes ratio  0.705:0.705  360 degrees 
	from  8.507 21.431 center at  7.802 21.431
%
%
\linethickness= 0.500pt
\setplotsymbol ({\thinlinefont .})
\plot  7.654 20.701  6.674 18.294 /
%
%
\plot  6.705 18.510  6.674 18.294  6.803 18.470 /
%
%
%
\linethickness= 0.500pt
\setplotsymbol ({\thinlinefont .})
\plot  5.427 20.925  6.625 18.371 /
%
%
\plot  6.487 18.540  6.625 18.371  6.583 18.585 /
%
%
%
\put{$c$} [lB] at  6.519 22.807
%
%
\linethickness= 0.500pt
\setplotsymbol ({\thinlinefont .})
\ellipticalarc axes ratio  2.963:2.963  360 degrees 
	from  9.631 21.289 center at  6.668 21.289
%
%
\put{$z$} [lB] at  7.732 20.288
%
%
\put{$ k = \partial c \cdot  g_1\cdot g_2$} [lB] at  5.2 24.661
%
%
\put{$t$} [lB] at  6.519 17.645
%
%
\put{$g_1$} [lB] at  4.737 21.484
%
\put{$g_2$} [lB] at  7.355 21.510
\linethickness=0pt
\putrectangle corners at  3.687 24.661 and  9.648 17.645
\endpicture}
$$(D_{\varepsilon, \mu,\nu};c, g_1,g_2)$$
\end{center}
(To assist the reader in the deciphering of these pictures here are some points.  The 2-cell is oriented clockwise as are the boundary components.  The start vertex is $t$.  Draw the corresponding surface polygon wth label $k$ on the outer large circle.  Read off: $\partial c = kg_2^{-1}g_1^{-1}$ giving $k = \partial c.g_1.g_2$ as claimed.  If you prefer an anticlockwise orientation - look in the mirror!)

This situation leads to an interesting relation. If we have a formal $\mathcal{C}$-map on  $D_{\varepsilon, \mu,\nu}$ in which, for simplicity, we assume that the 2-cell is assigned the element $1_C$ and then add suitable cylinders, labelled with $c_1$ and $c_2$ respectively, to the boundary components $Y$ and $Z$ then the resulting cobordism  can be rearranged to give a labelling with the 2-cell coloured $c_1\cdot {}^{g_1}c_2$ as shown in the following diagram:

\

\begin{center}
\font\thinlinefont=cmr5
\begingroup\makeatletter\ifx\SetFigFont\undefined
\def\x#1#2#3#4#5#6#7\relax{\def\x{#1#2#3#4#5#6}}%
\expandafter\x\fmtname xxxxxx\relax \def\y{splain}%
\ifx\x\y   
\gdef\SetFigFont#1#2#3{%
  \ifnum #1<17\tiny\else \ifnum #1<20\small\else
  \ifnum #1<24\normalsize\else \ifnum #1<29\large\else
  \ifnum #1<34\Large\else \ifnum #1<41\LARGE\else
     \huge\fi\fi\fi\fi\fi\fi
  \csname #3\endcsname}%
\else
\gdef\SetFigFont#1#2#3{\begingroup
  \count@#1\relax \ifnum 25<\count@\count@25\fi
  \def\x{\endgroup\@setsize\SetFigFont{#2pt}}%
  \expandafter\x
    \csname \romannumeral\the\count@ pt\expandafter\endcsname
    \csname @\romannumeral\the\count@ pt\endcsname
  \csname #3\endcsname}%
\fi
\fi\endgroup
\mbox{\beginpicture
\setcoordinatesystem units <0.70000cm,0.70000cm>
\unitlength=0.50000cm
\linethickness=1pt
\setplotsymbol ({\makebox(0,0)[l]{\tencirc\symbol{'160}}})
\setshadesymbol ({\thinlinefont .})
\setlinear
%
%
\linethickness= 0.500pt
\setplotsymbol ({\thinlinefont .})
\ellipticalarc axes ratio  0.705:0.705  360 degrees 
	from  5.442 22.092 center at  4.737 22.092
%
%
\linethickness= 0.500pt
\setplotsymbol ({\thinlinefont .})
\ellipticalarc axes ratio  0.705:0.705  360 degrees 
	from 14.783 22.077 center at 14.078 22.077
%
%
\linethickness= 0.500pt
\setplotsymbol ({\thinlinefont .})
\ellipticalarc axes ratio  0.705:0.705  360 degrees 
	from 17.319 21.986 center at 16.614 21.986
%
%
\linethickness= 0.500pt
\setplotsymbol ({\thinlinefont .})
\plot 16.463 21.256 15.483 18.849 /
%
%
\plot 15.514 19.065 15.483 18.849 15.612 19.025 /
%
%
%
\linethickness= 0.500pt
\setplotsymbol ({\thinlinefont .})
\plot 14.296 21.359 15.435 18.925 /
%
%
\plot 15.297 19.094 15.435 18.925 15.393 19.139 /
%
%
%
\put{$g_1$} [lB] at 13.661 22.049
%
%
\put{$g_2$} [lB] at 16.228 21.971
%
%
\put{$c_1\cdot {}^{g_1}c_2$} [lB] at 14.5 23.4
%
%
\linethickness= 0.500pt
\setplotsymbol ({\thinlinefont .})
\ellipticalarc axes ratio  0.028:0.028  360 degrees 
	from  3.067 21.554 center at  3.040 21.554
%
%
\linethickness= 0.500pt
\setplotsymbol ({\thinlinefont .})
\ellipticalarc axes ratio  2.963:2.963  360 degrees 
	from  8.956 21.882 center at  5.992 21.882
%
%
\linethickness= 0.500pt
\setplotsymbol ({\thinlinefont .})
\ellipticalarc axes ratio  0.705:0.705  360 degrees 
	from  7.811 22.121 center at  7.106 22.121
%
%
\linethickness= 0.500pt
\setplotsymbol ({\thinlinefont .})
\ellipticalarc axes ratio  2.963:2.963  360 degrees 
	from 18.443 21.844 center at 15.479 21.844
%
%
\linethickness= 0.500pt
\setplotsymbol ({\thinlinefont .})
\setdashes < 0.1270cm>
\ellipticalarc axes ratio  1.168:1.168  360 degrees 
	from  5.906 22.092 center at  4.737 22.092
%
%
\linethickness= 0.500pt
\setplotsymbol ({\thinlinefont .})
\ellipticalarc axes ratio  1.137:1.137  360 degrees 
	from  8.202 22.145 center at  7.065 22.145
%
%
\linethickness= 0.500pt
\setplotsymbol ({\thinlinefont .})
\setsolid
\plot  4.815 21.378  5.954 18.944 /
%
%
\plot  5.817 19.113  5.954 18.944  5.912 19.158 /
%
%
%
\linethickness= 0.500pt
\setplotsymbol ({\thinlinefont .})
\plot  7.065 21.404  6.085 18.997 /
%
%
\plot  6.116 19.213  6.085 18.997  6.214 19.173 /
%
%
%
\put{$c_1$} [lB] at  4.525 23.044
%
%
\put{$c_2$} [lB] at  6.879 23.019
%
%
\put{$g_1$} [lB] at  4.498 22.092
%
%
\put{$g_2$} [lB] at  6.720 22.172
%
%
\put{$\sim$} [lB] at 10.4 22.013
%
%
\put{$\partial c_1\cdot g_1\cdot \partial c_2\cdot g_2$} [lB] at  4.2 25.5
%
%
\put{$\partial c_1\cdot g_1\cdot \partial c_2\cdot g_2$} [lB] at 13.6 25.5
\linethickness=0pt
\putrectangle corners at  2.995 24.977 and 18.459 18.824
\label{pantsreduction}
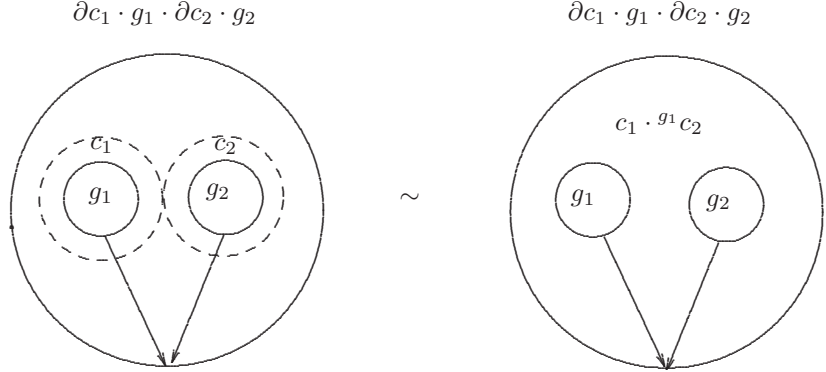
\endpicture}\\ Figure \ref{pantsreduction} : From the case $c=1$ to the general one.
\end{center}
The importance of this element $c_1\cdot {}^{g_1}c_2$ is that it is the $C$-part of the product of the two cylinder labels in the semidirect product, $C\rtimes P$, more exactly, the elements $(c_1,g_1)$ and $(c_2,g_2) \in C\rtimes P$ correspond to the two added cylinders and within that semi-direct product $(c_1,g_1)\cdot (c_2,g_2) = (c_1\cdot {}^{g_1}c_2, g_1g_2)$.

\section{Formal HQFTs}
As before we will restrict attention to modelling 1+1 HQFTs and so, here,
 will give a definition of a formal HQFT only for that case. First some notation and a convention:

If we have formal $\mathcal{C}$-cobordisms,
$$\mathbf{F} : \mathbf{g}_0\to \mathbf{g}_1, \quad \mathbf{G} : \mathbf{g}_1\to \mathbf{g}_2$$ 
then we will denote the composite $\mathcal{C}$-cobordism by $\mathbf{F}\#_{\mathbf{g}_1}\mathbf{G}$. 

For $\mathbf{g}$ as before, the trivial identity $\mathcal{C}$-cobordism on $\mathbf{g}$ will be denoted $1_{\mathbf{g}}$.

Finally, unless otherwise stated we will assume that all vector spaces, projective modules etc. will be of finite type.
\subsection{The definition \label{deffhqft}}
Fix, as before, a crossed module, $\mathcal{C} = (C,P, \partial)$, and also  fix a ground field, $\mathbb{K}$.

A \emph{formal HQFT} with background $\mathcal{C}$ assigns\begin{itemize}
\item to each formal $\mathcal{C}$-circuit, $\mathbf{g} = (g_1,\ldots, g_n)$,   a $\mathbb{K}$-vector space $\tau(\mathbf{g})$, and by extension, to each formal $\mathcal{C}$-map on a 1-manifold $S$, given by a list, $\mathbf{g} = \{\mathbf{g}_i\, | \, i = 1,2, \ldots, m\}$ of formal $\mathcal{C}$-circuits, a vector space  $\tau(\mathbf{g})$ and an identification, $$\tau(\mathbf{g}) = \bigotimes_{i = 1, \ldots, m}\tau(\mathbf{g}_i), $$giving  $\tau(\mathbf{g})$ as a tensor product;
\item to any formal $\mathcal{C}$-cobordism, $(M,\mathbf{F})$ between $(S_0,\mathbf{g}_0)$ and $(S_1,\mathbf{g}_1)$ , a $\mathbb{K}$-linear transformation
$$\tau(\mathbf{F}) : \tau(\mathbf{g}_0) \to \tau(\mathbf{g}_1).$$
\end{itemize} 
These assignments are to satisfy the following axioms:
\begin{enumerate}[(i)]
\item Disjoint union of formal $\mathcal{C}$-maps corresponds to tensor product of the corresponding vector spaces via specified isomorphisms:
$$\tau(\mathbf{g}\sqcup \mathbf{h}) \stackrel{\cong}{\to}\tau(\mathbf{g})\otimes \tau(\mathbf{h}),$$
$$\tau(\emptyset)\stackrel{\cong}{\to}\mathbb{K}$$
for the ground field $\mathbb{K}$, so that a) the diagram of specified isomorphisms
$$\xymatrix{\tau(\mathbf{g})\ar[r]^\cong\ar[d]_\cong &\tau(\mathbf{g}\sqcup \emptyset)\ar[d]^\cong \\
\tau(\mathbf{g})\otimes \mathbb{K}&\tau(\mathbf{g})\otimes \tau(\emptyset)\ar[l]^\cong}$$
for $\mathbf{g} \to \mathbf{g}\sqcup\emptyset$, commutes and similarly for $\mathbf{g} \to \emptyset\sqcup \mathbf{g}$, and b) the assignments are compatible with the associativity isomorphisms for $\sqcup$ and $\otimes$, (so that $\tau$ satisfies the usual axioms for a symmetric monoidal functor).
 \item For formal $\mathcal{C}$-cobordisms
 $$\mathbf{F} : \mathbf{g}_0\to \mathbf{g}_1, \quad \mathbf{G} : \mathbf{g}_1\to \mathbf{g}_2$$ with composite $\mathbf{F}\#_{\mathbf{g}_1}\mathbf{G}$, we have 
 $$\tau(\mathbf{F}\#_{\mathbf{g}_1}\mathbf{G}) = \tau(\mathbf{G})\tau(\mathbf{F}) : \tau(\mathbf{g}_0) \to  \tau(\mathbf{g}_2).$$
 \item For the identity formal $\mathcal{C}$-cobordism on  $\mathbf{g}$, 
 $$\tau(1_\mathbf{g}) = 1_{\tau(\mathbf{g})}.$$
 \item Interaction of cobordisms and disjoint union is transformed correctly by $\tau$, i.e., given formal $\mathcal{C}$-cobordisms
 $$\mathbf{F} : \mathbf{g}_0\to \mathbf{g}_1, \quad \mathbf{G} : \mathbf{h}_0\to \mathbf{h}_1$$ 
 the following diagram
$$ \xymatrix{\tau(\mathbf{g}_0\sqcup\mathbf{h}_0)\ar[r]^\cong\ar[d]_ {\tau(\mathbf{F}\sqcup\mathbf{G})}&\tau(\mathbf{g}_0)\otimes\tau(\mathbf{h}_0)\ar[d]^{\tau(\mathbf{F})\otimes\tau(\mathbf{G})}\\
 \tau(\mathbf{g}_1\sqcup \mathbf{h}_1)\ar[r]_\cong&\tau(\mathbf{g}_1)\otimes \tau(\mathbf{h}_1)}$$
 commutes, compatibly with the associativity structure.
\end{enumerate}
\subsection{Basic Structure}
We know that formal $\mathcal{C}$-maps could be specified by composing / combining the basic building blocks outlined in section \ref{building blocks}. As a formal 1+1 HQFT transforms the formal $\mathcal{C}$-maps to vector space structure compatibly with the combination rules of gluing and disjoint union, to specify a formal HQFT, we need only give it on the connected 1-manifolds (formal $\mathcal{C}$-circuits) and  on the building blocks mentioned before, and we can limit our specification to the cellular examples.  We assume $\tau$ is a formal HQFT with $\mathcal{C}$ as before.\footnote{Throughout this section it may help to refer to the corresponding discussion in the original paper, \cite{turaev:hqft1}.}

On a formal $\mathcal{C}$-circuit $\mathbf{g} = (g_1,\ldots, g_n)$, we can assume $n=1$, since the obvious formal $\mathcal{C}$-cobordism between $\mathbf{g}$ and $\{(g_1\ldots g_n)\}$, based on the cylinder yields an isomorphism 
$$\tau(\mathbf{g})\stackrel{\cong}{\to}\tau(g_1\ldots g_n).$$
For any element $g\in P$, we thus have the formal $\mathcal{C}$-circuit $\{(g)\}$ and a vector space $\tau(g)$, where we have shortened the notation in an obvious way. In fact, for later use it will be convenient to change notation to $L_g$ (or for the more general case $L_\mathbf{g}$) as otherwise we will end with far too many brackets!

For a general $\mathbf{g} = (\mathbf{g}_1,\ldots, \mathbf{g}_n)$, we now have 
\begin{align*}L_\mathbf{g} = \bigotimes_{i = 1}^n L_{\mathbf{g}_i}\end{align*}
The special case when $\mathbf{g}$ is empty gives $L_\emptyset = \mathbb{K}$ and the isomorphism in section \ref{deffhqft}, and above, are compatible with these assignments.

The basic formal $\mathcal{C}$-cobordisms give us various structural maps:
\begin{itemize}
\item the formal $\mathcal{C}$-disc with $c\in C$ gives 
$$\tau(Disc(c)) : \tau(\emptyset)\to \tau(\partial c),$$
that is, a linear map, which we will write as 
$$\ell_c : \mathbb{K} \to L_{\partial c}.$$
(The values of these linear maps on the 1 of the field, $\mathbb{K}$, will play a crucial role in the classification of these HQFTs; we write $\tilde{c} := \ell_c(1) \in L_{\partial c}.$)
\item the formal $\mathcal{C}$-annuli of the two basic types yield\\
(a) $(Cyl_{\varepsilon,\mu}; 1,g,h) : \{(g)\}\to \{(h^{-1}gh)\},$ and hence a linear isomorphism
$$L_g \to L_{h^{-1}gh}$$
(cf. \cite{turaev:hqft1}), or a related one, depending on the sign of $\mu$;\\
or\\
(b) $(Cyl_{\varepsilon,\mu}; c,g,1) : \{(g)\}\to\{ (\partial c\cdot g)\}$ and a linear isomorphism,
$$L_g\to L_{\partial c\cdot g},$$
again with variants for other signs.
\item the formal $\mathcal{C}$-disc with 2 holes,
$$(D_{\varepsilon,\mu,\nu};c,g_1,g_2) : \{(g_1),(g_2)\} \to \{(\partial c\cdot g_1\cdot g_2)\},$$
giving a bilinear map $$L_{g_1}\otimes L_{g_2}\to L_{\partial c\cdot g_1\cdot g_2}.$$
Again, the key case is $c = 1$ and consequently,
$$L_{g_1}\otimes L_{g_2} \to L_{g_1g_2}.$$
The general case can be obtained from that and a suitable formal $\mathcal{C}$-annulus,
$$L_{g_1}\otimes L_{g_2} \to L_{g_1g_2}\to L_{\partial c\cdot g_1g_2}.$$
This can be done, as here, by adding the annulus after the `pair of pants' or adding it on the first component, somewhat as in Figure \ref{pantsreduction}. The two formal $\mathcal{C}$-cobordisms are equivalent.

We can, of course, reverse the orientation to get
$$L_{g_1g_2} \to L_{g_1} \otimes L_{g_2},$$
a `comultiplication'.  It is fairly standard that this comultiplication is `redundant' as it can be recovered from the annuli and a suitable `positive pair of pants', see, for instance, the brief argument given in  section 5.1 of \cite{turaev:hqft1}.
\end{itemize}
To sum up, for a formal $\mathcal{C}$-HQFT, the passage from a `background' 1-type $B = K(P,1)$, and therefore from the known case with  a trivial top group in $\mathcal{C}$, to the model for a general 2-type $B = B \mathcal{C}$, we require merely the addition of extra linear isomorphisms $L_g \to L_{\partial c\cdot g}$ for $c$ in the top group of the crossed module, $\mathcal{C}$.  This structure is therefore very similar to that of a $\pi$-algebra (cf. \cite{turaev:hqft1}) for $\pi = P/\partial C$. We turn to this structure so as to be able to exhibit this new feature more fully.  This will give some insight into the connections between the structure  of the formal HQFT and  the background crossed module, $\mathcal{C}$.


\section{Crossed $\mathcal{C}$-algebras}

In \cite{turaev:hqft1}, the second author classified (1+1)-HQFTs with background a $K(\pi,1)$ in terms of crossed group-algebras. These were  generalisations of classical group algebras with many of the same features, but `twisted' by an action.  In \cite{B&T}, M. Brightwell and P. Turner examined the analogous case when the background is a $K(G,2)$ for $G$ an Abelian group, and classified them in terms of $G$-Frobenius algebras, that is, Frobenius algebras with a $G$-action.  In this section we will summarise both types of algebra before introducing a new type, `crossed $\mathcal{C}$-algebras', which combines features of both and which will classify formal HQFTs as above.
\subsection{Frobenius algebras}
First some background (adapted from \cite{rodrigues}) on Frobenius objects and Frobenius algebras.

Let $\mathcal{A}$ be a symmetric monoidal category with monoidal structure denoted $\otimes$ and with $\mathbb{K}$ as unit.  We say $\mathcal{A}$ has a (left) duality structure if for each object $A$, there is an object $A^*$, the dual of $A$, and morphisms
$$b_A : \mathbb{K}\to A \otimes A^*,$$
$$d_A :A^* \otimes A \to \mathbb{K}$$
such that \\
(i)\hspace{2.5cm} $ (A \stackrel{\cong}{\to}  \mathbb{K}\otimes A \stackrel{b_A\otimes A}{\to}A\otimes A^*\otimes A \stackrel{A\otimes d_A}{\to}A\otimes  \mathbb{K} \stackrel{\cong}{\to} A) = Id_{A}$\\
and\\
(ii) \hspace{2cm}$ (A^* \stackrel{\cong}{\to} A^*\otimes \mathbb{K} \stackrel{A^*\otimes b_A}{\to}A^*\otimes A\otimes A^* \stackrel{d_A\otimes A^*}{\to} \mathbb{K}\otimes A^* \stackrel{\cong}{\to} A^*) = Id_{A^*},$\\
where the unlabelled isomorphisms are the structural isomorphisms of $\mathcal{A}$ corresponding to $\mathbb{K}$ being a left and right unit for $\otimes$. 

The assignment of $A^*$ to $A$ extends to give a functor from $\mathcal{A}$ to $\mathcal{A}^{op}$, the opposite category.  If $f: A \to B$ is a morphism in $\mathcal{A}$, its dual or adjoint morphism $f^* : B^* \to A^*$ is given by the composition
$$B^*\stackrel{\cong}{\to}B^*\otimes \mathbb{K}\stackrel{B^*\otimes b_A}{\to} B^*\otimes A \otimes A^* \stackrel{B^*\otimes f \otimes A^*}{\to}B^*\otimes B\otimes A^*\stackrel{d_B\otimes A^*}{\to}\mathbb{K}\otimes A^*\stackrel{\cong}{\to}A^*.$$

If $\mathcal{A}$ has a duality structure as above, a \emph{Frobenius object} in $\mathcal{A}$ consists of\\
\vspace{-\baselineskip}

\begin{itemize}
\item an object $A$ of $\mathcal{A}$;
\item a `multiplication' morphism $\mu : A \otimes A \to A$;
\item a `unit' morphism $\eta : \mathbb{K} \to A$ such that $(A, \mu, \eta)$ is a monoid in $(\mathcal{A}, \otimes)$; \\
 \hspace*{-1.2cm} and
\item a symmetric `inner product' morphism,
$$\rho : A \otimes A \to \mathbb{K},$$
\end{itemize}
\vspace{-\baselineskip}
such that (i)
\hspace*{4cm}\xymatrix{A\otimes A \otimes A \ar[r]^{\quad A \otimes \mu}\ar[d]_{\mu \otimes A} &A \otimes A\ar[d]^\rho\\
A\otimes A \ar[r]_\rho & \mathbb{K}}\\
commutes (so writing $\mu(a,b) = ab$, $\rho(ab,c) = \rho(a,bc)$),\\
and \\
(ii) $\rho$ is non-degenerate, i.e., the following two induced maps from $A$ to $A^*$ are isomorphisms:
$$A\stackrel{\cong}{\to} A \otimes \mathbb{K} \stackrel{A\otimes b_A}{\to} A\otimes A\otimes A^* \stackrel{\rho \otimes A^*}{\to} \mathbb{K}\otimes A^* \stackrel{\cong}{\to}A ^*$$
and
$$A\stackrel{\cong}{\to} \mathbb{K}\otimes A \stackrel{\rho^* \otimes A}{\to} A^*\otimes A^* \otimes A \stackrel{A^* \otimes d_A}{\to}A^*\otimes \mathbb{K}\stackrel{\cong}{\to }A^*.$$
(This second composite tacitly uses the isomorphisms $(A\otimes A)^* \cong A^*\otimes A^*$, and $\mathbb{K}^* \cong \mathbb{K}$ which hold since $\mathcal{A}$ is assumed to be symmetric monoidal.)

\

\textbf{Examples}

Frobenius objects in the category $(Vect, \otimes)$ or more generally $(Mod\;R),\otimes )$ are Frobenius algebras in the usual sense.

\medskip

We next describe two variants of Frobenius algebras that arise when modelling HQFTs with $B$ a $K(\pi,1)$ and a $K(G,2)$ respectively.  (The detailed references for the two situations are \cite{turaev:hqft1} and \cite{B&T} or \cite{rodrigues}, respectively.)

\subsection{Crossed $\pi$-algebras}
(Based on parts of \cite{turaev:hqft1})

Here $\pi$ will be a group corresponding to $\pi_1(B)$ if $B$ is a 1-type. 
\medskip

\textbf{Definition}

A \emph{graded $\pi$-algebra} or \emph{$\pi$-algebra} over a field $\mathbb{K}$ is an associative algebra $L$ over $\mathbb{K}$ with a splitting
$$L = \bigoplus_{g\in \pi} L_g,$$ 
as a direct sum of projective $\mathbb{K}$-modules of finite type such that\\
(i) $L_gL_h \subseteq L_{gh}$ for any $g,h \in \pi$ (so, if $\ell_1$ is graded $g$, and $\ell_2$ is graded $h$, then $\ell_1\ell_2$ is graded $gh$),\\
and \\
(ii) $L$ has a unit $1 = 1_L\in L_1$ for 1, the identity element of $\pi$.

\

\textbf{Example}

The group algebra $\mathbb{K}[\pi]$ has a $\pi$-algebra structure as has $A[\pi]= A\otimes_\mathbb{K}\mathbb{K}[\pi]$ for any associative $\mathbb{K}$-algebra $A$.  Multiplication in $A[\pi]$ is given by $(ag)(bh) = (ab)(gh)$ for $a,b \in A$, $g,h \in \pi$, in the obvious notation.

\medskip

\textbf{Definition}

A \emph{Frobenius $\pi$-algebra} is a $\pi$-algebra $L$ together with a symmetric $\mathbb{K}$-bilinear form
$$\rho : L\otimes L \to \mathbb{K}$$
such that \\
(i) \quad $\rho(L_g\otimes L_h) = 0$ if $h \neq g^{-1}$;\\
(ii) the restriction of $\rho$ to $L_g \otimes L_{g^{-1}}$ is non-degenerate for each $g\in \pi$, (so $L_{g^{-1}} \cong  L_g^*$, the dual of $L_g$);\\
and\\
(iii) \quad $\rho(ab,c) = \rho(a,bc)$ for any $a,b,c \in L$.

\medskip

\textbf{Example: continued}

The group algebra, $L = \mathbb{K}[\pi]$, is a Frobenius $\pi$-algebra with $\rho(g,h) = 1$ if $gh = 1$, and 0 otherwise.

\medskip

Finally the notion of  crossed $\pi$-algebra combines the above with an action of $\pi$ on $L$, explicitly:

\textbf{Definition}

A \emph{crossed $\pi$-algebra} over $\mathbb{K}$ is a Frobenius $\pi$-algebra over $\mathbb{K}$ together with a group homomorphism
$$\phi: \pi \to Aut(L)$$
satisfying:\\
(i) if $g\in \pi$ and we write $\phi_g = \phi(g)$ for the corresponding automorphism of $L$, then $\phi_g$ preserves $\rho$, (i.e. $\rho(\phi_ga,\phi_gb) = \rho(a,b)$) and 
$$\phi_g(L_h) \subseteq L_{ghg^{-1}}$$
for all $h\in \pi$;\\
(ii) $\phi_g|_{L_g} = id$ for all $g\in \pi$;\\
(iii) for any $g,h \in \pi$, $a\in L_g$, $b\in L_h$, $\phi_h(a)b = ba$\label{Palg3};\\
(iv) for any $g,h \in \pi$ and $c \in L_{ghg^{-1}h^{-1}}$,
$$Tr(c\phi_h : L_g \to L_g) = Tr(\phi_{g^{-1}}c : L_h \to L_h),$$
where $Tr$ denotes the $\mathbb{K}$-valued trace of the endomorphism. (The homomorphism $c\phi_h$ sends $a\in L_g$ to $c\phi_h(a) \in L_g$, whilst $(\phi_{g^{-1}}c)(b) = \phi_{g^{-1}}(cb)$ for $c \in L_h$.)

\

\textbf{Example revisited:}  It is easily checked, see \cite{turaev:hqft1}, that $\mathbb{K}[\pi]$ us a crossed $\pi$-algebra.
\subsection{$G$-Frobenius algebras}
(See \cite{B&T,rodrigues})

(In this subsection, $G$ will denote an Abelian group.)

We have defined a Frobenius object in a symmetric monoidal category $\mathcal{A}$.  A \emph{$G$-Frobenius object} in $\mathcal{A}$ is a Frobenius object $A$ together with a homomorphism
$$G \to End(A).$$
In the cases $\mathcal{A} = (Vect, \otimes)$ or $(Mod R,\otimes)$, the resulting concept is that of a \emph{$G$-Frobenius algebra}.  Examination of the action shows that if we write $g\cdot a$ for the action of $g$ on an element $a \in A$,
$$a(g\cdot b) = g\cdot(ab) = (g\cdot a)b$$
and
$$\rho(a,g\cdot b) = \rho(g\cdot a,b).$$
As $A$ is a unital algebra, $$g\cdot v = g\cdot 1 v = (g\cdot 1)v,$$
so the action actually comes from a morphism of monoids $$G\to A$$
$$g\to g\cdot 1,$$
and $g\cdot 1$ is in the center of $A$.

\subsection{Crossed $\mathcal{C}$-algebras: the definition}
We now turn to the general case with $\mathcal{C}$, the crossed module $(C,P,\partial)$ as earlier.  Any specification of  formal $\mathcal{C}$-maps on  simplicial complexes must include formal maps  in which $C$ itself plays no part, corresponding to the 2-cells being all labelled $1_C$.  We thus should expect an associated crossed $P$-algebra underlying any crossed $\mathcal{C}$-algebra.  The additional structure is then  that given by the annuli or cylinders $(Cyl_{\epsilon,\mu}; c,g_1,g_2)$. We saw earlier that this collection of  operations could be reduced further to the case $g_2 = 1$ and $c \neq 1$, and, in fact, the only ones we actually need are with $g_1 = 1$ as well, the general case being a composite of this with the unit on the left and the `pair of pants' multiplication. (The general case gives an isomorphism
$$\theta_{(c,g)} : L_g \to L_{\partial c \cdot g}$$
and we can build this up by
$$L_g \to \mathbb{K}\otimes L_g \to L_1 \otimes L_g \stackrel{\theta_{(c,1)} \otimes L_g}{\longrightarrow} L_{\partial c}\otimes L_g\stackrel{\mu}{\to}L_{\partial c \cdot g},$$
where the third morphism is that given by that special case $g = 1$.  We say that $\theta_{(c,g)}$ is obtained by `translation' from $\theta_{(c,1)}$.)

The extra structure can therefore be thought of as a collection of isomorphisms
$$\Theta_{\mathcal{C}} = \{\theta_{(c,1)} : L_1 \to L_{\partial c} ~ : ~ c\in C\}.$$
It is worth noting that if $C = \{1\}$, the resulting structure reduces to that of a crossed $P$-algebra and if $P = 1$ and $C$ is just an Abelian group then the $\theta_{(c,1)} : L_1 \to L_1$ are just automorphisms of $L_1$, which is itself just a Frobenius algebra.

\

This structure of extra specified automorphisms does not immediately tell us how to retrieve the structure given by the $\mathcal{C}$-discs.  Those gave  linear maps 
$$\ell_c : \mathbb{K} \to L_{\partial c}.$$
We can, however, recover them from $\ell_1 : \mathbb{K} \to L_1$, which was part of the crossed $P$-algebra structure, together with $\theta_{(c,1)} : L_1 \to L_{\partial c}$, but conversely given the $\ell_c$, we can recover the $\theta_{(c,g)}$: 
\begin{proposition}\hspace*{1mm}\\
The composite
$$L_g \stackrel{\cong}{\to} \mathbb{K} \otimes L_g \stackrel{\ell_c\otimes L_g}{\to}L_{\partial c}\otimes L_g \stackrel{\mu}{\to} L_{\partial c \cdot g}$$
is equal to $\theta_{(c,g)}$.
\end{proposition}
\textbf{Proof}

We can realise this composite by a $\mathcal{C}$-cobordism

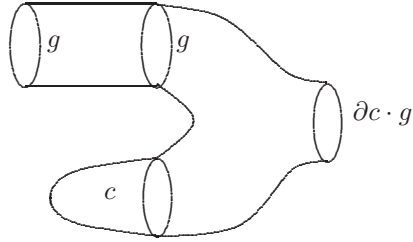
\begin{center}
\font\thinlinefont=cmr5
\begingroup\makeatletter\ifx\SetFigFont\undefined
\def\x#1#2#3#4#5#6#7\relax{\def\x{#1#2#3#4#5#6}}%
\expandafter\x\fmtname xxxxxx\relax \def\y{splain}%
\ifx\x\y   
\gdef\SetFigFont#1#2#3{%
  \ifnum #1<17\tiny\else \ifnum #1<20\small\else
  \ifnum #1<24\normalsize\else \ifnum #1<29\large\else
  \ifnum #1<34\Large\else \ifnum #1<41\LARGE\else
     \huge\fi\fi\fi\fi\fi\fi
  \csname #3\endcsname}%
\else
\gdef\SetFigFont#1#2#3{\begingroup
  \count@#1\relax \ifnum 25<\count@\count@25\fi
  \def\x{\endgroup\@setsize\SetFigFont{#2pt}}%
  \expandafter\x
    \csname \romannumeral\the\count@ pt\expandafter\endcsname
    \csname @\romannumeral\the\count@ pt\endcsname
  \csname #3\endcsname}%
\fi
\fi\endgroup
\mbox{\beginpicture
\setcoordinatesystem units <0.5000cm,0.5000cm>
\unitlength=0.5000cm
\linethickness=1pt
\setplotsymbol ({\makebox(0,0)[l]{\tencirc\symbol{'160}}})
\setshadesymbol ({\thinlinefont .})
\setlinear
%
%
\linethickness= 0.500pt
\setplotsymbol ({\thinlinefont .})
\plot  7.489 22.648  7.495 22.644 /
\plot  7.495 22.644  7.508 22.636 /
\plot  7.508 22.636  7.531 22.619 /
\plot  7.531 22.619  7.567 22.591 /
\plot  7.567 22.591  7.616 22.557 /
\plot  7.616 22.557  7.673 22.515 /
\plot  7.673 22.515  7.741 22.464 /
\plot  7.741 22.464  7.813 22.409 /
\plot  7.813 22.409  7.889 22.350 /
\plot  7.889 22.350  7.965 22.291 /
\plot  7.965 22.291  8.037 22.231 /
\plot  8.037 22.231  8.107 22.172 /
\plot  8.107 22.172  8.170 22.117 /
\plot  8.170 22.117  8.227 22.062 /
\plot  8.227 22.062  8.278 22.013 /
\plot  8.278 22.013  8.321 21.965 /
\plot  8.321 21.965  8.359 21.920 /
\plot  8.359 21.920  8.388 21.878 /
\plot  8.388 21.878  8.412 21.836 /
\plot  8.412 21.836  8.431 21.797 /
\plot  8.431 21.797  8.441 21.759 /
\plot  8.441 21.759  8.448 21.721 /
\plot  8.448 21.721  8.450 21.683 /
\plot  8.450 21.683  8.446 21.647 /
\plot  8.446 21.647  8.439 21.611 /
\plot  8.439 21.611  8.426 21.575 /
\plot  8.426 21.575  8.412 21.539 /
\plot  8.412 21.539  8.388 21.501 /
\plot  8.388 21.501  8.361 21.461 /
\plot  8.361 21.461  8.329 21.419 /
\plot  8.329 21.419  8.289 21.374 /
\plot  8.289 21.374  8.244 21.328 /
\plot  8.244 21.328  8.194 21.277 /
\plot  8.194 21.277  8.136 21.224 /
\plot  8.136 21.224  8.075 21.169 /
\plot  8.075 21.169  8.007 21.112 /
\plot  8.007 21.112  7.938 21.052 /
\plot  7.938 21.052  7.863 20.993 /
\plot  7.863 20.993  7.791 20.934 /
\plot  7.791 20.934  7.719 20.879 /
\plot  7.719 20.879  7.654 20.826 /
\plot  7.654 20.826  7.595 20.781 /
\plot  7.595 20.781  7.542 20.743 /
\plot  7.542 20.743  7.501 20.712 /
\plot  7.501 20.712  7.472 20.688 /
\plot  7.472 20.688  7.453 20.676 /
\plot  7.453 20.676  7.440 20.667 /
\plot  7.440 20.667  7.436 20.663 /
%
%
\linethickness= 0.500pt
\setplotsymbol ({\thinlinefont .})
\ellipticalarc axes ratio  0.373:1.04  360 degrees 
from  7.86 19.632 center at  7.487 19.632
%
%
\linethickness= 0.500pt
\setplotsymbol ({\thinlinefont .})
\ellipticalarc axes ratio  0.373:1.04  360 degrees 
from  4.369 23.68 center at  3.97 23.68
%
\linethickness= 0.500pt
\setplotsymbol ({\thinlinefont .})
\ellipticalarc axes ratio  0.373:1.04  360 degrees 
from  7.86 23.68 center at  7.46 23.68
%
%
\linethickness= 0.500pt
\setplotsymbol ({\thinlinefont .})
\putrule from  3.97 22.6 to  7.46 22.6
%
%
\linethickness= 0.500pt
\setplotsymbol ({\thinlinefont .})
\putrule from  3.97 24.8 to  7.46 24.8
%
%
\linethickness= 0.500pt
\setplotsymbol ({\thinlinefont .})
\ellipticalarc axes ratio  0.373:1.033  360 degrees 
	from 12.383 21.618 center at 12.010 21.618
%
%
\linethickness= 0.500pt
\setplotsymbol ({\thinlinefont .})
\plot  7.46 22.6  7.476 22.63 /
\putrule from  7.476 24.763 to  7.495 24.763
\plot  7.495 24.763  7.523 24.759 /
\plot  7.523 24.759  7.563 24.757 /
\plot  7.563 24.757  7.616 24.750 /
\plot  7.616 24.750  7.679 24.744 /
\plot  7.679 24.744  7.751 24.735 /
\plot  7.751 24.735  7.834 24.727 /
\plot  7.834 24.727  7.921 24.716 /
\plot  7.921 24.716  8.014 24.706 /
\plot  8.014 24.706  8.111 24.693 /
\plot  8.111 24.693  8.206 24.680 /
\plot  8.206 24.680  8.302 24.666 /
\plot  8.302 24.666  8.395 24.651 /
\plot  8.395 24.651  8.486 24.636 /
\plot  8.486 24.636  8.570 24.619 /
\plot  8.570 24.619  8.653 24.604 /
\plot  8.653 24.604  8.731 24.585 /
\plot  8.731 24.585  8.805 24.568 /
\plot  8.805 24.568  8.877 24.549 /
\plot  8.877 24.549  8.943 24.530 /
\plot  8.943 24.530  9.006 24.509 /
\plot  9.006 24.509  9.068 24.486 /
\plot  9.068 24.486  9.125 24.462 /
\plot  9.125 24.462  9.182 24.437 /
\plot  9.182 24.437  9.237 24.409 /
\plot  9.237 24.409  9.292 24.382 /
\plot  9.292 24.382  9.345 24.350 /
\plot  9.345 24.350  9.398 24.318 /
\plot  9.398 24.318  9.449 24.285 /
\plot  9.449 24.285  9.502 24.246 /
\plot  9.502 24.246  9.555 24.208 /
\plot  9.555 24.208  9.608 24.168 /
\plot  9.608 24.168  9.663 24.126 /
\plot  9.663 24.126  9.716 24.083 /
\plot  9.716 24.083  9.771 24.037 /
\plot  9.771 24.037  9.823 23.990 /
\plot  9.823 23.990  9.878 23.942 /
\plot  9.878 23.942  9.934 23.891 /
\plot  9.934 23.891  9.989 23.840 /
\plot  9.989 23.840 10.044 23.787 /
\plot 10.044 23.787 10.099 23.736 /
\plot 10.099 23.736 10.154 23.683 /
\plot 10.154 23.683 10.207 23.630 /
\plot 10.207 23.630 10.259 23.578 /
\plot 10.259 23.578 10.312 23.527 /
\plot 10.312 23.527 10.363 23.476 /
\plot 10.363 23.476 10.414 23.425 /
\plot 10.414 23.425 10.463 23.379 /
\plot 10.463 23.379 10.511 23.332 /
\plot 10.511 23.332 10.558 23.285 /
\plot 10.558 23.285 10.602 23.243 /
\plot 10.602 23.243 10.647 23.203 /
\plot 10.647 23.203 10.691 23.163 /
\plot 10.691 23.163 10.734 23.127 /
\plot 10.734 23.127 10.774 23.091 /
\plot 10.774 23.091 10.814 23.059 /
\plot 10.814 23.059 10.854 23.027 /
\plot 10.854 23.027 10.907 22.989 /
\plot 10.907 22.989 10.958 22.953 /
\plot 10.958 22.953 11.011 22.921 /
\plot 11.011 22.921 11.066 22.894 /
\plot 11.066 22.894 11.121 22.866 /
\plot 11.121 22.866 11.178 22.843 /
\plot 11.178 22.843 11.237 22.822 /
\plot 11.237 22.822 11.301 22.805 /
\plot 11.301 22.805 11.369 22.786 /
\plot 11.369 22.786 11.438 22.771 /
\plot 11.438 22.771 11.513 22.758 /
\plot 11.513 22.758 11.587 22.746 /
\plot 11.587 22.746 11.663 22.735 /
\plot 11.663 22.735 11.735 22.725 /
\plot 11.735 22.725 11.800 22.718 /
\plot 11.800 22.718 11.860 22.712 /
\plot 11.860 22.712 11.906 22.708 /
\plot 11.906 22.708 11.942 22.703 /
\putrule from 11.942 22.703 to 11.968 22.703
\plot 11.968 22.703 11.980 22.701 /
\putrule from 11.980 22.701 to 11.987 22.701
%
%
\linethickness= 0.500pt
\setplotsymbol ({\thinlinefont .})
\putrule from  7.461 18.599 to  7.465 18.599
\putrule from  7.465 18.599 to  7.476 18.599
\putrule from  7.476 18.599 to  7.497 18.599
\plot  7.497 18.599  7.527 18.601 /
\plot  7.527 18.601  7.571 18.603 /
\plot  7.571 18.603  7.626 18.605 /
\plot  7.626 18.605  7.696 18.608 /
\plot  7.696 18.608  7.777 18.612 /
\plot  7.777 18.612  7.868 18.616 /
\plot  7.868 18.616  7.969 18.620 /
\plot  7.969 18.620  8.075 18.627 /
\plot  8.075 18.627  8.185 18.633 /
\plot  8.185 18.633  8.299 18.641 /
\plot  8.299 18.641  8.412 18.650 /
\plot  8.412 18.650  8.524 18.658 /
\plot  8.524 18.658  8.632 18.669 /
\plot  8.632 18.669  8.738 18.680 /
\plot  8.738 18.680  8.837 18.690 /
\plot  8.837 18.690  8.932 18.703 /
\plot  8.932 18.703  9.023 18.716 /
\plot  9.023 18.716  9.108 18.730 /
\plot  9.108 18.730  9.188 18.745 /
\plot  9.188 18.745  9.265 18.762 /
\plot  9.265 18.762  9.335 18.779 /
\plot  9.335 18.779  9.402 18.798 /
\plot  9.402 18.798  9.466 18.817 /
\plot  9.466 18.817  9.525 18.838 /
\plot  9.525 18.838  9.582 18.862 /
\plot  9.582 18.862  9.637 18.887 /
\plot  9.637 18.887  9.690 18.915 /
\plot  9.690 18.915  9.741 18.942 /
\plot  9.741 18.942  9.798 18.978 /
\plot  9.798 18.978  9.853 19.016 /
\plot  9.853 19.016  9.908 19.056 /
\plot  9.908 19.056  9.961 19.097 /
\plot  9.961 19.097 10.014 19.143 /
\plot 10.014 19.143 10.065 19.190 /
\plot 10.065 19.190 10.116 19.238 /
\plot 10.116 19.238 10.166 19.291 /
\plot 10.166 19.291 10.217 19.344 /
\plot 10.217 19.344 10.266 19.399 /
\plot 10.266 19.399 10.315 19.456 /
\plot 10.315 19.456 10.363 19.514 /
\plot 10.363 19.514 10.410 19.571 /
\plot 10.410 19.571 10.456 19.630 /
\plot 10.456 19.630 10.501 19.689 /
\plot 10.501 19.689 10.545 19.746 /
\plot 10.545 19.746 10.588 19.804 /
\plot 10.588 19.804 10.630 19.861 /
\plot 10.630 19.861 10.670 19.916 /
\plot 10.670 19.916 10.710 19.969 /
\plot 10.710 19.969 10.748 20.019 /
\plot 10.748 20.019 10.784 20.068 /
\plot 10.784 20.068 10.820 20.115 /
\plot 10.820 20.115 10.856 20.159 /
\plot 10.856 20.159 10.890 20.199 /
\plot 10.890 20.199 10.924 20.240 /
\plot 10.924 20.240 10.956 20.276 /
\plot 10.956 20.276 10.990 20.309 /
\plot 10.990 20.309 11.032 20.350 /
\plot 11.032 20.350 11.074 20.386 /
\plot 11.074 20.386 11.117 20.417 /
\plot 11.117 20.417 11.163 20.447 /
\plot 11.163 20.447 11.210 20.472 /
\plot 11.210 20.472 11.261 20.494 /
\plot 11.261 20.494 11.314 20.515 /
\plot 11.314 20.515 11.371 20.532 /
\plot 11.371 20.532 11.432 20.546 /
\plot 11.432 20.546 11.496 20.561 /
\plot 11.496 20.561 11.563 20.572 /
\plot 11.563 20.572 11.633 20.582 /
\plot 11.633 20.582 11.703 20.591 /
\plot 11.703 20.591 11.771 20.597 /
\plot 11.771 20.597 11.834 20.602 /
\plot 11.834 20.602 11.889 20.606 /
\plot 11.889 20.606 11.936 20.608 /
\putrule from 11.936 20.608 to 11.970 20.608
\plot 11.970 20.608 11.993 20.610 /
\putrule from 11.993 20.610 to 12.006 20.610
\putrule from 12.006 20.610 to 12.012 20.610
%
%
\linethickness= 0.500pt
\setplotsymbol ({\thinlinefont .})
\putrule from  7.436 20.690 to  7.430 20.690
\plot  7.430 20.690  7.417 20.688 /
\plot  7.417 20.688  7.396 20.686 /
\plot  7.396 20.686  7.360 20.682 /
\plot  7.360 20.682  7.311 20.676 /
\plot  7.311 20.676  7.247 20.667 /
\plot  7.247 20.667  7.169 20.659 /
\plot  7.169 20.659  7.080 20.648 /
\plot  7.080 20.648  6.981 20.635 /
\plot  6.981 20.635  6.875 20.621 /
\plot  6.875 20.621  6.763 20.606 /
\plot  6.763 20.606  6.648 20.591 /
\plot  6.648 20.591  6.532 20.576 /
\plot  6.532 20.576  6.420 20.559 /
\plot  6.420 20.559  6.310 20.544 /
\plot  6.310 20.544  6.204 20.530 /
\plot  6.204 20.530  6.104 20.513 /
\plot  6.104 20.513  6.009 20.498 /
\plot  6.009 20.498  5.920 20.483 /
\plot  5.920 20.483  5.838 20.470 /
\plot  5.838 20.470  5.759 20.455 /
\plot  5.759 20.455  5.690 20.441 /
\plot  5.690 20.441  5.622 20.428 /
\plot  5.622 20.428  5.560 20.413 /
\plot  5.560 20.413  5.501 20.398 /
\plot  5.501 20.398  5.446 20.384 /
\plot  5.446 20.384  5.395 20.369 /
\plot  5.395 20.369  5.347 20.354 /
\plot  5.347 20.354  5.300 20.337 /
\plot  5.300 20.337  5.241 20.316 /
\plot  5.241 20.316  5.184 20.290 /
\plot  5.184 20.290  5.129 20.267 /
\plot  5.129 20.267  5.078 20.240 /
\plot  5.078 20.240  5.029 20.212 /
\plot  5.029 20.212  4.983 20.182 /
\plot  4.983 20.182  4.938 20.153 /
\plot  4.938 20.153  4.898 20.121 /
\plot  4.898 20.121  4.860 20.087 /
\plot  4.860 20.087  4.824 20.053 /
\plot  4.824 20.053  4.792 20.017 /
\plot  4.792 20.017  4.763 19.983 /
\plot  4.763 19.983  4.737 19.947 /
\plot  4.737 19.947  4.714 19.911 /
\plot  4.714 19.911  4.695 19.876 /
\plot  4.695 19.876  4.678 19.840 /
\plot  4.678 19.840  4.665 19.806 /
\plot  4.665 19.806  4.655 19.772 /
\plot  4.655 19.772  4.646 19.738 /
\plot  4.646 19.738  4.640 19.706 /
\plot  4.640 19.706  4.638 19.672 /
\plot  4.638 19.672  4.635 19.641 /
\putrule from  4.635 19.641 to  4.635 19.611
\plot  4.635 19.611  4.638 19.581 /
\plot  4.638 19.581  4.642 19.552 /
\plot  4.642 19.552  4.646 19.522 /
\plot  4.646 19.522  4.655 19.492 /
\plot  4.655 19.492  4.665 19.461 /
\plot  4.665 19.461  4.678 19.431 /
\plot  4.678 19.431  4.693 19.399 /
\plot  4.693 19.399  4.710 19.368 /
\plot  4.710 19.368  4.731 19.336 /
\plot  4.731 19.336  4.754 19.304 /
\plot  4.754 19.304  4.779 19.274 /
\plot  4.779 19.274  4.809 19.243 /
\plot  4.809 19.243  4.839 19.213 /
\plot  4.839 19.213  4.875 19.185 /
\plot  4.875 19.185  4.911 19.156 /
\plot  4.911 19.156  4.949 19.128 /
\plot  4.949 19.128  4.989 19.103 /
\plot  4.989 19.103  5.033 19.078 /
\plot  5.033 19.078  5.078 19.054 /
\plot  5.078 19.054  5.127 19.033 /
\plot  5.127 19.033  5.177 19.010 /
\plot  5.177 19.010  5.228 18.991 /
\plot  5.228 18.991  5.283 18.970 /
\plot  5.283 18.970  5.328 18.957 /
\plot  5.328 18.957  5.372 18.942 /
\plot  5.372 18.942  5.421 18.929 /
\plot  5.421 18.929  5.472 18.917 /
\plot  5.472 18.917  5.524 18.904 /
\plot  5.524 18.904  5.582 18.891 /
\plot  5.582 18.891  5.643 18.876 /
\plot  5.643 18.876  5.709 18.864 /
\plot  5.709 18.864  5.779 18.851 /
\plot  5.779 18.851  5.853 18.838 /
\plot  5.853 18.838  5.933 18.826 /
\plot  5.933 18.826  6.020 18.813 /
\plot  6.020 18.813  6.109 18.798 /
\plot  6.109 18.798  6.206 18.783 /
\plot  6.206 18.783  6.306 18.771 /
\plot  6.306 18.771  6.409 18.756 /
\plot  6.409 18.756  6.517 18.741 /
\plot  6.517 18.741  6.625 18.726 /
\plot  6.625 18.726  6.735 18.713 /
\plot  6.735 18.713  6.841 18.699 /
\plot  6.841 18.699  6.945 18.686 /
\plot  6.945 18.686  7.040 18.673 /
\plot  7.040 18.673  7.129 18.663 /
\plot  7.129 18.663  7.207 18.654 /
\plot  7.207 18.654  7.275 18.646 /
\plot  7.275 18.646  7.330 18.639 /
\plot  7.330 18.639  7.372 18.635 /
\plot  7.372 18.635  7.402 18.631 /
\plot  7.402 18.631  7.421 18.629 /
\plot  7.421 18.629  7.432 18.627 /
\putrule from  7.432 18.627 to  7.436 18.627
%
%
\put{$g$} [lB] at  4.551 23.654
%
%
\put{$g$} [lB] at  7.990 23.654
%
%
\put{$c$} [lB] at  6.060 19.605
%
%
\put{$\partial c \cdot g$} [lB] at 12.675 21.696
\linethickness=0pt
\putrectangle corners at  3.607 24.790 and 12.675 18.574
\endpicture}
\end{center}
but this is equivalent to the $\mathcal{C}$-annulus that gives us $\theta_{(c,g)}.$\hfill $\square$

\

As before we will write $\tilde{c} = \ell_c(1) \in L_{\partial c}.$
\begin{corollary}\hspace*{1mm}\\
For any $c\in C$, $g \in P$ and for $x \in L_g$,
$$\theta_{(c,g)}(x) = \tilde{c} \cdot x,$$
where $\cdot$ denotes the product in the algebra structure of $L = \bigoplus_{h\in P}L_h$.\hfill $\square$

\end{corollary}
Abstracting this extra structure, we get:

\

\textbf{Definition}

Let $\mathcal{C} = (C,P,\partial)$ be a crossed module.  A \emph{crossed $\mathcal{C}$-algebra} consists of a crossed $P$-algebra, $L = \bigoplus_{g\in P}L_g$, together with elements $\tilde{c} \in L_{\partial c}$, for $c \in C$, such that \begin{enumerate}[(a)]
\item $\tilde{1} = 1\in L_1$;
\item for $c,c^\prime \in C$, $\widetilde{(c^\prime c)} = \tilde{c^\prime}\cdot \tilde{c}$;
\item for any $h\in P$, $\phi_h(\tilde{c}) = \widetilde{{}^hc}$.
\end{enumerate}
We note for future use that the first two conditions make `tilderisation' into a group homomorphism $(\tilde{~}) : C \to U(L)$, the group of units of the algebra, $L$.

There is an obvious notion of morphism of crossed $\mathcal{C}$-algebras, which we will examine in more detail in section 6.2.  There is, of course, a linked notion of isomorphism of crossed $\mathcal{C}$-algebras which will enable us in section 5.1 to state our main theorem.

\

The two special cases with $C = 1$ and, for Abelian $C$, with $P = 1$ correspond, of course, to crossed $P$-algebras and $C$-Frobenius algebras respectively. An interesting special case of the general form is when $C$ is a $P$-module and $\partial$ sends every element in $C$ to the identity of $P$.  In this case we have an object that could be described as a $C$-crossed $P$-algebra! It consists of a crossed $P$-algebra together with a $C$-action by multiplication by central elements. This results in a very weak mixing of the two structures. The important thing to note is that the general form is more highly structured as the twisting in the crossed modules, in general, can result in non-central elements amongst the  $\tilde{c} $s. 

\section{A classification of formal $\mathcal{C}$-HQFTs}
\subsection{Main Theorem}
\begin{theorem}\hspace*{1cm}\\
There is a canonical bijection between isomorphism classes of formal  2-dimensional HQFTs based on a crossed module $\mathcal{C}$ and isomorphism classes of crossed $\mathcal{C}$-algebras.
\end{theorem}

More explicitly:
\begin{theorem}\hspace*{1mm}\\
a)  For any formal 2-dimensional HQFT, $\tau$, based on $\mathcal{C}$, the crossed $P$-algebra, $L = \bigoplus_{g\in P}L_g$, having $L_g = \tau(g)$, is a crossed $\mathcal{C}$-algebra, where for $c \in C$, $\tilde{c} = \ell_c(1)$ (notation as above).\\
b)  Given any crossed $\mathcal{C}$-algebra, $L = \bigoplus_{g\in P}L_g$, there is a formal 2-dimensional HQFT, $\tau$, based on $\mathcal{C}$ yielding $L$ as its crossed $\mathcal{C}$-algebra, up to isomorphism.
\end{theorem}
Before we launch into the proof of this result some comments are in order.  We will need to understand the combination of formal $\mathcal{C}$-cobordisms in some detail before the proof can be undertaken, however much of what we need will be an adaptation of the geometric ideas already used in the $K(G,1)$-case in \cite{turaev:hqft1}. 

\subsection{Combination of fragments  of $\mathcal{C}$-cobordisms}

We can schematically represent a fragment of a $\mathcal{C}$-cobordism by  a 2-cell
$$
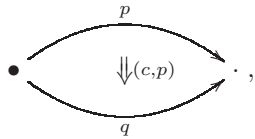
\xymatrix@+10pt{\bullet\rrtwocell<7>^p_q{~~~(c,p)} &&. },$$
with $q =\partial c\cdot p$ and initially $p$ and $q$ may be combinations of edge labels. The $\mathcal{C}$-cobordism is given by a cellular decomposition of the underlying surface, hence is made up of building blocks which are 2-cells. For instance, when we apply this in the analysis of the building blocks for $\mathcal{C}$-cobordisms, one case will correspond to the annulus or cylinder, $(Cyl_{\varepsilon,\mu};c,g,1)$ on page \pageref{CYL}. Thus 
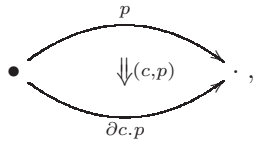
$$\xymatrix@+10pt{\bullet\rrtwocell<7>^p_{\partial c.p}{~~~(c,p)} &&. },$$
corresponds to the surface polygon:
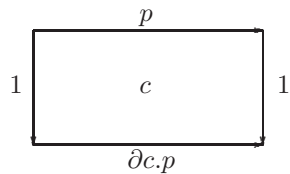
\begin{center}
\font\thinlinefont=cmr5

\begingroup\makeatletter\ifx\SetFigFont\undefined%
\gdef\SetFigFont#1#2#3#4#5{%
  \reset@font\fontsize{#1}{#2pt}%
  \fontfamily{#3}\fontseries{#4}\fontshape{#5}%
  \selectfont}%
\fi\endgroup%
\mbox{\beginpicture
\setcoordinatesystem units <.40000cm,.40000cm>
\unitlength=.40000cm
\linethickness=1pt
\setplotsymbol ({\makebox(0,0)[l]{\tencirc\symbol{'160}}})
\setshadesymbol ({\thinlinefont .})
\setlinear
%
%
\put{$\partial c.p$} [lB] at 8.2 17.
%
%
\linethickness= 0.500pt
\setplotsymbol ({\thinlinefont .})
{\putrule from 12.700 21.590 to 12.700 17.780
%
%
\plot 12.636 18.034 12.700 17.780 12.764 18.034 /
}%
%
%
\linethickness= 0.500pt
\setplotsymbol ({\thinlinefont .})
{\putrule from  5.080 17.780 to 12.700 17.780
%
%
\plot 12.446 17.716 12.700 17.780 12.446 17.844 /
}%
%
%
\linethickness= 0.500pt
\setplotsymbol ({\thinlinefont .})
{\putrule from  5.080 21.590 to  5.080 17.780
%
%
\plot  5.017 18.034  5.080 17.780  5.143 18.034 /
}%
%
%
\put{$p$} [lB] at  8.6 21.907
%
%
\put{$1$}[lB] at 13.176 19.5
%
%
\put{$c$}[lB] at  8.6 19.5
%
%
\put{$1$}[lB] at  4.286 19.5
%
%
\linethickness= 0.500pt
\setplotsymbol ({\thinlinefont .})
{\putrule from  5.080 21.590 to 12.700 21.590
%
%
\plot 12.446 21.526 12.700 21.590 12.446 21.654 /
}%
\linethickness=0pt
\putrectangle corners at  4.286 22.246 and 13.392 17.018
\endpicture}\end{center}
viewed as
\begin{center}
\font\thinlinefont=cmr5
\begingroup\makeatletter\ifx\SetFigFont\undefined
\def\x#1#2#3#4#5#6#7\relax{\def\x{#1#2#3#4#5#6}}%
\expandafter\x\fmtname xxxxxx\relax \def\y{splain}%
\ifx\x\y   
\gdef\SetFigFont#1#2#3{%
  \ifnum #1<17\tiny\else \ifnum #1<20\small\else
  \ifnum #1<24\normalsize\else \ifnum #1<29\large\else
  \ifnum #1<34\Large\else \ifnum #1<41\LARGE\else
     \huge\fi\fi\fi\fi\fi\fi
  \csname #3\endcsname}%
\else
\gdef\SetFigFont#1#2#3{\begingroup
  \count@#1\relax \ifnum 25<\count@\count@25\fi
  \def\x{\endgroup\@setsize\SetFigFont{#2pt}}%
  \expandafter\x
    \csname \romannumeral\the\count@ pt\expandafter\endcsname
    \csname @\romannumeral\the\count@ pt\endcsname
  \csname #3\endcsname}%
\fi
\fi\endgroup
\mbox
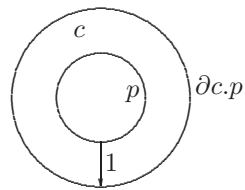
{\beginpicture
\setcoordinatesystem units <0.40000cm,0.40000cm>
\unitlength=0.40000cm
\linethickness=1pt
\setplotsymbol ({\makebox(0,0)[l]{\tencirc\symbol{'160}}})
\setshadesymbol ({\thinlinefont .})
\setlinear
%
%
\linethickness= 0.500pt
\setplotsymbol ({\thinlinefont .})
\ellipticalarc axes ratio  2.963:2.963  360 degrees 
	from  8.970 21.643 center at  6.007 21.643
%
%
\linethickness= 0.500pt
\setplotsymbol ({\thinlinefont .})
\putrule from  5.997 20.180 to  5.997 18.735
%
%
\plot  5.944 18.946  5.997 18.735  6.049 18.946 /
%
%
%
\put{$c$} [lB] at  5.093 23.664
%
%
\put{$p$} [lB] at  6.835 21.660
%
%
\put{$1$} [lB] at  6.128 19.2
%
%
\put{$\partial c.p$} [lB] at  9.129 21.749
%
%
\linethickness= 0.500pt
\setplotsymbol ({\thinlinefont .})
\ellipticalarc axes ratio  1.482:1.482  360 degrees 
	from  7.489 21.643 center at  6.007 21.643
\linethickness=0pt
\putrectangle corners at  3.027 24.621 and  9.129 18.665
\endpicture} \end{center}
If we have two such which are composable then, after conjugating if needs be, we can assume that the start vertex of the second cell will be a vertex on the first one.  Firstly we will look at the case where the two start vertices coincide and there is a common edge containing it, schematically:
$$
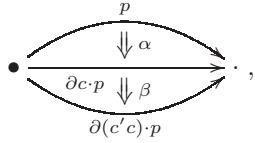
\xymatrix@+10pt{\bullet\rruppertwocell^p{~\alpha}\rrto_(.32){\partial c\cdot p}\rrlowertwocell_{\partial(c^\prime c)\cdot p}{~\beta}&&.},$$
where $\alpha = (c,p)$, and $\beta = (c^\prime,\partial c\cdot p)$. The obvious form for this `vertically'  composed $\mathcal{C}$-cell is:$$
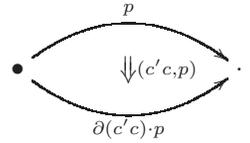
\xymatrix@+10pt{\bullet\rrtwocell<7>^p_{\partial(c^\prime c)\cdot p}{\hspace{6mm}(c^\prime c,p)} &&. }$$
This is justified from our earlier simplicial cases and a similar analysis of \emph{cellular} equivalence.  In other words, you build a `cylinder' over the first diagram with the lower diagram at its top, and then use the cocycle condition.

\

\textbf{Notation}

We will sometimes summarise this `vertical' composition as 
\begin{align*}
(c,p) \#_1 (c^\prime, \partial c \cdot p) = (c^\prime c,p).
\end{align*}
The 1 in $\#_1$ is there to indicate that the composite is formed across a shared 1-cell.

\

If the start vertex of the second $\mathcal{C}$-cell is another vertex of the first cell, we must use a `$\mathcal{C}$-path' from the first start vertex to the second.  There is a choice but it makes no difference.  Schematically we can reduce this to the case:
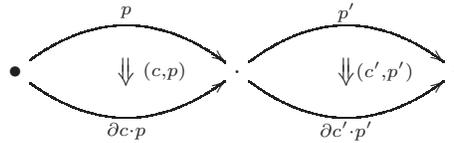
\begin{align*}\xymatrix@+10pt{\bullet\rrtwocell<7>^p_{\partial c\cdot  p}{\hspace{6mm}(c,p)} && . \rrtwocell<7>^{p^\prime}_{\partial c^\prime \cdot  p^\prime}{\hspace{6mm}(c^\prime,p^\prime)} &&.} \end{align*}
This happens, for instance, when combining two cylinders together, where both have start vertex on the inner circle or in the pair of pants cobordism, which can be represented as:
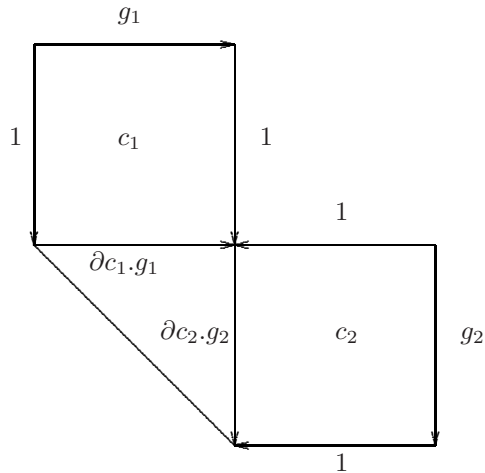
\begin{center}
\font\thinlinefont=cmr5
\begingroup\makeatletter\ifx\SetFigFont\undefined%
\gdef\SetFigFont#1#2#3#4#5{%
  \reset@font\fontsize{#1}{#2pt}%
  \fontfamily{#3}\fontseries{#4}\fontshape{#5}%
  \selectfont}%
\fi\endgroup%
\mbox{\beginpicture
\setcoordinatesystem units <0.70000cm,0.70000cm>
\unitlength=0.70000cm
\linethickness=1pt
\setplotsymbol ({\makebox(0,0)[l]{\tencirc\symbol{'160}}})
\setshadesymbol ({\thinlinefont .})
\setlinear
%
%
\put{$1$} [lB] at  9.525 10.001
%
%
\linethickness= 0.500pt
\setplotsymbol ({\thinlinefont .})
{\color[rgb]{0,0,0}\putrule from 11.430 14.287 to 11.430 10.478
%
%
\plot 11.366 10.731 11.430 10.478 11.494 10.731 /
}%
%
%
\linethickness= 0.500pt
\setplotsymbol ({\thinlinefont .})
{\color[rgb]{0,0,0}\putrule from  7.620 14.287 to  7.620 10.478
%
%
\plot  7.557 10.731  7.620 10.478  7.683 10.731 /
}%
%
%
\linethickness= 0.500pt
\setplotsymbol ({\thinlinefont .})
{\color[rgb]{0,0,0}\putrule from  7.620 18.098 to  7.620 14.287
%
%
\plot  7.557 14.541  7.620 14.287  7.683 14.541 /
}%
%
%
\linethickness= 0.500pt
\setplotsymbol ({\thinlinefont .})
{\color[rgb]{0,0,0}\putrule from  3.810 18.098 to  3.810 14.287
%
%
\plot  3.747 14.541  3.810 14.287  3.873 14.541 /
}%
%
%
\linethickness= 0.500pt
\setplotsymbol ({\thinlinefont .})
{\color[rgb]{0,0,0}\putrule from  3.810 18.098 to  7.620 18.098
%
%
\plot  7.366 18.034  7.620 18.098  7.366 18.161 /
}%
%
%
\linethickness= 0.500pt
\setplotsymbol ({\thinlinefont .})
{\color[rgb]{0,0,0}%
%
\plot  7.874 14.351  7.620 14.287  7.874 14.224 /
\putrule from  7.620 14.287 to 11.430 14.287
}%
%
%
\linethickness= 0.500pt
\setplotsymbol ({\thinlinefont .})
{\color[rgb]{0,0,0}%
%
\plot  7.874 10.541  7.620 10.478  7.874 10.414 /
\putrule from  7.620 10.478 to 11.430 10.478
}%
%
%
\linethickness= 0.500pt
\setplotsymbol ({\thinlinefont .})
{\color[rgb]{0,0,0}\plot  3.810 14.287  7.620 10.478 /
}%
%
%
\put{$c_2$} [lB] at  9.525 12,5
%
%
%
\put{$1$} [lB] at  3.334 16.192
%
%
\put{$1$} [lB] at  9.525 14.764
%
%
\put{$1$} [lB] at  8.096 16.192
%
%
\put{$g_1$} [lB] at  5.4 18.574
%
%
\put{$c_1$} [lB] at  5.4 16.192
%
%
%
\put{$g_2$} [lB] at 11.906 12.5
%
%
\put{$\partial c_2.g_2$} [lB] at  6.2 12.5
%
%
\put{$\partial c_1.g_1$} [lB] at  4.85 13.811
%
%
\linethickness= 0.500pt
\setplotsymbol ({\thinlinefont .})
{\color[rgb]{0,0,0}\putrule from  3.810 14.287 to  7.620 14.287
%
%
\plot  7.366 14.224  7.620 14.287  7.366 14.351 /
}%
\linethickness=0pt
\putrectangle corners at  3.334 19.025 and 12.266  9.874
\endpicture}
\end{center}
If we start the overall diagram at the top left corner, we have to move the top left (start) vertex of the $c_2$-cell back to that top left, acting on that 2-cell  with $g_1$ in the process.  

\

In the general picture, if we move the start vertex of the second cell back along the upper 1-cell we get.
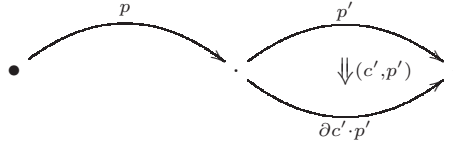
\begin{align*}\xymatrix@+10pt{\bullet\ar@/^1.5pc/[rr]^p & &. \rrtwocell<7>^{p^\prime}_{\partial c^\prime \cdot p^\prime} {\hspace{6mm}(c^\prime,p^\prime)} && .} \end{align*}
and clearly this is be the $\mathcal{C}$-cell
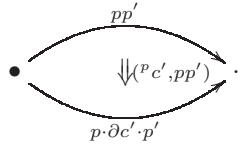
\begin{align*}\xymatrix@+10pt{\bullet\rrtwocell<7>^{pp^\prime}_{p\cdot\partial c^\prime \cdot p^\prime}{\hspace{8mm}({}^pc^\prime,pp^\prime)} &&. }\end{align*}
The other half of our data can be fitted to the bottom of this by `whiskering' on the right. This just shifts the second vertex of our first `glob' to that of the second and adds in a cancellable sub-path:
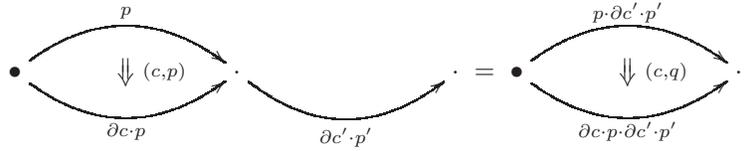
\begin{align*}\xymatrix@+10pt{\bullet\rrtwocell<7>^p_{\partial c\cdot  p}{\hspace{6mm}(c,p)} && .\ar@/_1.5pc/[rr]_{\partial c^\prime \cdot  p^\prime} & & .} = \xymatrix@+10pt{\bullet\rrtwocell<7>^{p\cdot\partial c^\prime \cdot  p^\prime}_{\partial c\cdot p\cdot\partial c^\prime \cdot p^\prime}{\hspace{6mm}(c,q)} &&. },\end{align*}
where $q = p. \partial c^\prime \cdot  p^\prime$. These $\mathcal{C}$-cells now have matching edges and so can be composed `vertically' to get a $\mathcal{C}$-cell labelled $( c\cdot({}^pc^\prime),pp^\prime).$
One might object that  this combination or composition algorithm looks as if it depends on choices being made and the obvious key choice here was the way we whiskered the cells starting from our initial data. We could equally well have decided to decompose this as:
$$
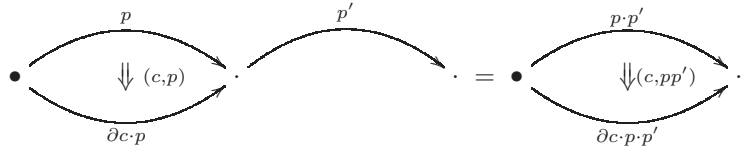
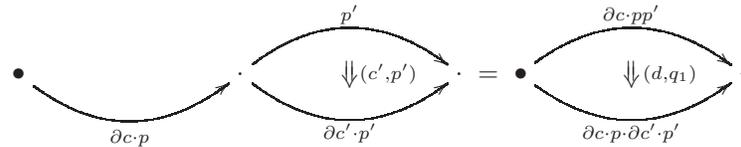
\xymatrix@+10pt{\bullet\rrtwocell<7>^p_{\partial c\cdot  p}{\hspace{6mm}(c,p)} && .\ar@/^1.5pc/[rr]^{ p^\prime} & & .} = \xymatrix@+10pt{\bullet\rrtwocell<7>^{p \cdot  p^\prime}_{\partial c\cdot p \cdot p^\prime}{\hspace{6mm}(c,pp^\prime)} &&. },$$
and
$$\xymatrix@+10pt{\bullet\ar@/_1.5pc/[rr]_{\partial c\cdot p} & &. \rrtwocell<7>^{p^\prime}_{\partial c^\prime \cdot p^\prime} {\hspace{6mm}(c^\prime,p^\prime)} && .} = \xymatrix@+10pt{\bullet\rrtwocell<7>^{\partial c\cdot pp^\prime}_{\partial c\cdot p\cdot \partial c^\prime\cdot p^\prime}{\hspace{6mm}(d,q_1)}&&.}$$
with $(d,q_1) = ({}^{\partial c\cdot p}c^\prime,{\partial c\cdot pp^\prime})$. Now we should form a `vertical' composite
of the two cells.  This will give us an apparently different composite $\mathcal{C}$-cell, this time labelled by $({}^{\partial c\cdot p}c^\prime\cdot c, {pp^\prime}).$
These two $\mathcal{C}$-cobordisms have the same upper and, less obviously, lower parts, so we need to compare the two $C$-parts.  This uses the Peiffer identity, i.e. the second axiom for crossed modules:
$
{}^{\partial c\cdot p}c^\prime\cdot c = {}^{\partial c}({}^pc^\prime)\cdot c = c\cdot ({}^pc^\prime)\cdot c^{-1}\cdot c = c\cdot ({}^pc^\prime).$
The two composite $\mathcal{C}$-cells are thus the same.

Clearly there is something going on here that has not been revealed in detail.  Readers who know some 2-category theory will have noticed that the above is a manifestation of the `interchange law' of the theory of 2-categories.  The Peiffer identity is an instance of that law.  (This corresponds closely  to the monoidal category structure on the category of cobordisms.)  Consistently with our previous notation, it may be useful to use $\#_0$ for this second `horizontal' composition, so 
\begin{align*}
(c,p)\#_0 (c^\prime,p^\prime) = (c\cdot ({}^pc^\prime),pp^\prime).
\end{align*}
The 0 in $\#_0$ indicates that the composite can be formed because of a shared 0-cell.


\subsection{Proof of Main Theorem}

We start by identifying the geometric behaviour of the isomorphisms $\theta_{(c,g)}.$ We know that, from the special case of $C=1$,  $L$ is a crossed $P$-algebra, so we need to look at the extra structure:
\begin{itemize}
\item \textbf{Influence of composition of $\mathcal{C}$-cobordism fragments.}

\emph{`Vertical' compositions.}

The $\mathcal{C}$ structure will reflect the composition of such $\mathcal{C}$-fragments.  Firstly we handle $\#_1$, the vertical composition:
$$(c,p)\#_1(c^\prime, \partial c,p) = (c^\prime c, p).$$
The basic condition is thus that the composite
$$L_g \stackrel{\theta_{(c,g)}}{\to}L_{\partial c\cdot g}\stackrel{\theta_{(c^\prime,\partial c\cdot g)}}{\to}L_{\partial(c^\prime c)\cdot g}$$
 is $\theta_{(c^\prime c,g)}$:
\begin{align*}\fbox{\quad$\theta_{(c^\prime c,g)} = \theta_{(c^\prime,\partial c\cdot g)}\circ \theta_{(c,g)} : L_g \to L_{\partial(c^\prime c)\cdot g}
$\quad}\end{align*}
since $\tau$ must be compatible with the `vertical composition' of $\mathcal{C}$-fragments.  Evaluating this on an element gives$$\widetilde{(c^\prime c)} = \tilde{c^\prime}\cdot \tilde{c},$$
where $\tilde{c^\prime}\cdot \tilde{c} = \mu(\tilde{c^\prime}, \tilde{c})$.  Similarly $\tilde{1} = 1$.

\emph{`Horizontal' composition  of $\mathcal{C}$-fragments.}

Using the interchange law / Peiffer rule or, equivalently, the `pair of pants' to give the multiplication, we thus have
$$(1,g_1)\#_0(c,g_2)=  ({}^{g_1}c,g_1g_2) = ({}^{g_1}c,g_1)\#_0(1,g_2).$$
(Here the useful notation $\#_0$ corresponds to the horizontal composition in the associated strict 2-group of $\mathcal{C}$.)   We thus have two composite $\mathcal{C}$-cobordisms giving the same result and hence
\begin{align*}
L_{g_1}\otimes L_{g_2}\stackrel{L_{g_1}\otimes \theta_{(c, g_2)}}{\longrightarrow}& L_{g_1}\otimes L_{\partial c \cdot g_2}\stackrel{\mu}{\longrightarrow} L_{g_1\cdot\partial c \cdot g_2}\\
& =L_{g_1}\otimes L_{g_2}\stackrel{\mu}{\longrightarrow} L_{g_1 g_2}\stackrel{\theta_{({}^{g_1}c, g_2)}}{\longrightarrow}L_{g_1\cdot\partial c \cdot g_2}.
 \end{align*}
In general for $d \in C$, the second type of composite will be 
\begin{align*}
L_{g_1}\otimes L_{g_2}\stackrel{\theta_{(d,g_1)} \otimes L_{g_2}}{\longrightarrow} &L_{\partial d\cdot g_1} \otimes L_{g_2} \stackrel{\mu}{\longrightarrow} L_{\partial d \cdot g_1\cdot g_2}\\
& =L_{g_1}\otimes L_{g_2}\stackrel{\mu}{\longrightarrow} L_{g_1 g_2}\stackrel{\theta_{(d,g_1g_2)}}{\longrightarrow}L_{\partial d \cdot g_1\cdot g_2},
 \end{align*}
and we need this for $d = {}^{g_1}c$ for which the corresponding composite cobordisms are equal.
Geometrically these rules correspond to a pair of pants with $g_1$, $g_2$ on the trouser cuffs and the 2-cell colored $c$.  We can push $c$ onto either leg, but in so doing may have to conjugate by $g_1$, somewhat as in Figure \ref{pantsreduction} .

Summarising, for given $c\in C$, $g_1,g_2\in G$,
\begin{align*}\fbox{\quad$\mu ( id_{L_{g_1}}\otimes \theta_{(c,g_2)}) = \theta_{({}^{g_1}c,g_1g_2)} \circ \mu = \mu ( \theta_{({}^{g_1}c,g_1)} \otimes id_{L_{g2}}).
$}\end{align*}

As we have reduced $\#_0$ to `whiskering' and the vertical composition, $\#_1$, and have already checked the interpretation of $\#_1$, we might expect this pair of equations to follow from our earlier calculations, however we have invoked here the interchange law and that was not used earlier.  The above equations reduce, and simplify,  to give
\begin{align*}x\cdot \tilde{c} = \widetilde{{}^gc}\cdot x,\end{align*}
but this is implied by axiom c) of a crossed $\mathcal{C}$-algebra, since we have  
$$\widetilde{{}^gc}\cdot x = \phi_g(\tilde{c})x = x\cdot \tilde{c},$$
using the third axiom (page \pageref{Palg3}) of the crossed $P$-algebra structure on $L$. Thus the combination of these two rules corresponds in part to the Interchange Law. Conversely this rule in either form  is clearly implied by the axioms for a formal HQFT. 

We still have to check that  the inner product structure of $L$ and action of $P$ via $\phi$ are compatible with the new structure. The compatibility of the isomorphisms $\theta_{(c,g)}$ defined via the $\tilde{c}$ will  follow, both from the geometry of the HQFT and from the axioms of  crossed $\mathcal{C}$-algebras.

\item \textbf{Inner product.}

The inner product $$\rho : L\otimes L \to \mathbb{K}$$
restricts, for any $g \in P$, to $$L_g\otimes L_{g^{-1}} \to \mathbb{K}.$$ Now consider the two possible composite pairings
$$L_g\otimes L_{(\partial c\cdot g)^{-1}} \stackrel{\theta_{(c,g)}\otimes  L_{(\partial c\cdot g)^{-1}}}{\longrightarrow}L_{\partial c\cdot g}\otimes L_{(\partial c\cdot g)^{-1}}\stackrel{\rho}{\to} \mathbb{K}$$
and the alternative
$$L_g\otimes L_{(\partial c\cdot g)^{-1}} \stackrel{L_g\otimes \theta_{({}^{g^{-1}}c,g^{-1}\partial c^{-1})}}{\longrightarrow}L_g\otimes L_{g^{-1}}\stackrel{\rho}{\to} \mathbb{K}.$$
These correspond to two composite $\mathcal{C}$-cobordisms that are equivalent, as is clear from the geometry.  They thus imply an equality
\begin{align*}\fbox{\quad $\rho(\tilde{c}\cdot x,y) = \rho(x, \widetilde{{}^{g^{-1}}c}\cdot y)$}
\end{align*}
for $x\in L_g$, $y\in L_{(\partial c\cdot g)^{-1}}$.  We need to check that the inner product property follows from the axioms of a crossed $\mathcal{C}$-algebra.

From the third axiom for the $\tilde{c}$s, we get
$$\tilde{c}=\widetilde{{}^{gg^{-1}}c}= \phi_{g}(\widetilde{{}^{g^{-1}}c}),$$
but then $$\tilde{c}\cdot x = \phi_g(\widetilde{{}^{g^{-1}}c})\cdot x = x\cdot\widetilde{ {}^{g^{-1}}c}$$and 
$$\rho(\tilde{c}\cdot x,y) = \rho(x\cdot\widetilde{ {}^{g^{-1}}c},y) = \rho(x,\widetilde{{}^{g^{-1}}c}\cdot y),$$
as required.

 \item \textbf{$P$-action via $\phi$.}

This, geometrically, is clearly the 3rd condition on `tilderisation'
$$\phi_h(\tilde{c}) = \widetilde{{}^hc}.$$
Composing a formal $\mathcal{C}$-disc, $Disc(c)$, with a cylinder $(Cyl_{-,+};1,\partial c,h)$ is equivalent to $Disc({}^{h^{-1}}c)$:
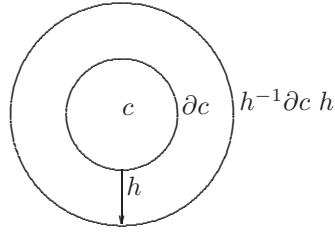
\begin{center}
\font\thinlinefont=cmr5
\begingroup\makeatletter\ifx\SetFigFont\undefined
\def\x#1#2#3#4#5#6#7\relax{\def\x{#1#2#3#4#5#6}}%
\expandafter\x\fmtname xxxxxx\relax \def\y{splain}%
\ifx\x\y   
\gdef\SetFigFont#1#2#3{%
  \ifnum #1<17\tiny\else \ifnum #1<20\small\else
  \ifnum #1<24\normalsize\else \ifnum #1<29\large\else
  \ifnum #1<34\Large\else \ifnum #1<41\LARGE\else
     \huge\fi\fi\fi\fi\fi\fi
  \csname #3\endcsname}%
\else
\gdef\SetFigFont#1#2#3{\begingroup
  \count@#1\relax \ifnum 25<\count@\count@25\fi
  \def\x{\endgroup\@setsize\SetFigFont{#2pt}}%
  \expandafter\x
    \csname \romannumeral\the\count@ pt\expandafter\endcsname
    \csname @\romannumeral\the\count@ pt\endcsname
  \csname #3\endcsname}%
\fi
\fi\endgroup
\mbox{\beginpicture
\setcoordinatesystem units <0.50000cm,0.50000cm>
\unitlength=0.50000cm
\linethickness=1pt
\setplotsymbol ({\makebox(0,0)[l]{\tencirc\symbol{'160}}})
\setshadesymbol ({\thinlinefont .})
\setlinear
%
%
\linethickness= 0.500pt
\setplotsymbol ({\thinlinefont .})
\ellipticalarc axes ratio  2.963:2.963  360 degrees 
	from  8.970 21.643 center at  6.007 21.643
%
%
\linethickness= 0.500pt
\setplotsymbol ({\thinlinefont .})
\putrule from  5.997 20.180 to  5.997 18.735
%
%
\plot  5.944 18.946  5.997 18.735  6.049 18.946 /
\put{$c$} [lB] at  6.007 21.643

%
%
%
\put{$\partial c$} [lB] at  7.6 21.660
%
%
\put{$h$} [lB] at  6.128 19.490
%
%
\put{$h^{-1}\partial c~ h$} [lB] at  9.129 21.749
%
%
\linethickness= 0.500pt
\setplotsymbol ({\thinlinefont .})
\ellipticalarc axes ratio  1.482:1.482  360 degrees 
	from  7.489 21.643 center at  6.007 21.643
\linethickness=0pt
\putrectangle corners at  3.027 24.621 and  9.129 18.665
\endpicture} \end{center}\end{itemize}
The formal details of the reconstruction of $\tau$ from $L$ follow the same pattern as for the case $C = 1$ and, for the most part, are exactly the same, as the only extra feature is the `tilde' operation.  The details are not hard and left to the reader.\hfill$\square$

\

\textbf{Remark}

It is sometimes useful to have the extra rules of the $\tilde{c}$s written in the intermediate language of the family of isomorphisms
$$\theta_{(c,g)} : L_g \to L_{\partial c\cdot g}.$$
The first two conditions are easily so interpreted and the last corresponds to the compositions given earlier and also to
the equality of $$L_g \stackrel{\theta_{(c,g)}}{\to} L_{\partial c \cdot g} \stackrel{\phi_h}{\to} L_{h\cdot\partial c \cdot g\cdot h^{-1}},$$
and$$L_g \stackrel{\phi_h}{\to} L_{hgh^{-1}}\stackrel{\theta_{({}^hc,{}^hg)}}{\to}L_{h\partial c \cdot g\cdot h^{-1}},$$and  thus to 
 \begin{align*}\fbox{\quad$\phi_h\circ \theta_{(c,g)} = \theta_{({}^hc,{}^hg)}\circ \phi_h.$ \quad }\end{align*}
Collectively these boxed equations in their various forms give compatibility conditions for the various structures.  They help express the extra structure coming from the non-trivial `2-cells' in an algebraic form.

\medskip

There is a very neat interpretation of these conditions.  Let $L$ be an associative algebra and $U(L)$ be its group of units. There is a homomorphism of groups $\delta = \delta_L : U(L) \to Aut(L)$ given by $\delta(u) (x) = u\cdot x\cdot u^{-1}$.
\begin{lemma}\label{AutL} \hspace*{1mm}\\
With the obvious action of $Aut(L)$ on the group of units, $(U(L),Aut(L),\delta)$ is a crossed module.\hspace*{.1mm}
\hfill$\square$\end{lemma}
The proof is simple, although quite instructive, and will be left to the reader. We will denote this crossed module by $\mathfrak{Aut}(L)$. If $L$ has extra structure such as being a Frobenius algebra or being graded, the result generalises to have the automorphisms respecting that structure.
\begin{proposition}\hspace*{1mm}\\
Suppose that $L$ is a crossed $\mathcal{C}$-algebra.
The diagram
\begin{align*}\xymatrix{C\ar[d]_\partial \ar[r]^{(\tilde{~})}& U(L)\ar[d]^\delta\\
P\ar[r]_\phi &Aut(L) }\end{align*}is a morphism of crossed modules from $\mathcal{C}$ to  $\mathfrak{Aut}(L)$.
\end{proposition}
\textbf{Proof}

First we check commutativity of the square in the statement of the proposition.  Let $c\in C$, going around clockwise gives $\delta(\tilde{c})$ and on an element $x\in L$, this gives $\tilde{c}\cdot x \cdot \tilde{c}^{-1}$.  We compare this with
the other composite, again acting on $x\in L$. If we multiply $\phi_{\partial c}(x)$ by $\tilde{c}$, then we get $\phi_{\partial c}(x)\tilde{c} = \tilde{c}\cdot x$,  but therefore  $\phi_{\partial c}(x) = \tilde{c}\cdot x \cdot \tilde{c}^{-1}$ as well.

The other thing to check is that the maps are compatible with the actions of the bottom groups on the top ones, but this is exactly what the third condition on the `tilde' gives.\hfill$\square$
\section{Constructions on formal HQFTs and  crossed $\mathcal{C}$-algebras}
As formal HQFTs correspond to crossed $\mathcal{C}$-algebras by our main result above, the category of crossed $\mathcal{C}$-algebras needs to be understood better if we are to understand the relationships between formal HQFTs. We clearly also need some examples of crossed $\mathcal{C}$-algebras. 

First we note that the usual constructions of direct sum and tensor product of graded algebras extends to crossed $\mathcal{C}$-algebras in the obvious way.
\subsection{Examples of crossed $\mathcal{C}$-algebras}
As usual we will fix a crossed module $\mathcal{C} = (C,P,\partial)$. We assume, for convenience, that $\ker \partial $ is a finite group, although this may not always be strictly necessary.

\

\textbf{The group algebra, $\mathbb{K}[C]$, as a crossed $\mathcal{C}$-algebra.}

We take $L = \mathbb{K}[C]$ and will denote the generator corresponding to $c\in C$ by $e_c$ rather than merely using the symbol $c$ itself, as we will need a fair amount of precision when specifying various types of related elements in different settings.  Define $L_p = \mathbb{K}\langle\{ e_c : \partial c = p\}\rangle$, so, if $p\in P\setminus \partial C$, this is the zero dimensional $\mathbb{K}$-vector space, otherwise it has dimension the order of $\ker \partial$ (whence our requirement that this be finite).
\begin{lemma}\hspace*{1mm}\\
With this grading structure, $L$ is a crossed $P$-algebra.
\end{lemma}
\textbf{Proof}
\begin{itemize}
\item $L$ is $P$-graded: this follows since \quad $e_c \cdot e_{c^\prime} = e_{cc^\prime}$,  \quad $\partial$ is a group homomorphism and $e_1\in L_1$.
\item There is an inner product:
$$\rho : L \otimes L \to \mathbb{K}$$
$$\rho(e_c\otimes e_{c^\prime}) = \left\{\begin{array}{ll}
0 \quad& \quad \textrm{if } c^{-1} \neq c^\prime\\
1 & \quad\textrm{otherwise}\end{array}\right.$$
and this is clearly non-degenerate.  Moreover $$\rho(e_{c_1}e_{c_2}\otimes e_{c_3}) = \rho(e_{c_1c_2}\otimes e_{c_3}) = 0$$unless $c_3 = c_2^{-1}c_1^{-1}$ when it is 1, whilst
$$\rho(e_{c_1}\otimes e_{c_2c_3}) = 0$$
unless $c_1^{-1} = c_2c_3$, etc., so the inner product satisfies the third condition for a Frobenius $P$-algebra.
\item Finally there is a group homomorphism
$$\phi : P \to Aut(L),$$
given by $\phi_g(e_c) = e_{{}^gc}$, which permutes the basis, compatibly with the multiplication and innerproduct structures.
\end{itemize}
As $\partial({}^gc) = g\cdot \partial c\cdot  g^{-1}$, $\phi$ clearly satisfies $\phi_g(L_h) \subseteq L_{ghg^{-1}}$, and the Peiffer identity implies ${}^{\partial c}c = c$, so $\phi_g | L_g$ is the identity.  The Peiffer identity in general gives
$${}^{\partial c}c^\prime = cc^\prime c^{-1},$$
so $e_ce_{c^\prime} = e_{({}^{\partial c}c^{\prime})}e_c$, i.e., $\phi_h(a)b = ba$ if $a\in L_g, b\in L_h$.\hfill$\square$
\medskip

As we want this to be a crossed $\mathcal{C}$-algebra, the remaining structure we have to specify is the `tildefication'
$$\tilde{\,} : C \to \mathbb{K}[C].$$
The obvious mapping gives $\tilde{c} = e_c$, and, of course,
$$\delta(\tilde{c})(e_{c^\prime}) = e_ce_{c^\prime}e_{c^{-1}} = e_{cc^\prime c^{-1}} = \phi_{\partial c}(e_{c^\prime}),$$
as above. We thus have
\begin{proposition}\hspace*{1mm}\\
With the above structure, $\mathbb{K}[C]$ is a crossed $\mathcal{C}$-algebra.\hfill $\square$
\end{proposition}
By its construction $\mathbb{K}[C]$ records little of the structure of $P$ itself, only the way the $P$-action permutes the elements of $C$, but, of course, it records $C$ faithfully.  The next example give another extreme.

\vspace{5mm}

\textbf{The group algebra $\mathbb{K}[P]$ as a crossed $\mathcal{C}$-algebra.}

We first note the following result from \cite{turaev:hqft1}:
\begin{lemma}\hspace*{1mm}\\
$\mathbb{K}[P]$ has the structure of a crossed $P$-algebra with $(\mathbb{K}[P])_p = \mathbb{K}e_p$, the subspace generated by the basis element labelled by $p\in P$.\hfill$\square$
                            \end{lemma}
The one thing to note is that the axiom
$$\phi_h(a)b = ba$$
for any $g,h \in P$, $a \in L_g,b\in L_h$ implies that 
$$\phi_h(e_g) = e_h e_ge_{h^{-1}} = e_{hgh^{-1}},$$
since $e_h$ is a unit of $\mathbb{K}[P]$ with inverse $e_{h^{-1}}$.
\begin{proposition}\hspace*{1mm}\\
For $c\in C$, defining $\tilde{c} = e_{\partial c}$, gives $\mathbb{K}[P]$ the additional structure of a crossed $\mathcal{C}$-algebra.
\end{proposition}
\textbf{Proof}

The grading is as expected and $\delta(\tilde{c}) = \phi_{\partial c}$, by construction.\hfill$\square$

\

Of course, $\mathbb{K}[P]$ does not encode anything about the kernel of $\partial : C\to P$.  In fact, it basically remains a crossed $P$-algebra as the extra crossed $\mathcal{C}$-structure is derived from that underlying algebra. 

We will give further examples of crossed $\mathcal{C}$-algebras shortly.

\subsection{Morphisms of crossed algebras}
We clearly need to have a notion of morphism of crossed $\mathcal{C}$-algebras. We start with a fixed crossed module $\mathcal{C} = (C,P,\partial)$.

\textbf{Definition}

Suppose $L$ and $L^\prime$ are two crossed $\mathcal{C}$-algebras.  A $\mathbb{K}$-algebra morphism $\theta : L\to L^\prime$ is a \emph{morphism of crossed $\mathcal{C}$-algebras} if it is compatible with the extra structure.  Explicitly:
\begin{align*}
\theta(L_p)&\subseteq L_p^\prime\\
\rho^\prime(\theta a,\theta b) &= \rho(a,b),\\
\phi^\prime_h(\theta a) &= \theta(\phi_h (a)),\\
\theta(\tilde{c}) &= \tilde{c}
\end{align*}
for all $a,b \in L$, $h\in P$, $c\in C$, where primes indicate the structure in $L^\prime$. 

\medskip

We know that a given crossed module represents a homotopy 2-type, but that different crossed modules can give equivalent 2-types, so it will also be necessary to compare crossed algebras over different crossed modules. We need this not just to move within a 2-type, but for various constructions linking different 2-types. We therefore put forward the following definition.  First some preliminary notation:\\
Suppose $f: \mathcal{C}\to \mathcal{D}$ is a morphism of crossed modules.  The morphism $f$ gives a commutative square of group homomorphisms
$$\xymatrix{C\ar[r]^{f_1}\ar[d]_\partial &D\ar[d]^{\partial^\prime}\\
P\ar[r]_{f_0}&Q.}$$
We want to define a morphism of crossed algebras over $f$, i.e., an algebra morphism, $\theta : L\to L^\prime$, where $L$ is a crossed $\mathcal{C}$-algebra and $L^\prime$, a crossed $\mathcal{D}$-algebra.

\

\textbf{Definition}

Suppose $L$ and $L^\prime$ are two crossed algebras over $\mathcal{C}$ and $\mathcal{D}$, respectively.  A $\mathbb{K}$-algebra morphism $\theta : L\to L^\prime$ is a \emph{morphism of crossed algebras over $f$} if it is compatible with the extra structure.  Explicitly:
\begin{align*}
\theta(L_p)&\subseteq L_{f_0(p)}^\prime\\
\rho^\prime(\theta a,\theta b) &= \rho(a,b),\\
\phi^\prime_{f_0(h)}(\theta a) &= \theta(\phi_h (a)),\\
\theta(\tilde{c}) &= \widetilde{f_1(c)}
\end{align*}
for all $a,b \in L$, $h\in P$, $c\in C$, where primes indicate the structure in $L^\prime$.

\subsection{Pulling back a crossed $\mathcal{C}$-algebra}

A morphism, as above, over $f$ can be replaced by a morphism of crossed $\mathcal{C}$-algebras, $L\to  f_0^*(L^\prime)$, where $f_0^*(L^\prime)$ is obtained by pulling back $L^\prime$ along $f$.  We will consider this construction independently of any particular $\theta$.

If $f_0 : P\to Q$ is a group homomorphism, we know, from \cite{turaev:hqft1} that given a crossed $Q$-algebra, $L$, we obtain a crossed $P$-algebra $f_0^*(L)$, by pulling back using $f_0$. The structure of $f_0^*(L)$ is given by:\\
\begin{itemize}\item $(f_0^*(L))_p$ is $L_{f_0(p)}$,
by which we mean that $(f_0^*(L))_p$ is a copy of $L_{f_0(p)}$ with grade $p$ and we note that if $x\in L_{f_0(p)}$, it can  be useful to write it  $x_{f_0(p)}$ with $x_p$ denoting the corresponding element of $(f_0^*(L))_p$;
\item if $x$ and $y$ have inverse grades, say $x \in (f_0^*(L))_p$, $y \in (f_0^*(L))_{p^{-1}}$, then $\rho(x,y)$ is the same as in $L$, but if $x$ and $y$ have non-inverse grades then $\rho(x,y) = 0$;
\item $\phi_h (x_p) := \phi_{f_0(h)}(x_{f_0(p)})$.
\end{itemize}
If, in addition, we consider the crossed $\mathcal{C}$-structure assuming that $L^\prime$ is a crossed $\mathcal{D}$-algebra, then defining $\tilde{c}:= \widetilde{f_1(c)}_{\partial c}$ gives us:
\begin{proposition}\hspace*{1mm}\\
The crossed $P$-algebra $f_0^*(L)$ has a crossed $\mathcal{C}$-algebra structure given by the above.\hfill$\square$
\end{proposition}
The construction of $f_0^*(L)$, then, makes it clear that
\begin{proposition}\hspace*{1mm}\\
There is a bijection between the set of crossed algebra morphisms from $L$ to $L^\prime$ over $f$ and the set of crossed $\mathcal{C}$-algebra morphisms from $L$ to $f_0^*(L^\prime)$.\hfill$\square$
\end{proposition}
Of course, as with most such operations, this pullback construction gives a functor from the category of crossed $\mathcal{D}$-algebras to that of crossed $\mathcal{C}$-algebras (up to isomorphism in the usual way).
\subsection{Applications of pulling back}

Consider our crossed module $\mathcal{C} = (C,P,\partial)$ and let $G = P/\partial C$.  We can realise this as a morphism of crossed modules:
$$\xymatrix{C\ar[r]\ar[d]_\partial &1\ar[d]\\
P\ar[r]_{q}&G.}$$
If $\partial$ was an inclusion then this would be a weak equivalence of crossed modules as then both the kernel and cokernels of the crossed modules would be mapped isomorphically by the induced maps. In that case, thinking back to our original motivations for introducing formal $\mathcal{C}$-maps, we would really be in a situation corresponding to a HQFT with background a $K(G,1)$ and by \cite{turaev:hqft1}, we know such theories are classified by crossed $G$-algebras. Thus it is of interest to see what the pullback algebra of a crossed $G$-algebra along this morphism will be. We will look at the obvious example of $\mathbb{K}[G]$, the group algebra of $G$ with its usual crossed $G$-algebra structure (cf., \cite{turaev:hqft1}). We will assume that the crossed module, $\mathcal{C}$, is finite.

Writing $N = \partial C$, for convenience, we have an extension
$$\xymatrix{N\ar[r] & P \ar[r]^q&G}.$$
Pick a section $s$ for $q$ and define the corresponding cocycle $f(g,h) = s(g)s(h)s(gh)^{-1}$, so $f : G\times G\to N$ is naturally normalised, $f(1,h)= f(g,1) = 1$ and satisfies the cocycle condition:
\begin{align*}f(g,h)f(gh,k) = {}^{s(g)}f(h,k)f(g,hk).\end{align*}
Take $L = \mathbb{K}[G]$, the group algebra of $G$ considered with its crossed $G$-algebra structure and form the crossed $P$-algebra, $q^*(L)$. We will give a cohomological proof of the following to illustrate some of the links between cohomology and constructions on crossed algebras.
\begin{proposition}\hspace*{1mm}\\
The two crossed $\mathcal{C}$-algebras $\mathbb{K}[P]$ and $q^*(\mathbb{K}[G])$ are isomorphic.
\end{proposition}
\textbf{Proof}

We first note that $$q^*(L)_p = L_{q(p)} = \mathbb{K}e_{q(p)}.$$ We will write $g = q(p)$, so $p\in P$ has the form $p = n s(g)$. (We will need to keep check of which $e_{q(p)}$ is which and will later introduce notation which will handle this.)

Recall the description of the product in $P$ in terms of the cocycle and the section:
\begin{eqnarray}
n_1s(g_1)\cdot n_2s(g_2) &=& n_1{}^{s(g_1)}n_2 s(g_1)s(g_2)\\
			       &=& (n_1{}^{s(g_1)}n_2 f(g_1,g_2))s(g_1g_2).
\end{eqnarray}
Each unit, $e_g$, of $\mathbb{K}[G]$ gives $\#(N)$ copies in $q^*(L)$.  Write $(e_g)_n$ for the copy of $e_g$ in $q^*(L)_{ns(g)}$ and examine the multiplication in $q^*(L)$ in this notation:
$$(e_{g_1})_{n_1}\cdot (e_{g_2})_{n_2} = (e_{g_1g_2})_{(n_1{}^{s(g_1)}n_2 f(g_1,g_2))}.$$
(That this gives an associative multiplication corresponds to the cocycle condition above.)

We next have to ask : what is $\phi_p$? Of course as $p = ns(g)$, we can restrict to examining $\phi_n$ and $\phi_{s(g)}$.\begin{itemize}
\item $\phi_n$ links the two copies $ q^*(L)_p$ and $q^*(L)_{npn^{-1}}$ of $L_{q(p)}$ via what is essentially the identity map between the two copies;
\item $\phi_{s(g)}$ restricts to $\phi_{s(g)} : q^*(L)_p \to q^*(L)_{s(g)ps(g)^{-1}}$, but on identifying these two subspaces as $L_{q(p)}$ and $L_{gq(p)g^{-1}}$, this is just $\phi_g$.
\end{itemize}
In fact we can be more explicit if we look at the basic  units and, as these do form a basis, behaviour on them determines the automorphisms:
$$\phi_n((e_{g_1})_{n_1}) (e_1)_n = (e_1)_n(e_{g_1})_{n_1},$$
so 
\begin{eqnarray}
\phi_n((e_{g_1})_{n_1}) &=& (e_{g_1})_{nn_1}(e_1)_n^{-1} \\
			& = & (e_{g_1})_{nn_1}(e_1)_{n^{-1}}\\
			& = & (e_{g_1})_{nn_1{}^{s(g_1)}n^{-1}}
\end{eqnarray}
that is, conjugation by $(e_1)_n$.

This leads naturally on to noting that $\tilde{c} = (e_1)_{\partial c}$, so we have explicitly given  the crossed $\mathcal{C}$-algebra structure on $q^*(L)$. Sending $e_p$ to $(e_{q(p)})_n$ (using the same notation as before) establishes the isomorphism of the statement without difficulty. \hfill $\square$

\medskip

\textbf{Remark}

In this identification of $q^*(\mathbb{K}[G])$ as $\mathbb{K}[P]$, it is worth noting that $$q^*(L)_1 = L_1 = \mathbb{K}e_1 \cong \mathbb{K},$$ as a vector space, but also that $q^*(L)_n \cong \mathbb{K}$ for each $n \in N$.  The notation $(e_g)_n$ used and the behaviour of these basis elements suggests that $q^*(\mathbb{K}(P))$ behaves like some sort of  twisted tensor product with basis $e_n \otimes e_g$, with that element corresponding to $(e_g)_n$, and with multiplication
$$(e_{n_1} \otimes e_{g_1})(e_{n_2} \otimes e_{g_2}) = (e_{n_1{}^{s(g_1)}n_2 f(g_1,g_2)} \otimes e_{g_1g_2}).$$
We have not yet investigated how general this construction may be.

\subsection{Pushing forward}
We have shown that, given $f : \mathcal{C}\to \mathcal{D}$ and a $\theta : L \to L^\prime$ over $f$, we can pull $L^\prime$ back over $\mathcal{C}$ to get a map from $L$ to $f^*(L^\prime)$ that encodes the same information as $L^\prime$ (provided $f$ is an epimorphism and all crossed modules are finite).  An obvious question to ask is whether there is an `adjoint' push-forward construction with $\theta$ corresponding to some morphism from $f_*(L)$ to $L^\prime$ over $\mathcal{D}$.  This is what we turn to next, keeping the same assumptions of finiteness, etc.

\medskip

Given such a context, setting, as before, $N = \ker f_0$, $B = \ker f_1$, we have
\begin{align*}\fbox{\quad$\phi_n(a) - a \in \ker \theta,
$\quad}\end{align*}
as $\theta(\phi_n(a)) = \phi^\prime_{f_0(n)}(\theta (a)) = \phi^\prime_1(\theta(a))= \theta(a)$.  Similarly, since $$\theta (\tilde{c}) =\widetilde{f_1(c)},$$if $b \in B = \ker f_1$, $$\theta(\tilde{b}) = \tilde{1},$$ so 
\begin{align*}\fbox{\quad$\tilde{b}-1\in \ker \theta .$\quad}\end{align*}
We therefore form the ideal $\mathcal{K}$ generated by elements of these forms, above.  Note this is not a $P$-graded ideal, but, in fact, that is exactly what is needed.  We have that $L/\mathcal{K}$  is an associative algebra and we give it a $Q$-graded algebra structure as follows.

For each $q \in Q$, let $$\mathcal{L}_q = \oplus_p\{L_p ~|~ f_0(p)=q\},$$
and $$\mathcal{K}_q = \mathcal{L}_q \cap \mathcal{K}.$$  The underlying $Q$-graded vector space of $f_*(L)$ will be 
$$f_*(L) = \oplus_{q \in Q} \mathcal{L}_q/\mathcal{K}_q.$$
This is an associative algebra as it is exactly $L/\mathcal{K}$, but we have to check that this grading is compatible with that multiplication.

Suppose $a + \mathcal{K}\in f_*(L)_{q_1}$, and $b + \mathcal{K}\in f_*(L)_{q_2}$, then $a\in L_{p_1}$ and $b\in L_{p_2}$ for some $p_1,p_2 \in P$ with $f_0(p_i) = q_i$, for $i = 1,2$, but then $ab + \mathcal{K} \in f_*(L)_{q_1q_2}$ as required.

We next define the bilinear form giving the inner product. Clearly, with the same notation,
\begin{align*}\rho(a+\mathcal{K},b+\mathcal{K}) := 0\textrm{ \quad if }q_1\neq q_2^{-1}.\end{align*}
If $q_1 = q_2^{-1}$, then we can assume that $p_1 = p_2^{-1}$, and,  if necessary, after changing the element $b$ representing $b + \mathcal{K}$ that $b \in L_{p_2}$.  Finally we set 
\begin{align*}\rho(a+\mathcal{K}, b+\mathcal{K}) := \rho(a,b).
\end{align*}
This is easily seen to be independent of the choices of $a$ and $b$, since, once we have a suitable pair $(a,b)$ with $a\in L_{p_1}$ and $b\in L_{p_1^{-1}}$, any other will be related by isometries induced by composites of $\phi$s and $\tilde{b}$s.  Clearly $\rho$ thus defined is a symmetric bilinear form and restricting to $f_*(L)_{q_1}\otimes f_*(L)_{q_2}$, it is essentially the original inner product restricted to $L_{p_1}\otimes L_{p_2}$, so is non-degenerate and satisfies
$$\rho(ab + \mathcal{K},c+\mathcal{K}) = \rho(a+\mathcal{K},bc+\mathcal{K}).$$

The next structure to check is the crossed $Q$-algebra action
$$\phi : Q \to Aut(f_*(L)).$$
The obvious formula to try is
\begin{align*}\phi_q(a+\mathcal{K}) := \phi_p(a) +\mathcal{K} 
\end{align*}
where  $f_0(p) = q$.  It is easy to reduce the proof that this is well defined to checking independence of the choice of $p$, but if $p^\prime$ is another element of $f_0^{-1}(q)$, then $p^\prime = np$ for some $n\in N$ and $\phi_{p^\prime}(a) = \phi_n\phi_p(a) \equiv_\mathcal{K} \phi_p(a)$, so it is well defined.  Of course, this definition will give us immediately that the $\phi_q(a)b = ba$ axiom holds and that $\phi_q|_{f_*(L)_q} = id$, etc.

The trace axiom follows from this definition by arguments similar to that used in the corresponding result for crossed $\pi$-algebras in \cite{turaev:hqft1} \S10.3; the requirement, there, that the kernel be central is avoided since $\phi_n(a)-a$ is defined to be in $\mathcal{K}$.

\begin{proposition}\hspace*{1mm}\\
With the above structure, $f_*(L)$ is a crossed $\mathcal{D}$-algebra.
\end{proposition}
\textbf{Proof}

The above argument shows it is a crossed $Q$-algebra, so we only have to define the tilde.  The obvious definition
is 
\begin{align*}\tilde{d} : = \tilde{c} + \mathcal{K},
\end{align*}
where $f_1(c) = d$.  This works.  It is well defined as each $\tilde{b} -1$ is in $\mathcal{K}$, and the equation
$$\phi_q(\tilde{d}) = \widetilde{{}^qd}$$
follows from the corresponding one in $L$.\hfill $\square$
\begin{proposition}
There is a natural bijection between the set of crossed algebra morphisms from $L$ to $L^\prime$ over $f$ and the set of crossed $\mathcal{D}$-algebra morphisms from $f_*(L)$ to $L^\prime$.\hfill $\square$
\end{proposition}
The proof is obvious given our construction of $f_*(L)$. We note that this, with its companion result on pulling back, give a pairs of adjoint functors determined by $f : \mathcal{C}\to \mathcal{D}$ between the categories of crossed $\mathcal{C}$-algebras and crossed $\mathcal{D}$-algebras.  We expect these to prove very useful when exploring in more depth the structure of crossed $\mathcal{C}$-algebras in future papers in this area, especially when looking at the relationship between such categories when $\mathcal{C}$ and $\mathcal{D}$ are weakly equivalent crossed modules which therefore model the same 2-type.


\noindent Timothy Porter,\\
School of Computer Science,\\
University of Wales Bangor,\\
Bangor, Gwynedd LL57 1UT, U.K.

\medskip

\noindent Vladimir Turaev,\\
Institut de Recherche Math\'{e}matique Avanc\'{e}e de Strasbourg,\\
7 rue Ren\'{e}-Descartes,\\ 67084 Strasbourg Cedex,\\ France.
\end{document}